\documentclass[fontsize=11pt,a4paper,headings=normal,
twoside=false,leqno,parskip=half-,abstract=true]{scrartcl}
\usepackage[T2A,T1]{fontenc}
\usepackage[utf8]{inputenc}
\usepackage{CJKutf8} 
\usepackage[russian,greek,english]{babel}
\usepackage{enumerate}
\usepackage{csquotes}

\usepackage[hyphens]{url}
\usepackage{hyperref}

\hypersetup{
 pdftitle={Kuramoto},
 pdfauthor={Bernold Fiedler}, 
 colorlinks=true,
 linkcolor=blue,
 citecolor=blue,
 filecolor=blue,
 urlcolor=blue}

 \usepackage{afterpage}
\usepackage{graphicx}
\usepackage{wrapfig} 
\usepackage[format=plain,labelfont=bf,font=small]{caption}
\usepackage{xcolor}
\usepackage[arrow, matrix, curve]{xy}
\usepackage{float}

\usepackage{caption}
\usepackage{tabulary}
\usepackage{array}
\newcolumntype{N}[1]{>{\centering\arraybackslash}m{#1}}

\usepackage{amsmath,amsfonts,amsthm,amssymb}
\swapnumbers 
\usepackage{eurosym} 
\newcommand{\transv}{\mathrel{\text{\tpitchfork}}}
\makeatletter
\newcommand{\tpitchfork}{%
  \vbox{
    \baselineskip\z@skip
    \lineskip-.52ex
    \lineskiplimit\maxdimen
    \m@th
    \ialign{##\crcr\hidewidth\smash{$-$}\hidewidth\crcr$\pitchfork$\crcr}
  }%
}
\makeatother
\usepackage{latexsym}
\usepackage{enumerate}

\usepackage[notref,notcite,color,final 
]{showkeys}

\definecolor{refkey}{rgb}{1,0,0}
\definecolor{labelkey}{rgb}{1,0,0}

\usepackage{cancel}

\usepackage{tikz}

  \mathchardef\ordinarycolon\mathcode`\:
  \mathcode`\:=\string"8000
  \begingroup \catcode`\:=\active
    \gdef:{\mathrel{\mathop\ordinarycolon}}
  \endgroup

\theoremstyle{plain}
\newtheorem{thm}{Theorem}[section]
\newtheorem{lem}[thm]{Lemma}
\newtheorem{prop}[thm]{Proposition}
\newtheorem{cor}[thm]{Corollary}

\newtheorem{rem}[thm]{Remark}

\newcommand\eps{\varepsilon}
\newcommand\mi{\mathrm{i}}
\newcommand\vth{\vartheta}
\newcommand\Th{\Theta}
\renewcommand\rho{\varrho}
\renewcommand\phi{\varphi}
\renewcommand\Re{\mathrm{Re}\,}

\newcommand{\bZ}{\mathbb{Z}}

\newcommand{\bR}{\mathbb{R}}

\newcommand{\bS}{\mathbb{S}}
\newcommand{\bT}{\mathbb{T}}

\newcommand{\bO}{\mathbb{O}}

\newcommand{\bE}{\mathbb{E}}

\newcommand{\bG}{\mathbb{G}}
\newcommand{\bJ}{\mathbb{J}}
\newcommand{\bK}{\mathbb{K}}

\def\cA{\mathcal{A}}
\def\cC{\mathcal{C}}

\def\bolth{{\boldsymbol{\vartheta}}}
\def\boTh{{\boldsymbol{\Theta}}}
\def\boxi{{\boldsymbol{\xi}}}
\def\boXi{{\boldsymbol{\Xi}}}
\def\bopsi{{\boldsymbol{\psi}}}
\def\boPhi{{\boldsymbol{\Phi}}}

\begin{document}

\title{\Large{Transient rebellions in the Kuramoto oscillator: \\[2mm]
Morse--Smale structural stability and \\[2mm]
connection graphs of finite 2-shift type\\[1cm]}}
\vspace{1cm}
\subtitle{
 \emph{Dedicated to Professor Yoshiki Kuramoto}\\
 \emph{in celebration of the 50th anniversary} \\
 \emph{of the classical Kuramoto oscillator} \\
 }

\author{
\\
Jia-Yuan Dai*, Bernold Fiedler**, Alejandro L\'{o}pez-Nieto***
\\[1cm]
}

\date{\small{version of \today}}
\maketitle
\thispagestyle{empty}
\vfill
*\\
Department of Mathematics \\
National Tsing Hua University \\
No. 101, Sec. 2, Kuang-Fu Road \\ 
Hsinchu, Taiwan; \\
National Center for Theoretical Sciences \\
National Taiwan University \\
No. 1, Sec. 4, Roosevelt Road \\ 
Taipei, Taiwan \\[5mm]
**\\
Institut für Mathematik\\
Freie Universität Berlin\\
Arnimallee 3\\ 
14195 Berlin, Germany\\[5mm]
***\\
Department of Mathematics\\
National Taiwan University\\
No. 1, Sec. 4, Roosevelt Road,\\
10617 Taipei, Taiwan \\


\newpage
\pagestyle{plain}
\pagenumbering{roman}
\setcounter{page}{1}

\begin{abstract}
\noindent
\noindent
The celebrated 1975 \emph{Kuramoto model} of $N$ oscillators with phase vector $\bolth=(\vth_1,\ldots,\vth_N)$, frequencies $\omega_j \in \bR$, and all-to-all coupling constant $K>0$ reads
\begin{equation}
\label{*}
\dot\vth_j\,=\,\omega_j+\tfrac{K}{N}\sum_{k=1}^N\,\sin(\vth_k-\vth_j)\,.
\end{equation}
The model is highly accessible to rigorous mathematical analysis, and has been studied as a paradigm for effects like total and partial synchronization.

\smallskip\noindent
We assume identical frequencies $\omega_j=\omega$, pass to co-rotating coordinates, and normalize the time scale in \eqref{*} to
\begin{equation}
\label{**}
\dot\vth_j\,=\,\tfrac{1}{N}\sum_{k=1}^N\,\sin(\vth_k-\vth_j)\,. 
\end{equation}
Most initial conditions $\bolth$ lead to total synchronization. 
The large number of asymptotically $\sim 2^{N\!-\!1}$ (circles of) partially synchronized states, however, is unstable.
Their precise global role in the transitions to synchrony seems to have eluded description.

\smallskip\noindent
In the present paper, we address this gap.
By the gradient structure of \eqref{**}, the global dynamics decompose into equilibria and heteroclinic orbits among them.
Except for the extremes of total instability and total synchrony, all equilibria are 2-cluster solutions: their phase angles $\vth_j$ attain only two values, with a phase difference of $\pm\pi$ between them.

\smallskip\noindent
Any heteroclinic orbit between 2-cluster equilibria, or towards synchrony, can be realized as a 3-cluster rebellion.
The rebellion splits the smaller, ``slim'', minority cluster of the source equilibrium and sends the rebellious part to join the bigger, ``fat'', majority cluster.
This identifies the heteroclinic \emph{connection graph} among 2-cluster equilibria, and towards synchrony, as the lattice of fat index sets, partially ordered by inclusion.

\smallskip\noindent
All heteroclinic orbits are in fact geometrically transverse.
This reveals the normalized Kuramoto model as a structurally stable Morse--Smale system.
In particular, heteroclinic orbits can be concatenated by shadowing orbits.
The options involved in successive rebellious splittings amount to finite symbol sequences of 2-shift type.
Moreover, the connection graph persists for any 
small deviations from the totally symmetric Kuramoto paradigm.
We include numerical illustrations and a brief comparison with related phenomena in Stuart-Landau oscillators under identical all-to-all coupling. 
\end{abstract}

\newpage 
\vspace{1cm}
\tableofcontents
\listoffigures



\newpage
\pagenumbering{arabic}
\setcounter{page}{1}
\setcounter{equation}{0}

\section{Introduction} \label{Int}

\numberwithin{equation}{section}
\numberwithin{figure}{section}
\numberwithin{table}{section}

\subsection{The classical Kuramoto model}\label{Kurclass}

More than 50 years ago, Kuramoto introduced his now classical model of $N$ oscillators in a famous 3-page proceedings note \cite{KurOsc}; see also \cite[section 5.4]{Kurbook}.
With phase angles $\bolth=(\vth_1,\ldots,\vth_N)$ on the $N$-torus $\bT^N\!=\bR^N/2\pi\bZ^N$, frequencies $\omega_j \!\in\! \bR$, and all-to-all coupling constant $K\!>\!0$, the \emph{classical Kuramoto model} reads
\begin{equation}
\label{omegaj}
\dot\vth_j\,=\,\omega_j+\tfrac{K}{N}\sum_{k=1}^N\,\sin(\vth_k-\vth_j)\,, 
\end{equation}
for $j=1,\ldots,N$.
The model has been designed towards accessibility by rigorous mathematical analysis.
Therefore, and for its broad applied relevance, it has widely been studied as a paradigm for effects like total and partial synchronization.
 Google Scholar attests to more than 2,500 citations for \cite{KurOsc} alone, and more than 10,000 for \cite{Kurbook}.
In the continuum limit $N\rightarrow\infty$, the celebrated results by Chiba \cite{Chib} have established \emph{phase transition}, away from unique stable distributions, under sufficiently large perturbations of homogeneous frequency distributions.
The partial synchronization phenomenon of \emph{chimeras} is widely studied in variants of \eqref{omegaj}, and has originally been discovered in a setting of one-dimensionally coupled phase oscillators \cite{KBatt}.
For broader surveys see also \cite{Str, ABPRS, PiRoSurv} and the many references there.
A neat web app illustration is \cite{Kuramotocycle}.

Our focus, in contrast, is on elucidating qualitative details of \emph{transients} towards synchronization, for $N\!<\!\infty$.
We mainly address the simplest case of identical frequencies $\omega_j\!=\!\omega$.
This allows us to pass to co-rotating coordinates $\vth_j(t)\!-\!\omega t$ and normalize the coupling constant $K$ in \eqref{omegaj} to
\begin{equation}
\label{kurODE}
\dot\vth_j\,=\,\tfrac{1}{N}\sum_{k=1}^N\,\sin(\vth_k-\vth_j)\,. 
\end{equation}
Depending on context, here and below, we use the same notation for angles $\bolth=(\vth_1,\ldots,\vth_N)$ in the $N$-torus $\bT^N=\bR^N/2\pi\bZ^N$ and in the covering space $\bR^N$.
Most initial conditions $\bolth$ lead to total synchronization; see for example \cite{DoXu, HKP}.
The plethora of partially synchronized states, however, which rule the detailed transition to synchrony, seems to have eluded detailed description.
Our main objective will be to address this gap.

\subsection{Euclidean group actions}\label{Group}
A first general feature of the Kuramoto model \eqref{kurODE}, which substantially facilitates a global description of the dynamics, is equivariance under a large symmetry group $\bG$.
We therefore recall some basic concepts concerning groups and group actions, in appropriate language.
For some of the enormous local and global consequences of equivariance on dynamics see for example \cite{Vander,FieHopf,GoStSymm,ChoLau, GoStPersp} and the references there.
In section \ref{kurGroup}, we return to the co-rotating Kuramoto model \eqref{kurODE} and its specific group $\bG$; see \eqref{groupS}, \eqref{groupT}.

On standard Euclidean space $\bolth\!\in\!\bR^N$, we consider arbitrary bounded $C^1$ vector fields
\begin{equation}
\label{ODEf}
\dot\bolth\,=\,\mathbf{f}(\bolth)\,,
\end{equation}
in the present section.
Let $\boPhi^t(\bolth_0):=\bolth(t),\ t\!\in\!\bR$, denote the \emph{global flow} of the associated initial value problem $\bolth(0)\!=\!\bolth_0$\,.
See for example \cite{ArnoldGeoODE,ArnoldODE,HSD}.
In other words, the ODE \eqref{ODEf} defines a nonlinear group action of time $t\!\in\!\bR$ on $\bR^N$, i.e.\ a group homomorphism $t\mapsto\boPhi^t$ from the group $(\bR,+)$ to the group of diffeomorphisms under composition.

Let $\bE(N)=\bO(N)\!\ltimes\!\bR^N$ denote the \emph{Euclidean group} of isometries of $\bR^N$, i.e.\ the semidirect product of the \emph{orthogonal group} $\sigma\!\in\!\bO(N)$ with the \emph{translations} $\bopsi\!\in\!\bR^N$.
Composition $g_1g_2$ of $g_\iota=(\sigma_\iota\,,\,\bopsi_\iota)\in\bE(N)$ for $\iota\!=\!1,2$ is defined as
\begin{equation}
\label{EN}
(\sigma_1\,,\,\bopsi_1)\circ(\sigma_2\,,\,\bopsi_2) := (\sigma_1\sigma_2\,,\,\bopsi_1+\sigma_1\bopsi_2)\,.
\end{equation}
The product is semidirect, because the projection $P_1$ onto the orthogonal component,
\begin{equation}
\label{P1}
\begin{aligned}
   P_1\!:\  \bE(N)&\rightarrow\bO(N)\,,\\
   (\sigma,\bopsi)&\mapsto\sigma\,,
\end{aligned}
\end{equation}
is a surjective homomorphism, with the normal subgroup of translations $\bopsi\!\in\!\bR^N$ as kernel.
The standard action of $g\!=\!(\sigma,\bopsi)\!\in\!\bE(N)$ on $\bolth\!\in\!\bR^N$ is by the \emph{affine linear} map
\begin{equation}
\label{g}
(g\bolth):=\sigma\bolth+\bopsi\,.
\end{equation}
Let $\bG\!\leq\!\bE(N)$ denote a subgroup.
The set $\bG\bolth:=\{g\bolth\,|\,g\!\in\!\bG\}$ is called the \emph{group orbit} of $\bolth$ under $\bG$.
We call a flow $\boPhi^t$ \emph{equivariant} under a subgroup $\bG\!\leq\!\bE(N)$, if $\boPhi^t$ commutes with $g\!\in\!\bG$, for all $t\!\in\!\bR$.
In particular, an equivariant flow $\boPhi^t$ maps $\bG$-orbits onto $\bG$-orbits.
Equivalently to equivariance of $\boPhi^t$, the ODE vector field $\mathbf{f}=\partial_t\vert_{t\!=\!0}\boPhi^t$ commutes with the orthogonal action of $\sigma\!\in\!\bO(N)$, and is invariant under translations by $\bopsi\!\in\!\bR^N$\,:
\begin{equation}
\label{fg}
\mathbf{f}(\sigma\bolth+\bopsi)= \sigma\mathbf{f}(\bolth)\,,
\end{equation}
for all $g\!=\!(\sigma,\bopsi) \!\in\! \bG$ and all $\bolth \!\in\! \bR^N$.

\emph{Relative equilibria} $\bolth(0)\!=\!\bolth_0$, for example, require the time orbit to be contained in the group orbit, i.e.
\begin{equation}
\label{releq}
\bolth(t)\in\bG\bolth_0\,,\quad\textrm{for all } t\in\bR\,.
\end{equation}
\emph{Equilibria} $\boTh$ satisfying $\mathbf{f}(\boTh)\!=\!0$ are a special case of relative equilibria.
By equivariance \eqref{fg}, the whole group orbit $\bG\boTh$ of any equilibrium $\boTh$ consists of equilibria, as well.

Let us identify the orthogonal group $\sigma\!\in\!\bO(N)$ with the subgroup $\bO(N)\!\times\!\{0\}\leq\bE(N)$, notationally.
For extensive studies of equivariant dynamics under such \emph{linear} group actions see the above references, again.
We just mention a few  powerful  elementary concepts. 
Let $\bJ\!\leq\!\bK\!\leq\!\bO(N)$ be subgroups and assume $\mathbf{f}$ is $\bK$-equivariant.
Define the $\bJ$-\emph{fixed vectors} as
\begin{equation}
\label{FixJ}
\mathrm{Fix}(\bJ):=\left\{\bolth\in\bR^N\;|\; \bJ\bolth\!=\!\{\bolth\}\right\}\,.
\end{equation}
Then the linear space $\mathrm{Fix}(\bJ)$ is invariant under the flow $\boPhi^t$ of \eqref{ODEf}.

Conversely, fix any $\bolth\!\in\!\bR^N$ and define the \emph{isotropy group} or \emph{symmetry} of $\bolth$ as
\begin{equation}
\label{Kth}
\bK_\bolth :=\{\sigma\!\in\!\bK\;|\; \sigma\bolth\!=\!\bolth\}\,.
\end{equation}
Choosing $\bJ\!=\!\bK_\bolth$ in \eqref{FixJ}, for example, we obtain a flow-invariant subspace $\mathrm{Fix}(\bK_\bolth)$ where the isotropy is at least $\bK_\bolth$\,. 
Along any particular solution, therefore, the isotropy
\begin{equation}
\label{Ktht}
\bK_{\bolth(t)} = \bK_{\bolth(0)} 
\end{equation}
remains invariant. 
Similarly, for convergent sequences $\bolth_n\,,\, n\!\rightarrow\!\infty$, of constant isotropy $\bK_*$\,, we observe lower semicontinuity, i.e.
\begin{equation}
\label{Kthn}
\bK_{\bolth_n} = \bK_*\,,\  \bolth_n\rightarrow\bolth\quad\Rightarrow\quad \bK_\bolth\geq\bK_*\,.
\end{equation}
Indeed, linear subspaces are closed.
Analogous statements hold for general subgroups $\bJ,\bK$ and isotropy subgroups $\bK_\bolth$ of $\bE(N)$, if we admit affine, rather than just linear, flow-invariant subspaces $\mathrm{Fix}(\bJ)$ of $\bR^N$.
For example, passage to $\boldsymbol{\alpha}$- or $\boldsymbol{\omega}$-limit sets never decreases isotropy.
For explicit examples, compare clustering at $\bolth(t_n)\rightarrow\boTh_\pm$\,, for $t_n\rightarrow\pm\infty$, in theorems \ref{thmmon}, \ref{thm3clu}, \ref{thmMclu}.

For $\bG$-equivariant $\mathbf{f}$, the \emph{linearization} $\mathbf{f}'(\bolth)$
at any $\bolth\!\in\!\bR^N$ satisfies
\begin{equation}
\label{f'}
\mathbf{f}'(\sigma\bolth+\bopsi)\,\sigma \,=\, \sigma\,\mathbf{f}'(\bolth)\,,
\end{equation}
for all $g\!=\!(\sigma,\bopsi) \!\in\! \bG$ and all $\bolth \!\in\! \bR^N$.
This follows by differentiation of \eqref{fg} with respect to $\bolth$.
In particular, $A:=\mathbf{f}'(\bolth)$ commutes with the orthogonal part $P_1\bG_\bolth$ of the isotropy $\bG_\bolth$\,.
Therefore, the $A$-invariant decomposition of $\bR^N$ into irreducible linear representations of $P_1\bG_\bolth$ assists the spectral decomposition of $A$, by Schur's Lemma.
See for example section \ref{Part}.

\subsection{Kuramoto equivariance: period lattices and clusters}\label{kurGroup}

We now apply these concepts to the specific Kuramoto system \eqref{kurODE}, written in the form \eqref{ODEf}.
For group aspects of general networks see also \cite{GoStNet}.
\emph{Equivariance} holds under the group $g\!=\!(\sigma,\bopsi)\!\in\!\bG=\mathbf{S}\!\ltimes\!\mathbf{T}\leq\bE(N)=\bO(N)\!\ltimes\!\bR^N$, with
\begin{align}
\label{groupS}
\sigma\in\mathbf{S}&:=\{\pm\mathbf{id}_N\}\times S_N\,,\\
\label{groupT}
\bopsi\in\mathbf{T}&:=\Lambda+\bR\mathbf{1}_N\,. 
\end{align}
Here the \emph{period lattice} is abbreviated as $\Lambda:=2\pi\bZ^N$, the identity matrix on $\bR^N$ is $\mathbf{id}_N$\,,  and the column vector $\mathbf{1}_N:=(1,\ldots,1)$ indicates uniform phase shift.
Equivariance under $-\mathbf{id}_N\!\in\!\mathbf{S}$ only emerges in co-rotating coordinates \eqref{kurODE}, after elimination of the uniform rotation frequency $\omega_j\!=\!\omega$.
The permutations $\sigma\!\in\!S_N$ act by
\begin{equation}
\label{sigma}
(\sigma\bolth)_{\sigma(j)}:=\vth_j\,.
\end{equation} 
As a first consequence, distinct $\vth_j(t)$ cannot cross in finite time.
This follows from time-invariance \eqref{FixJ} of $\mathrm{Fix}(\bJ)$ under the $\bZ_2$-group $\bJ\!\leq\! S_N$ generated by any transposition $j_1\leftrightarrow j_2$ of components.

We address affine translations $\bolth\mapsto\bolth+\bopsi$ by $\bopsi\!\in\!\mathbf{T}$ next.
The action of $\bopsi\!\in\!\mathbf{T}$ is \emph{free}: any $\bolth\!\in\!\bR^N$ is only fixed by $\bopsi\!=\!0$.
Since the equivariant flow $\boPhi^t$ maps $\bG$-orbits to $\bG$-orbits, we seek to eliminate the action of $\mathbf{T}$ by suitable coordinates.

To eliminate the uniform shifts $\bR\mathbf{1}_N$, we use the orthogonal projection
\begin{equation}
\label{P}
\begin{aligned}
\mathbf{P}\!: \bR^N  &\rightarrow \boXi:=\mathbf{1}_N^\bot\,,   \\
    \bolth&\mapsto \boxi:=\bolth-\tfrac{1}{N} (\mathbf{1}_N^\mathsf{T}\cdot\bolth)\, \mathbf{1}_N\,,
\end{aligned}
\end{equation}
onto the linear section $\boxi\!\in\!\boXi$ of zero averages, $\xi_1+\ldots+\xi_N=0$.
In short, $\mathbf{P}=\mathbf{id}_N-\tfrac{1}{N}\mathbf{1}_N\mathbf{1}_N^\mathsf{T}$\,.
The orthogonal projection $\mathbf{P}$ is equivariant under $\mathbf{S}$ and eliminates $\bR\mathbf{1}_N$\,.
By \eqref{kurODE}, the $\bR\mathbf{1}_N$-section $\boXi\!=\!\mathbf{1}_N^\bot$ is flow-invariant:
\begin{equation}
\label{Xiinv}
\tfrac{d}{dt}(\mathbf{1}_N^\mathsf{T}\cdot\bolth)=\mathbf{1}_N^\mathsf{T}\cdot\dot\bolth=0\,.
\end{equation}
In fact, we have an integral of motion: $\mathbf{1}_N^\mathsf{T}\cdot\bolth=\mathrm{const}$\,.
The restriction of the flow $\boPhi^t$ to the flow-invariant section $\boXi$ defines the flow on the space of group orbits $\boXi=\bR^N/\bR\mathbf{1}_N$\,.

Equivariance of \eqref{kurODE}, \eqref{ODEf} on $\bR^N$, under the period lattice $\Lambda=2\pi\bZ^N\!\leq\!\mathbf{T}$ of \eqref{groupT}, defines a dynamics on the flat quotient $N$-torus
\begin{equation}
\label{TN}
\bolth\in\bT^N:=\bR^N/\Lambda\,,
\end{equation}
for which we abuse identical notation.
Here $\vth_j\!\in\!\bR$ recover their original meaning as angles $\vth_j\!\in\!\bS^1\!=\!\bR/2\pi\bZ$.

We collect a few facts on the relation between $\Lambda$ and the flow-invariant $\bR\mathbf{1}_N$-section $\boXi\!=\!\mathbf{1}_N^\bot$\,.
\begin{prop}[Lattices and cross section]\label{propLambda}
When we quotient out the full translation group $\mathbf{T}=\Lambda+\bR\mathbf{1}_N$ of \eqref{groupT}, we arrive at an $(N\!-\!1)$-torus covered by  the $\bR\mathbf{1}_N$-section $\boXi\!=\!\mathbf{1}_N^\bot$\,:
\begin{equation}\label{TN-1}
\begin{aligned}
   \bT^{N-1}&:=\bR^N/\mathbf{T}=\boXi/\Lambda'\,,\quad\textrm{with}   \\
   \Lambda'&:=\mathbf{P}\Lambda\leq\boXi\,.  
\end{aligned}
\end{equation}
The lattice $\Lambda\cap\boXi$ is a subgroup of index $N$ in the projected lattice $\Lambda'\!=\!\mathbf{P}\Lambda$\,. More precisely,
\begin{equation}
\label{LambdaN}
\Lambda'/(\Lambda\cap\boXi) \ \cong\  2\pi(\bZ/N\bZ)\,.
\end{equation}
\end{prop}

\emph{\textbf{Proof.}}\quad
To prove claim \eqref{TN-1}, let $\boxi,\tilde{\boxi}\!\in\!\boXi$ be in the same $\mathbf{T}$-orbit.
Then $\bopsi:=\tilde{\boxi}\!-\!\boxi\in\mathbf{T}\cap\boXi$, i.e.\ $\bopsi\!=\!\mathbf{P}\bopsi\!\in\!\mathbf{P}\Lambda\!=\!\Lambda'$.
Conversely, let $\boxi,\tilde{\boxi}\!\in\!\boXi$ be such that $\bopsi:=\tilde{\boxi}\!-\!\boxi\in\Lambda'\!=\!\mathbf{P}\Lambda$.
Then $\ker\mathbf{P}=\bR\mathbf{1}_N$ implies $\psi\!\in\!\mathbf{T}$, i.e.\ $\boxi$ and $\tilde{\boxi}$ lie in the same $\mathbf{T}$-orbit.
The projected period lattice $\Lambda'\!\leq\!\boXi$ possesses full rank, because $\Lambda\!\leq\!\bR^N$ does. 
This proves \eqref{TN-1}.

We prove claim \eqref{LambdaN} for $\Lambda\!=\!\bZ^N$, i.e.\ omitting all factors $2\pi$ without loss of generality.
Let $\mathbf{P}\vert_\Lambda$ denote the restriction of the projection $\mathbf{P}$ to $\Lambda$, and define $\mathbf{1}_N^*\lambda:= \mathbf{1}_N^\mathsf{T}\cdot\lambda\ \mathrm{mod}\,N$, for all $\lambda\!\in\!\Lambda$.
Then we can define a homomorphism $H$ by the diagram
\begin{equation}
\label{figlambda}
\vcenter{\hbox{\includegraphics[width=0.2\textwidth]{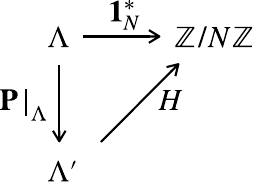}}}
\end{equation}
Indeed, the kernel of $\mathbf{P}\vert_\Lambda$ is $\ker \mathbf{1}_N^*\cap \bR\mathbf{1}_N$\,.
In particular the homomorphism $H$ is well-defined.
Inserting integer multiples $\lambda:=k\,\mathbf{e}_1$ of the first unit vector $\mathbf{e}_1\!\in\!\Lambda$ shows surjectivity of $H$.
Since $\ker H=\Lambda\cap\boXi$, this proves claim \eqref{LambdaN}, and the proposition.
\qed

We address the orthogonal permutation action \eqref{sigma} of $\sigma\!\in\!S_N\!\leq\!\mathbf{S}\!\leq\!\bO(N)$ next. 
For $S_N$-equivariance and coupled oscillator settings including Stuart--Landau scenarios, more specifically, see also \cite{Elm1,DiSt,StElm,Elm2,Krischer1,Krischer2,GoStNet}.
The linear orthogonal representation \eqref{sigma} decomposes $\bR^N$ into the $S_N$-invariant subspaces $\bR\mathbf{1}_N$, with the trivial irreducible representation, and the irreducible orthogonal complement $\boXi\!=\!\mathbf{1}_N^\bot$\,.
The irreducible representation on $\boXi$ is called the \emph{standard representation} of $S_N$\,.

Time-invariance of \emph{cluster subspaces} $\mathbf{X_J}\!\leq\!\bR^N$ is an immediate consequence of $S_N$-equivariance.
To define these subspaces, let $\mathbf{J}\!=\!(J_m)$ denote any partition
\begin{equation}
\label{J}
\{1,\ldots,N\}=J_1\,\dot\cup\,\ldots\,\dot\cup\, J_M
\end{equation}
into nonempty disjoint index sets $J_m$\,, for $m=1,\ldots,M$.
Let $\sigma\!\in\! J_m\!\leq\! S_N$ denote the permutations $\sigma$ which only permute elements of $J_m$\,. 
Then $\sigma\in S_\mathbf{J}:=S_{J_1}\times\ldots\times S_{J_M}\leq S_N$ denotes the subgroup of those permutations $\sigma$ which permute elements of each $J_m$\,, separately.
In other words, $S_\mathbf{J}$ respects the partition $\mathbf{J}$ in \eqref{J}.
With \eqref{FixJ} for $\bJ:=S_\mathbf{J}$, we define the cluster subspace  
as the space
\begin{equation}
\label{XJ}
\mathbf{X_J}:=\mathrm{Fix}(S_\mathbf{J})
\end{equation}
of those $\bolth\!\in\!\bR^N$ which remain fixed under the group $S_\mathbf{J}$\,.
I.e., the components $\vth_j$ are identical, separately on each index set $J_m$\,:
\begin{equation}
\label{xm}
\vth_j=x_m\,, \quad\textrm{for all } j\in J_m\,.
\end{equation}
We call the reduced coordinate vector $\mathbf{x}=(x_1,\ldots,x_M)\in \mathbf{X_J}\cong\bR^M$ an $M$-\emph{cluster}.
By construction \eqref{FixJ}, \eqref{XJ}, the cluster spaces are flow-invariant, of isotropy \eqref{Kth} at least $S_\mathbf{J}$\,.

Let $N_m$ denote the sizes of the clusters $J_m$ and abbreviate the size fractions as $\alpha_m:=N_m/N \!>\! 0$.
Note $\sum \alpha_m=1$.
The $M$-cluster dynamics of the Kuramoto model \eqref{kurODE} on the time-invariant $M$-cluster subspace $\mathbf{X_J}$ then reads
\begin{equation}
\label{cluODE}
\dot x_m\,=\,\sum_{m'=1}^M\,\alpha_{m'}\, \sin(x_{m'}-x_{m})\,.
\end{equation}
Although all $\alpha_m$ are rational, with shared denominator $N$, we may just as well discuss \eqref{cluODE} for real $\alpha_m\!>\!0$.
The original case of $M\!=\!N$ single-oscillator $N$-``clusters'' $J_m\!=\!\{m\}$, of course, is subsumed in the above notation.
Then $x_m\!=\!\vth_m$ with $\alpha_m\!=\!1/N,\ S_\mathbf{J} \!=\! \{\mathbf{id}_N\},\ \mathbf{X_J}\!=\! \mathbf{X}:=\bR^N$, and no ``clustering'' at all.

Since the orthogonal projection $\mathbf{P}\!:\bR^N\rightarrow\boXi$ is $S_N$-equivariant, the above discussion extends to $\boxi\!\in\!\boXi$, verbatim, defining the flow-invariant subspaces
\begin{equation}
\label{XiJ}
\boXi_\mathbf{J}:=\mathrm{Fix}(S_\mathbf{J})\cap\boXi\,.
\end{equation}
Working with $\boxi\!\in\!\boXi$, instead of $\bolth\!\in\!\bR^N$, we now have $\xi_j\!=\!x_m$ in \eqref{xm}, with the constraint
\begin{equation}
\label{S1}
0\,=\,\tfrac{1}{N}\,\sum_{j=1}^N \xi_j \,=\, \sum_{m=1}^M \alpha_m\, x_m \,.
\end{equation}
The period lattices $\Lambda$ and $\Lambda'\!=\!\mathbf{P}\Lambda$ are $S_N$-invariant.
Therefore, $S_N$ acts on the tori $\bT^N\!=\!\bR^N/\Lambda$ and $\bT^{N\!-\!1}/\Lambda'$, as well.
In particular, restriction defines covering projections to the restricted tori
\begin{align}
\label{TJM}
   \bT_\mathbf{J}^M&:=\mathbf{X_J}/\Lambda_\mathbf{J}\,,\quad\textrm{with } \Lambda_\mathbf{J}:= \Lambda\cap\mathbf{X_J}\,; \\
\label{TJM-1}
    \bT_\mathbf{J}^{M-1}&:=\boXi_\mathbf{J}/\Lambda_\mathbf{J}'\,,\quad\textrm{with } \Lambda_\mathbf{J}':= \Lambda'\cap\boXi_\mathbf{J}=\mathbf{P_J}\Lambda_\mathbf{J}\,.
\end{align}
We caution our reader that the restricted projection $\mathbf{P_J}:=\mathbf{P}\vert_\mathbf{X_J}$ does not quite look ``orthogonal'' in coordinates $\mathbf{x}$ on $\mathbf{X_J}$; see \eqref{PJ}.
In fact, the proper Euclidean length $|\mathbf{x}|_\mathbf{J}$ on $\mathbf{X_J}$ satisfies
\begin{equation}
\label{||J}
\tfrac{1}{N}|\mathbf{x}|_\mathbf{J}^2:=\tfrac{1}{N}|\bolth|^2=\alpha_1x_1^2+\ldots+\alpha_Mx_M^2\,.
\end{equation}

Caution is also required, because clustering \eqref{xm} takes on a different meaning, when projecting from the covering section $\bolth\!\in\!\boXi\!=\!\mathbf{1}_N^\bot$\,, in proposition \ref{propLambda}, down to the torus $\bT^{N\!-\!1}\!=\!\boXi/\Lambda'$\,.
Clusters in $\boXi$ are defined by \eqref{xm}, as is.
Already in $\bT^N\!=\!\bR^N/2\pi\bZ^N$, however, we have to define clustering \eqref{xm} $\mathrm{mod}\ 2\pi$.
We use the same definition, i.e.\ componentwise $\mathrm{mod}\ 2\pi$, in $\bT^{N-1}=\boXi/\Lambda'$.
Indeed, $\vth_j=\vth_k\ \mathrm{mod}\ 2\pi$ holds, for some $\bolth\!\in\!\boXi$, if and only if $\vth_j'=\vth_k'$ holds for the same components of some representative $\bolth'\in\bolth+\Lambda'\in\boXi/\Lambda'$. 
This componentwise observation works for tori $\bT_\mathbf{J}^{M\!-\!1}\!=\!\boXi_\mathbf{J}/\Lambda'_\mathbf{J}$\,, as well.
Therefore we mostly work with definition \eqref{xm}, directly, for such representatives $\bolth'$.
Occasionally, on the other hand, we will have to emphasize lifts which break clustering -- e.g.\ for continuity reasons.
The targets $\boTh_+$ versus $\tilde{\boTh}_+$ of right versus left rebellions, in section \ref{Main} below, will be notable examples; see in particular the phase resets by $2\pi$ in \eqref{x2pm}, in remark \ref{resetTheta}, and in figure \ref{fig52}.

\subsection{Example: Two-cluster dynamics}\label{2Clu}

As a most elementary example, consider 2\emph{-clusters} $\mathbf{J}=J_1\dot\cup J_2$\,, i.e.\ $M\!=\!2$.
Later, we will consider 2-cluster dynamics embedded into the dynamics of $M\!\geq\!3$ clusters.
To avoid confusion, in the cluster notation \eqref{xm}, we will therefore distinguish the 2-cluster case by superscripts $J^b:=J_1$ and $J^c:=J_2$ of sizes $N^b$ and $N^c$, respectively.
The index set $J^c=J_1^c:=\{1,\ldots,N\}\setminus J^b$ is just complementary to $J^b\!=\!J_1$\,.
Size fractions $\alpha_1,\alpha_2$ and cluster coordinates $x_1,x_2$ of $J_1, J_2$ will be denoted as  $\alpha^b,\alpha^c$ and $x^b,x^c$, respectively.
Without loss of generality, we label the size fractions such that $0\!<\!\alpha^c\!\leq\!\alpha^b\!<\!1$.
In other words, the superscript ``$b$'' refers to the ``\underline{\textbf{b}}ulky'', fat cluster $J^b$, and ``c'' refers to the \underline{\textbf{c}}omplementary, slim cluster $J^c$.
Then \eqref{cluODE} reads
\begin{equation}
\label{2cluODE}
\begin{aligned}
\dot x^b \,&= \alpha^c\,\sin(x^c-x^b)\,,   \\
\dot x^c &= \,\alpha^b\,\,\sin(x^b-x^c)\,,
\end{aligned}
\end{equation}
in coordinates $\mathbf{x}\!=\!(x^b,x^c)$ on $\mathbf{X_J}$\,, with period lattice $\Lambda_\mathbf{J}\!=\!2\pi\bZ^2$. 
See also \eqref{TJM}.
Normalization \eqref{S1} focuses attention to the time-invariant line 
\begin{equation}
\label{2line}
\boXi_\mathbf{J}\!: \quad \alpha^b x^b+\alpha^cx^c=0
\end{equation}
of rational slope $-\alpha^c/\alpha^b$ with $\alpha^b\!+\!\alpha^c\!=\!1$\,.
In coordinate $x^c$, the resulting dynamics on $\boXi_\mathbf{J}$ reads
\begin{equation}
\label{2cluODEJ}
\dot x^c=-\alpha^b\,\sin(x^c/\alpha^b)\,.
\end{equation}
For $\dot x^b$\,, we just swap all superscripts ``$b$'' and ``$c$''.

We claim that the projected period lattice of \eqref{TJM-1} takes the form $\Lambda_\mathbf{J}' \!=\!\mathbf{P_J}\Lambda\!\cong\! 2\pi\alpha^b\,\bZ$, in co\-ordinate $x^c$ on the line $\boXi_\mathbf{J}$\,.
To see this, we write the restriction $\mathbf{P_J}$ to $\mathbf{X_J}$ of the projection $\mathbf{P}$ from \eqref{P} in coordinates $\mathbf{x}\!=\!(x^b, x^c)$ as
\begin{equation}
\label{PJ}
\mathbf{P_J\,x}\ =\ 
\mathbf{x}-\frac{1}{N}
\begin{pmatrix}
    N^b  &  N^c  \\
    N^b  & N^c 
\end{pmatrix}\mathbf{x}\ =\ 
(x^b-x^c)\begin{pmatrix}
   \alpha^c   \\
   -\alpha^b  
\end{pmatrix}\,.
\end{equation}
Here we have used $\alpha^b\!=\!N^b/N,\ \alpha^c\!=\!N^c/N$, and again $\alpha^b\!+\!\alpha^c\!=\!1$.
Inserting $\mathbf{x}\!\in\!\Lambda_\mathbf{J}\!=\!2\pi\bZ^2$ in \eqref{TJM-1} proves our claim.
Similarly, $\Lambda_\mathbf{J}'\!\cong\!2\pi\alpha^c\,\bZ$ in coordinate $x^b$\,.

In view of \eqref{TJM-1}, the ODE \eqref{2cluODEJ} also defines the dynamics on the 1-torus $x^c\in\bS^1=\bT_\mathbf{J}^1=\boXi_\mathbf{J}/\Lambda_\mathbf{J}'=\bR/2\pi\alpha^b\bZ$\,.
The angle-equivalent
\begin{equation}
\label{endpoints}
\emph{\textrm{endpoints:}}\quad x^b\!=\!x^c\!=\!0 \quad\textrm{and}\quad x^b\!=\!-2\pi\alpha^c\,,\ x^c\!=\!2\pi\alpha^b\,, 
\end{equation}
of the fundamental interval on the covering line $\boXi_\mathbf{J}$\,, mark an exponentially attracting, synchronous 1-cluster equilibrium $\bolth\!\equiv\!\boTh_0\!=\!0\ \mathrm{mod}\,2\pi$.
Indeed, $x^c\!-\!x^b\!=\!2\pi$ at the second endpoint, because $\alpha^b\!+\!\alpha^c\!=\!1$.
The 
\begin{equation}
\label{midpoint}
\emph{\textrm{midpoint:}} \quad x^b\!=\!-\pi\alpha^c\,,\ x^c\!=\!\pi\alpha^b\,,
\end{equation}
of the fundamental interval on $\boXi_\mathbf{J}$\,, with $x^c\!=\!x^b\!+\!\pi$, identifies the \emph{unique} 2\emph{-cluster equilibrium} $\bolth\!\equiv\!\boTh_-\!\neq\!0$ in $\bS^1\!=\!\bT_\mathbf{J}^1$\,, for any given 2-cluster partition $\mathbf{J}$.
Each half of the fundamental interval consists of a heteroclinic orbit from the exponentially repelling 2-cluster equilibrium $\boTh_-$ to total synchrony $\boTh_0$\,.
Trivial as this one-dimensional 2-cluster dynamics may appear, it will become an important ingredient of our global 3-cluster analysis; see proposition \ref{Equil}, section \ref{3Clu}.
The lower sides $x_3\!=\!x_1$ of the fundamental triangle in 3-cluster figures \ref{fig51} and \ref{fig52}, for example, illustrate the resulting 1-torus dynamics, in coordinates $x_1\!=\!x_3\!=\!x^b$ (vertical) and $x_2\!=\!x^c$ (horizontal), for cluster ratios $\alpha^b\!=\!\alpha_1\!+\!\alpha_3\!=\!4/5$ and $\alpha^c\!=\!\alpha_2\!=\!1/5$.

\subsection{The order parameter as Lyapunov function}\label{Grad}

A second general feature of the Kuramoto model \eqref{kurODE} is the \emph{gradient structure}; compare for example \cite{vanHemmen}.
This is manifest in the \emph{order parameter} $R\!\in\![0,1]$ and $\Psi\!\in\!\mathbb{S}^1$ which are defined as polar coordinates of the averages, alias mean fields
\begin{equation}
\label{RPsi}
R\exp(\mi\Psi)\,:=\,\tfrac{1}{N}\,\sum_{k=1}^N \exp(\mi\vth_k)\,=\,\sum_{m=1}^M \alpha_m\exp(\mi x_m)\,.
\end{equation}
Here we have already subsumed the $M$-cluster case \eqref{xm}.
The order parameter $R$ is invariant under the symmetry group $\bG=\mathbf{S}\!\ltimes\!\mathbf{T}$ of \eqref{groupS}, \eqref{groupT}.
The angle parameter $\Psi$ is invariant under the period lattice $\Lambda$ and the permutation group $S_N$, but co-rotates with angle translations $\bR\mathbf{1}_N$\,.
Multiplication by $\exp(-\mi\Psi)$ and comparison of real and imaginary parts yields
\begin{align}
\label{R=}
R\,&=\,\,\tfrac{1}{N}\,\sum_{k=1}^N \cos(\vth_k-\Psi)\,=\,\sum_{m=1}^M \alpha_m\cos(x_m-\Psi)\,;\\
\label{0=}
0\,&=\,\tfrac{1}{N}\,\sum_{k=1}^N \sin(\vth_k-\Psi)\,=\,\sum_{m=1}^M \alpha_m\sin(x_m-\Psi)\,.
\end{align}
At $R\!=\!0$, of course, these hold for any angle $\Psi$.
For all $R\!\geq\!0$ and all $j$, multiplication of \eqref{RPsi} by $\exp(-\mi\vth_j)$ and comparison of the resulting imaginary parts with \eqref{kurODE} reveals 
\begin{equation}
\label{thjdot}
\dot\vth_j=R\sin(\Psi-\vth_j)\,.
\end{equation}
For clusters \eqref{xm}, of course, we may as well substitute $\vth_j$ by $x_m$\,.
To determine $\dot R$ for $R\!>\!0$, we differentiate \eqref{R=} along solutions.
Eliminating $\dot\Psi$ by \eqref{0=} provides the \emph{gradient structure}
\begin{equation}
\label{Lyap}
\dot R\,=\,R\cdot\tfrac{1}{N} \,\sum_{k=1}^N\sin^2(\vth_k-\Psi)\,=\,R\cdot\,\,\sum_{m=1}^M\alpha_m\sin^2(x_m-\Psi)\,\geq\, 0\,,
\end{equation}
with the order parameter $0\!<\!R\!\leq\!1$ as a \emph{nondecreasing Lyapunov function}.

\subsection{Equilibria and (normal) hyperbolicity}\label{Equi}

Equilibria $\bolth\!=\!\boTh$ are defined by $\dot\bolth\!=\!\mathbf{f}(\boTh)\!=\!0$.
To eliminate redundancy by the translation group $\mathbf{T}$ of \eqref{groupT}, we describe equilibria in the sectional torus $\bolth\!\in\!\bT^{N\!-\!1}\!=\!\boXi/\Lambda'$ of \eqref{TN-1}.

To determine all equilibria, we first claim that $\dot R$ vanishes at equilibria, only.
Indeed, the set $\{\boTh \!\in\! \boXi\,|\, R \!=\! 0\}$ consists of equilibria, only, by \eqref{thjdot}.
For $R\!>\!0$, $\dot R \!=\! 0$ in \eqref{Lyap} implies $\sin^2(\Th_j\!-\!\Psi)=0$, for all $j$.
Again, \eqref{thjdot} implies we are at equilibrium.
Moreover, the vanishing $\sin^2$-terms imply that any component $\Th_j\!\in\!\{\Psi,\Psi\!+\!\pi\}$ can attain at most two values $\mathrm{mod}\ 2\pi$.
Therefore, $\boTh$ is either a synchronous 1-cluster equilibrium, which defines the set where $R\!=\!1$ is maximal, or a 2-cluster equilibrium, at $0<R=2\alpha^b\!-\!1<1$.
In section \ref{2Clu} we have already discussed 2-cluster equilibria $\boTh\!=\!\boTh_-$\, on the line $\boXi_\mathbf{J}$\,.
In particular, $\alpha\!=\!\alpha^b\!=\!N^b/N\!>\!1/2$ indicates the \emph{fat cluster} of size $N^b\!>\!N/2$ at cluster coordinate $x^b=(\alpha^b\!-\!1)\pi$.
The \emph{slim cluster} of complementary size $N^c\!=\!N\!-\!N^b\!<\!N/2$ is located at $x^c\!=\!\alpha^b\pi$, by the phase constraint \eqref{S1}.
The symmetric case $\alpha^b\!=\!\alpha^c\!=\!1/2$ of $R\!=\!0$ is included in the following remarks. 

The equilibrium set of all  $\boTh$ in covering space $\boXi$ or in $\bT^{N\!-\!1}\!=\!\boXi/\Lambda'$ such that $R\!=\!0$, i.e.
\begin{equation}
\label{R=0}
\sum_{j=1}^{N}\exp(\mi\Th_j)=0\,,
\end{equation}
defines \emph{closed planar} $N$\emph{-bar linkages}, a ubiquitously important device in mechanical engineering. 
In our case, the length of each joint in the complex plane is 1, and the directional angles are $\Th_j$\,.
They define an algebraic variety (a manifold for odd $N$) of codimension two in $\boXi\!\cong\!\bR^{N\!-\!1}$.
The linkage $N$-gons, defined by the partial sums of \eqref{R=0}, may be non-convex, and may even intersect themselves.
By permutation and phase-shift to increasing order $0\!=\!\Th_0\!\leq\!\Th_1\!\leq\!\ldots\!\leq\!\Th_N\!<\!2\pi$,
however, they can be represented by (possibly degenerate) convex polygons.
General $M$-bar linkages with rational joint lengths occur as $M$-cluster variants.
See the monograph \cite{Nbar} for a wealth of fascinating details, from an engineering point of view.

A neighborhood of sufficiently small $R\!\geq\!0$ can elegantly be studied via a Möbius transform of coordinates, which are called \emph{Watanabe--Strogatz variables}; see for example \cite{PiRoBunch} and the survey \cite{PiRoSurv}.
Unfortunately, these variables do not cover the whole phase space.
For example, the required normalization $\sum \exp(\mi\xi_j)=0$, in their notation, misses any 2-cluster equilibrium with order parameter $R\!>\!0$.
Indeed, equal cluster-angles $\vth_j\!=\!x_m$ lead to equal $\xi_j$\,, and clustering with $\alpha_1\!>\!1/2\!>\!\alpha_2$ contradicts the required Watanabe--Strogatz normalization of $\xi_j$\,.

An equilibrium $\boTh\!\in\!\bT^{N\!-\!1}$ is called \emph{hyperbolic}, if its linearization $\mathbf{f}'(\boTh)$ neither possesses purely imaginary nor zero eigenvalues, when restricted to the invariant linear section $\boXi$\,:
\begin{equation}
\label{hyp}
0\not\in\Re\,\textrm{spec}\ \mathbf{f}'(\boTh)\vert_\boXi\,.
\end{equation}
The \emph{unstable dimension}, or \emph{Morse index}, $i\!=\!i(\boTh)$ counts the eigenvalues with strictly positive real part, by their algebraic multiplicities. 
The case $i\!=\!0$ implies local asymptotic stability.

Our preliminary observations on equilibria will lead to the following trichotomy.

\begin{prop}[Classification of equilibria]\label{Equil}
In the sectional torus $\bolth\!\in\!\bT^{N\!-\!1}\!=\!\boXi/\Lambda'$ of \eqref{TN-1}, the equilibria $\bolth\!=\!\boTh$ of the Kuramoto model \eqref{kurODE} fall into three classes.
The three classes can be characterized by decreasing order parameter $R$ and decreasing stability as follows.
\begin{enumerate}[(i)]
  \item $R=1\mathrm{:}$\\ 
  the unique stable synchronous equilibrium $\boTh\!=\!0$ of unstable dimension $i\!=\!0$;
  \item $0<R=2\alpha^b\!-\!1<1\mathrm{:}$\\ the hyperbolic 2-cluster equilibria $\boTh$ of unstable dimension $i\!=\!|J^c|$\,, with rational fat and slim size fractions 
  \begin{equation}
\label{2clualpha}
\begin{aligned}
    |J^b|/N&= \,\alpha^b\, \in(\tfrac{1}{2},1)\,,     \\
   1-\alpha^b\!=\! |J^c|/N&=\,\alpha^c\,\in(0,\tfrac{1}{2})\,,
\end{aligned}
\end{equation}
   respectively.
   In $\bT^{N\!-\!1}$, each such equilibrium is determined uniquely by its fat, ``bulky'' cluster $J^b$ or, equivalently, by its complementary slim cluster $J^c$.
   See section \ref{2Clu}, where these equilibria are represented by their cluster coordinate $x^c\!=\!\pi\alpha^b$ at the midpoint \eqref{midpoint} of the fundamental domain of $x^c\in\bR/2\pi\alpha^b\,\bZ$;
  \item $R=0\mathrm{:}$\\ 
  the totally unstable set of closed planar $N$-bar linkages $\sum_{j=1}^N \exp(\mi\Th_j)=0$ of unit joint lengths, rotated such that $\sum_{j=1}^N \Th_j=0$ is normalized; see \eqref{R=0}.
\end{enumerate}
\end{prop}

We do not claim much novelty here.
The fundamental distinction between these types of equilibria already underlies \cite[section 5.4]{Kurbook}.
The trichotomy has already been developed in \cite[Theorem 1]{IEEE}, except for hyperbolicity and the precise Morse indices $i\!=\!i(\boTh)$.
Instead, the clusters (ii) were shown to be ``of saddle type'' with respect to a Lyapunov potential proportional to $R^2$.

The (in-)stability claims (i) and (iii) follow from the Lyapunov property \eqref{Lyap} and the absence of further equilibria near $R\!=\!1$ and $R\!=\!0$, respectively.
The set $R\!=\!0$ of $N$-bar linkage equilibria is a total \emph{repeller}. 
Indeed, the order parameter $R$ is a strictly increasing Lyapunov function, locally, upon any small perturbation of $R\!=\!0$ towards $R\!>\!0$.
Similarly, total synchrony $R\!=\!1$ is a local attractor.
We have already discussed the 2-cluster equilibria of claim (ii); see section \ref{2Clu}.
This establishes all claims of the proposition, except for the precise Morse indices $i\!=\!i(\boTh)\!=\!|J^c|$ which we determine in section \ref{Part}.

In the full torus of angle vectors $\boTh\in\bT^N=\bR^N/2\pi\bZ^N$ of \eqref{TN}, each of the above isolated equilibrium representatives gives rise to a unique circle of equilibria, by the line $\boTh+\bR\mathbf{1}_N$ in covering space $\bR^N$.
Hyperbolicity carries over to each of the flow-invariant sections $\boXi+c$\,; see \eqref{f'}, \eqref{Xiinv}.
The associated circles of equilibria are therefore \emph{normally hyperbolic} and the Morse indices $i\!=\!i(\boTh)$ are their \emph{strong unstable dimensions} in $\bR^N$.
At any given critical level $0<R=2\alpha^b\!-\!1<1$, we therefore encounter $\binom{N}{\alpha^b N} = \binom{N}{|J^b|}$ distinct circles of 2-cluster equilibria.
The cases correspond to the partitions of  $\{1,\ldots,N\}$ into pairs of nonempty complementary index sets $J^b$ and $J^c$ of sizes $|J^b|=\alpha^b N\!>\!N/2$ and $|J^c|=(1\!-\!\alpha^b)N<N/2$\,, respectively.

The LaSalle invariance principle for the Lyapunov function $R$ implies $\dot R\!=\!0$ on $\boldsymbol{\alpha}$- and $\boldsymbol{\omega}$-limit sets of any bounded solution \cite{Hale}.
We have seen how $\dot R\!=\!0$ holds at equilibria, only, in the co-rotating setting of \eqref{kurODE}.
The Lyapunov function $R$ therefore implies that $\boldsymbol{\alpha}$- and $\boldsymbol{\omega}$-limit sets consist of equilibria, only.
This prevents time-periodic orbits, and any more complicated recurrent dynamics.
In the original Kuramoto model \eqref{omegaj} with identical nonzero rotation frequencies $\omega_j\!=\!\omega\!\neq\!0$, however, the equilibrium circles become circles of relative equilibria \eqref{releq} under the circle action of $\mathbf{1}_N$ on $\bT^N$.
In fact, they become hyperbolic time-periodic orbits.
This is not a contradiction, because $R$ is invariant under the circle action.
In particular, \eqref{Lyap} remains valid, and $R$ remains a Lyapunov function.

\subsection{Invariant manifolds and heteroclinic time orbits}\label{Het}

We return to co-rotating Kuramoto dynamics \eqref{kurODE}.
By the Lyapunov function $0\!<\!R\!\leq\!1$ and hyperbolicity, non-equilibrium solutions $\bolth_*(t)$ in the compact torus $\bT^{N\!-\!1}\!=\!\boXi/\Lambda'$ of \eqref{TN-1}, or its covering space $\boXi$, become heteroclinic between distinct equilibria $\boTh_\pm$\,, for $t\rightarrow\pm\infty$.
In symbols,
\begin{equation}
\label{leadsto}
\bolth_*(t): \quad\boTh_-\leadsto \boTh_+ \neq \boTh_- \,.
\end{equation}
The only exception are the repelling $N$-bar linkage source equilibria $\boTh_-$\,, for which we only asserted local convergence to the equilibrium set \eqref{R=0} of unit $N$-bar linkages in reverse time $t\rightarrow-\infty$.
Henceforth, we only address heteroclinicity among 1- and 2-cluster equilibria $\boTh_\pm$ with $0<R_\pm:=R(\boTh_\pm)\!\leq\!1$, including convergence to the stable 1-cluster equilibrium $R_+\!=\!1$ of total synchrony.
Passing to $\boTh_\pm\!\in\!\bT^{N\!-\!1}\!=\!\boXi/\Lambda'$ of \eqref{TN-1} picks one of $N$ hyperbolic representatives for each equilibrium circle; see proposition \ref{propLambda}.
The \emph{connection graph} will consist of these equilibria, as vertices, and directed edges which represent the sets of rebellious heteroclinic orbits \eqref{leadsto} between such equilibria.
See corollary \ref{corconngraph}.
We occasionally conflate the $N$-bar linkages $R\!=\!0$, as a single top repeller.

Let $W_+^s$ and $W_-^u$ denote the stable and unstable manifolds, respectively, of the two hyperbolic equilibria $\boTh_\pm$\,, in the flow-invariant covering space $\bolth\!\in\!\boXi\!\cong\!\bR^{N\!-\!1}$ of angle normalization \eqref{S1}.
As usual, we subsume the synchronous 1-cluster equilibrium $\boTh_+\!=\!0$ at $R\!=\!1$ here.
For some general background on invariant manifolds, see for example \cite{HPS, ChowHale, PalisdeMelo}. 
\emph{Any} heteroclinicity $\bolth_*(t){:}\, \ \boTh_-\!\leadsto\! \boTh_+$ is then equivalent to a nonempty intersection
\begin{equation}
\label{cap}
\bolth_*(t)\in W_-^u\cap W_+^s\,.
\end{equation}
We say that $W_-^u$ and $W_+^s$ are \emph{transverse}, or \emph{intersect transversely}, along the heteroclinic orbit $\bolth_*(t)\!\in\!\boXi$, if the sum of their tangent spaces spans all of $\boXi$, i.e.
\begin{equation}
\label{trv}
T_{\,\bolth_*(t)}W_-^u\, +\, T_{\,\bolth_*(t)}W_+^s\ =\ \boXi
\end{equation}
for some, and hence all, $t\!\in\!\bR$.
In symbols, we also express such transversality as
\begin{equation}
\label{trvcap}
W_-^u\ \transv_{\,\bolth_*(t)} W_+^s\,,\quad\textrm{alias}\quad \boTh_-\mathrel{\mathop{\leadsto}_{\transv}^{\bolth_*(t)}}\boTh_+\,.
\end{equation}

Let $g=(\sigma,\bopsi)\in\bG=\mathbf{S}\!\ltimes\!\mathbf{T}$ be an equivariance of \eqref{kurODE}, as in section \ref{kurGroup}.
Then equivariance properties \eqref{g}, \eqref{f'} imply
\begin{align}
\label{Wg}
    W^{s,u}(g\boTh)\ &=\ \quad\ \sigma\, W^{s,u}(\boTh)+\bopsi \,,  \\
\label{TWg}
    T_{g\bolth}\,W^{s,u}(g\boTh)\ &=\ \sigma\,T_\bolth\, W^{s,u}(\boTh)\,,
\end{align}
for any equilibrium $\boTh$ and any $\bolth\!\in\! W^{s,u}(\boTh)$.
In particular, $W^{s,u}(\boTh)$ and their transversality in $\boXi$ may be studied in $\boXi/\Lambda'$, without loss of generality.

In theorem \ref{thmtrv} we will establish \emph{automatic} transversality of stable and unstable manifolds along arbitrary heteroclinic orbits $\bolth_*(t){:}\, \ \boTh_-\!\leadsto\! \boTh_+$ with $0\!<\!R_-\!<\!R_+\!\leq\!1$.
As an immediate consequence of transversality at $\bolth_*(t)$, the time-invariant intersection of the stable and unstable manifolds is itself a manifold, of dimension
\begin{equation}
\label{ipm}
\dim (W_-^u\cap W_+^s) \ =\ i_--i_+=|J_-^c|-|J_+^c|\,.   
\end{equation}
Here $i_\pm\!=\!i(\boTh_\pm)\!=\!|J_\pm^c|$ denote the unstable dimensions or Morse indices of the hyperbolic equilibria $\boTh_\pm$ in $\boXi$, respectively.
See proposition \ref{Equil}.
Factoring out the time direction, e.g.~by a Poincaré cross section transverse to the flow on $\boXi$, this provides a local manifold of heteroclinic orbits of dimension $i_-\!-\!i_+\!-\!1\geq 0$.
Transverse heteroclinic orbits between adjacent Morse indices $i_\pm$\,, in particular, are isolated and hence locally unique.
Transverse homoclinic orbits, or saddle-saddle connections with $i_-\!\leq\! i_+$, are excluded, \emph{per se}.

\subsection{Outline}\label{Out}

In the next section, we collect our \emph{main results} on heteroclinic orbits \eqref{leadsto}, in the time-invariant section space $\boXi\!=\!\mathbf{1}_N^\bot$ of \eqref{P}, between equilibria $\boTh_\pm$ with order parameters $0\!<\!R_-\!<\!R_+\!\leq\!1$.
Specifically, theorem \ref{thmmon} establishes necessary relations between source equilibrium $\boTh_-$ and target equilibrium $\boTh_+$\,, for any heteroclinic orbit.
Central claim (iii) formulates strict monotonicity of fat clusters.
The more elementary aspects (i) and (ii) of theorem \ref{thmmon} are proved right there.
Theorem \ref{thm3clu} then shows how any heteroclinic edge in the connection graph can already be realized by a rebellious 3-cluster heteroclinic orbit, which only splits the slim cluster of the 2-cluster source equilibrium $\boTh_-$\,.
Towards $\boTh_+$\,, i.e.\ for $t\rightarrow +\infty$, the rebellious cluster joins the fat majority $J_+^b$\,, whereas the remaining slim third cluster recovers antipodal opposition at angle difference $\pi$.
Corollary \ref{corxshift} distinguishes left and right rebellions of the rebellious cluster, when it switches allegiance.
Theorem \ref{thmtrv} establishes transversality, in full space, of general heteroclinic orbits which need not be 3-clusters.
In corollaries \ref{strucstab}--\ref{shift-type}, we collect some consequences of these results.
In particular, we mention structural stability for $0\!<\!R_-\!<\!R_+\!\leq\!1$, and transitivity of the connection graph.
Theorem \ref{thmMclu} characterizes general heteroclinic rebellions $\bolth_*(t){:}\, \ \boTh_-\!\leadsto\! \boTh_+$ and their realization by multi-cluster swarms.

Section \ref{Part} introduces \emph{partitioned coordinates} for a deeper analysis of the interplay between 3-cluster heteroclinicity and general perturbations of initial conditions.
The proof of proposition \ref{Equil} on 2-cluster equilibria is an immediate illustration.
In section \ref{Wsu}, clustered coordinates provide a global description of stable and unstable manifolds.
In particular, lemma \ref{lemWpm}(ii) proves strict monotonicity of fat clusters, as claimed in theorem \ref{thmmon}(iii).
Theorem \ref{thm3clu} on 3-cluster realization of heteroclinic orbits is proved in section \ref{3Clu}, and transversality theorem \ref{thmtrv} follows in section \ref{Trans}.
For comments on the numerical illustrations of our results see section \ref{Numerics}.
We conclude with an updated discussion, as outlined at the beginning of section \ref{Disc}.

\subsection{Acknowledgment}\label{Ack}

We are grateful to Chien-Yu Chen for pointing out the connection problem between partially synchronous equilibria to Jia-Yuan Dai in 2023.
Arkady Pikovsky and Michael Rosenblum kindly shared their insight on Watanabe--Strogatz variables.
The referees have provided very diligent and insightful comments, thankfully,  during a 300-day refereeing period.
In particular, they motivated us to write expanded sections \ref{Group} and \ref{kurGroup} of the introduction and to clarify the role of the period lattice $\Lambda'$ in proposition \ref{propLambda} and \eqref{TJM}, \eqref{TJM-1}.
Moreover, they provided reference \cite{IEEE} for part of proposition \ref{Equil}.
They also pointed out the closely related later preprint by Burns et al., which we had come across, meanwhile, and which we discuss in our added section \ref{REV}.
Jia-Yuan Dai has been supported by NSTC (National Science and Technology Council) grant 113-2628-M-007-005-MY4. 
Alejandro López-Nieto was supported by NSTC grant 113-2123-M-002-009 of Chun-Hsiung Hsia.
Jia-Yuan Dai and Bernold Fiedler are indebted for mutual hospitality during in-person working visits, and for the pertinent grant support, in particular by the Deutsche Forschungsgemeinschaft.

\section{Main results}\label{Main}

In this section we formulate our main results on the global dynamics of the co-rotating Kuramoto model \eqref{kurODE}.
Throughout, and without loss of generality, we work in the flat torus $\bT^{N\!-\!1}\!=\!\boXi/\Lambda'$ and its flow-invariant covering space $\boXi\!=\!\mathbf{1}_N^\bot$\,; see proposition \ref{propLambda}.
This eliminates shift-equivariance through the time-invariant angle normalization \eqref{S1}.

We consider heteroclinic orbits $\bolth(t){:}\, \ \boTh_-\!\leadsto\! \boTh_+$ among synchronous 1-cluster and 2-cluster equilibria $\boTh_\pm$ at levels $0\!<\!R_- \!<\! R_+ \!\leq\! 1$ of the order parameter.
See proposition \ref{Equil} and \eqref{leadsto}.
We recall the notation $J_\pm^{b,c}$ to denote fat and slim clusters at $\boTh_\pm$\,, respectively, and similarly for their sizes $N_\pm^{b,c}\!=\!|J_\pm^{b,c}|$\,, size fractions $\alpha_\pm^{b,c}\!=\!N_\pm^{b,c}/N\,,\ \alpha_\pm^b\!+\!\alpha_\pm^c\!=\!1$, and phase angles $x_\pm^{b,c}$.
Theorem \ref{thmmon} collects necessary conditions for such heteroclinic rebellions.
Based on these conditions, theorem \ref{thm3clu} shows that, conversely, general heteroclinic connections can be represented by 3-cluster variants. 
Theorem \ref{thmtrv} establishes transversality of general heteroclinic orbits. 
Corollaries \ref{corconcat}--\ref{shift-type} collect some consequences for concatenation of heteroclinic rebellions and the connection graph.
Theorem \ref{thmMclu} addresses general swarms of heteroclinic $M$-clusters.

\begin{thm}[Angles and strict monotonicity of equilibrium  $2$-clusters]\label{thmmon} 
Consider any nonstationary time orbit $\bolth(t)\!\in\!\boXi,\ t\!\in\!\bR$, with a positive lower bound $0<R_-:= \inf_{t \in \bR} R(t)$ of the increasing order parameter $R(t)\!=\!R(\bolth(t))$.
Then $\bolth(t){:}\, \ \boTh_-\!\leadsto\! \boTh_+$ is heteroclinic between two equilibria $\boTh_\pm$ with order parameters $0\!<\!R_-\!<\!R_+\!\leq\! 1$.
Moreover, when viewed in $\bT^{N\!-\!1}\!=\!\boXi/\Lambda'$,
\begin{enumerate}[(i)]
  \item $\boTh_\pm$ are either 1- or 2-cluster equilibria, with size fractions $1/2 
 <\alpha_\pm^b=(R_\pm\!+\!1)/2\leq 1$ of their fat clusters $J_\pm^b$\,, respectively;
  \item the phase angles of their fat clusters $J_\pm^b$ are $x_\pm^b=-\alpha_\pm^c\pi$ in $\bS_\pm^1=\bR/2\pi\alpha_\pm^c\bZ$\,; 
  \item their fat clusters are supported on strictly growing index sets $J_-^b\subsetneq J_+^b$\,.
\end{enumerate}
We emphasize that (i)--(iii) refer to clusters at the source and target equilibria $\boTh_\pm$\,, only, and not to the heteroclinic orbit $\bolth(t)$.
The exceptional case $\boTh_+\!=\!0$ of a totally synchronous target 1-cluster is subsumed here as $R_+\!=\!\alpha_+^b\!=\!1$.
\end{thm}

\emph{\textbf{Proof of claims (i) and (ii).}}\quad
Convergence to equilibria in $\bT^{N\!-\!1}\!=\!\boXi/\Lambda'$ follows from relative compactness of the orbit $\bolth(t)\!\in\!\bT^{N\!-\!1}$, the Lyapunov property \eqref{Lyap}, and the LaSalle invariance principle \cite{Hale}. 
See section \ref{Equi}.
For $0\!<\!R_-\!<\!R_+\!\leq\!1$, convergence to single equilibria in $\bT^{N\!-\!1}$ follows because these equilibria are isolated; moreover $\boldsymbol{\alpha}$- and $\boldsymbol{\omega}$-limit sets in compact $\bT^{N\!-\!1}$ are connected.
See proposition \ref{Equil}(i)--(ii), and the explicit representation of equilibria in section \ref{2Clu}.
Convergence to single equilibria in covering space $\boXi$ then follows because the lattice $\Lambda'$ is discrete.
The remaining parts of claims (i)--(ii) of the theorem have already been established in proposition \ref{Equil}(i),(ii).
\qed

Claim (iii) on strict monotonicity of fat clusters will be established in lemma \ref{lemWpm}(ii).
This claim is somewhat surprising: once a cluster has gained absolute majority $\alpha_-^b\!>\!1/2$ of its supporters $J_-^b$\,, it will never lose any of its adherents by further heteroclinic rebellion.
Rather, rebellions only split the complementary slim minority cluster $J_-^c$\,.

\begin{figure}[t]
\centering \includegraphics[width=0.9\textwidth]{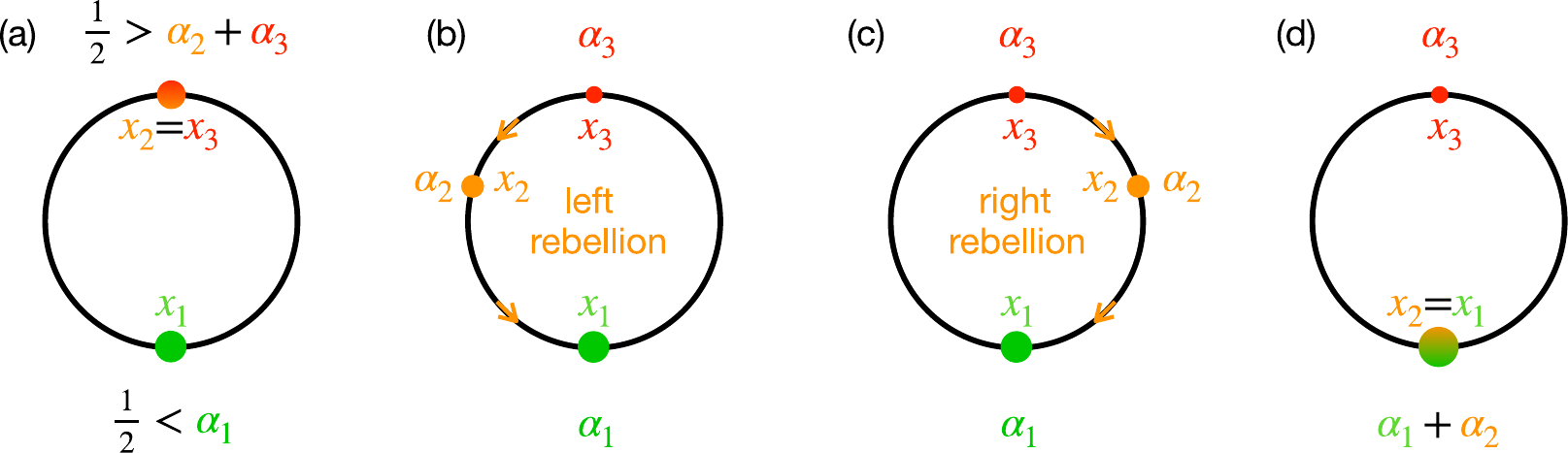}
\caption[Left and right rebellions]{\emph{
Heteroclinic left and right 3-cluster rebellions $\bolth_*(t){:}\, \ \boTh_-\!\leadsto\! \boTh_+$\,; see theorem \ref{thm3clu}.
Antipodal angles $x_1\!=\!x_\pm^b$ (green) and $x_3\!=\!x_\pm^b\!+\!\pi$ (red) of the 2-cluster source equilibrium $\boTh_-$ and the target equilibrium $\boTh_+$ are marked by solid dots on the phase circle $\mathbb{S}^1$, in (a) and (d). 
In (b) and (c), rebels $x_2(t)$ (orange) switch allegiance from (a), $x_2\!=\!x_3$ at $\bolth_*(-\infty)\!=\!\boTh_-$\,, towards (d), $x_2\!=\!x_1$ at $\bolth_*(+\infty)\!=\!\boTh_+$\,, in two different ways. 
Angles $x_2$ parametrize the positively (left, counter-clockwise) oriented circle.
Areas of dots refer to cluster size.\\
In (b) and (c): heteroclinic rebellions of cluster $x_2(t)$ (orange), which splits off from the slim $\boTh_-$-cluster at $x_2\!=\!x_3\!=\!x_-^c$ to join the fat $\boTh_+$-cluster at $x_2\!=\!x_1\!=\!x_+^b$.
See \eqref{x2pm}.
Left rebellion of  cluster $x_2$ is illustrated in (b), and right rebellion in (c); see cases \eqref{left} and \eqref{right}, respectively.
Cases (b) and (c) are swapped under $\sigma\!=\!-\mathbf{id}_N\!\in\!\mathbf{S}$; see \eqref{groupS}.
}}
\label{fig21}
\end{figure}

We now assert that, conversely, all heteroclinicity with $\inf_{t \in \bR} R(\bolth(t))\!>\!0$ is realized, either by a 3-cluster heteroclinic orbit between 2-cluster equilibria, or else by a 2-cluster heteroclinic orbit towards synchrony $R_+\!=\!1$.
This identifies 3-cluster heteroclinics as an elementary building brick for all heteroclinicity.

\begin{thm}[Realization by 3-cluster heteroclinic orbits]\label{thm3clu}
Let  $\boTh_\pm\!\in\!\bT^{N\!-\!1}\!=\!\boXi/\Lambda'$ denote any pair of equilibria with order parameters $0\!<\!R_-\!<\!R_+\!<\!1$.
Assume that the index sets $J_\pm^b$ of their fat clusters satisfy strict monotonicity $J_-^b\subsetneq J_+^b$.
Then there exists a 3-cluster heteroclinic orbit $\bolth_*(t){:}\, \ \boTh_-\!\leadsto\! \boTh_+$ in $\bT^{N\!-\!1}$.

More precisely,  the above pairs $\boTh_\pm\!\in\!\boXi/\Lambda'$ of 2-cluster equilibria are in one-to-one correspondence to certain index partitions $\mathbf{J}$ by $J_1\dot\cup J_2\dot\cup J_3$\,, as in \eqref{J}, with size fractions $0<\alpha_2,\alpha_3<1/2<\alpha_1<1$ of cluster sizes $N_1\!+\!N_2\!+\!N_3\!=\!N$.
Specifically, 
\begin{equation}
\label{J123}
\begin{aligned}
J_1&=J_-^b\,,\quad& J_2&=J_+^b\setminus J_-^b\,,&\quad  J_3&=J_+^c\,,\\
\alpha_1&=\alpha_-^b\,,& \alpha_2&=\alpha_+^b-\alpha_-^b\,, &\alpha_3&=1-\alpha_+^b\,; 
\end{aligned}
\end{equation}
see theorem \ref{thmmon}.
Let $\mathbf{x}(t)\!=\!(x_1,x_2,x_3)$ denote the 3-cluster angles \eqref{xm} of this partition.
Then $\mathbf{x}(t)\in \bT_\mathbf{J}^2=\boXi_\mathbf{J}/\Lambda_\mathbf{J}'$\,, under the normalization $\alpha_1x_1\!+\!\alpha_2x_2\!+\!\alpha_3x_3=0$ in covering space $\boXi_\mathbf{J}$.
See \eqref{S1}, \eqref{TJM-1}.
Since the cluster angles $x_m(t)$ for $m \!=\!1,2,3$ cannot cross, for finite $t$, their order relation must be preserved.
Consider any fundamental region of $\bT_\mathbf{J}^2$ in the strip $x_1\leq x_2, x_3\leq x_1\!+\!2\pi$.

\begin{figure}[t]
\centering \includegraphics[width=0.8\textwidth]{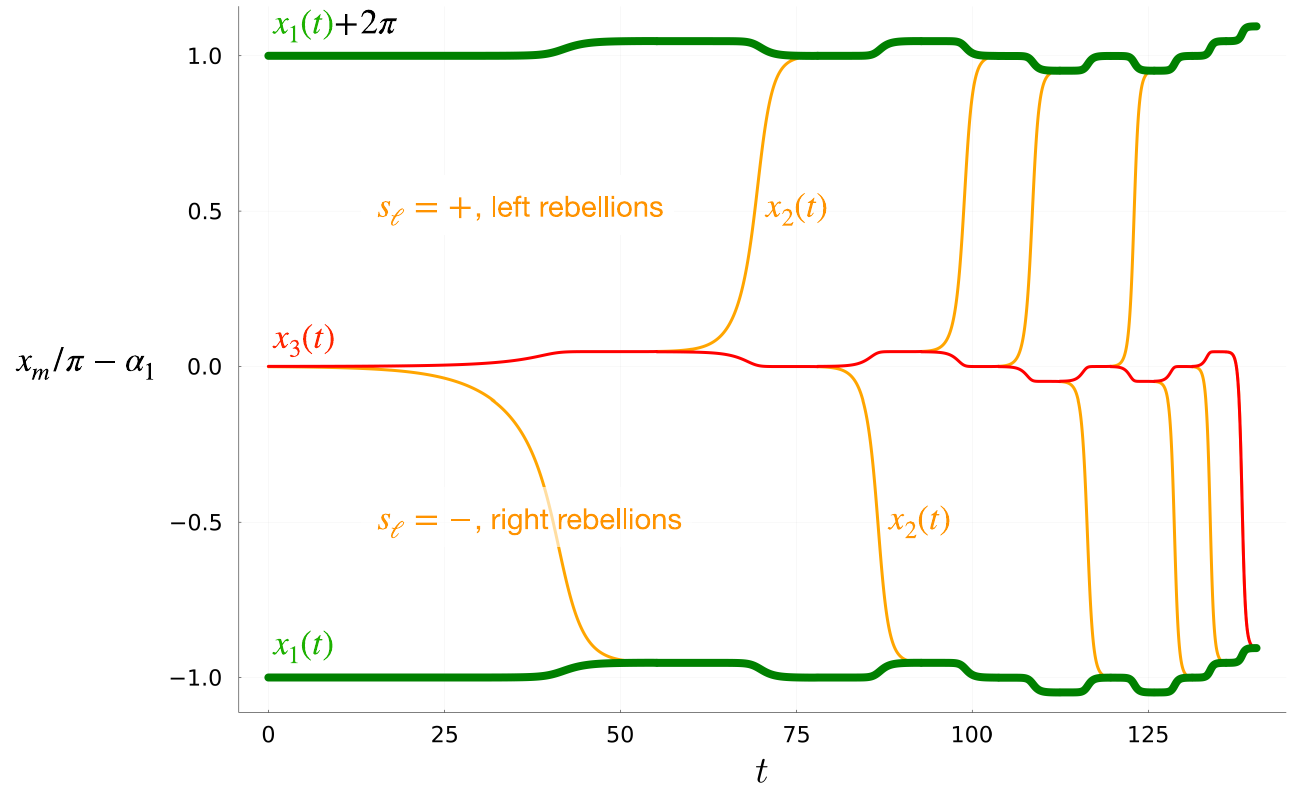} 
\caption[Concatenation of rebellions I]{\emph{
Time plot of $n\!=\!10$ successive transient one-man rebellions $\bolth_{\ell} (t)\!\in\!\boXi$ towards synchrony, concatenated to a single heteroclinic orbit $\bolth_*(t){:}\, \ \boTh_0\!\leadsto\!\boTh_n$\,. 
See corollaries \ref{corxshift}, \ref{corconcat}, \ref{shift-type}, and section \ref{NumConcat}.
The three cluster phases (vertical) of 3-cluster rebellions $\bolth_{\ell} (t)$ are $x_1(t)$ (fat, green), $x_2(t)$ (one-man rebellions $N_2\!=\!1$, orange) and $x_3(t)$ (slim, red).
Compare figure \ref{fig21} for left (here: upward) versus right rebellions (here: downward).
Initial cluster sizes are $N_1\!=\!11,\ N_2\!=\!1,\ N_3\!=\!9$ for $N\!=\!21$ oscillators.
Vertical coordinates $x_m$ have been scaled to $x_m/\pi\!-\!\alpha_1$\,, with $\alpha_1$ fixed at the initial fat size fraction $\alpha_1\!=\!11/21$.
Initial phases $\mathbf{x}(t)$ are near equilibrium, i.e.\ $x_1/\pi\!-\!\alpha_1=\pm1$ and $x_2/\pi\!-\!\alpha_1=x_3/\pi\!-\!\alpha_1=0$, for $t\rightarrow -\infty$.
Note the nine successive one-man rebellions (orange) of faster and faster ejection from the slimming red cluster, and convergence to the fattening green cluster.
Left rebellions $s\!=\!+$, as in figure \ref{fig21}(b), appear above the red graph of slim $x_3(t)$, and right rebellions $s\!=\!-$, as in (c), appear below.
In symbols \eqref{sequence}, $\mathbf{s}= -+-++-+---$.
The final tenth rebellion (red, down) is of 2-cluster type.
It vacates the vestigial near-antipodal red slim ``cluster'' of ``last man standing'' opposition at $J_2\!=\!J_3\!=\!\{21\}$, entirely, and submits to total synchrony, at last.
See \emph{\cite{FieTatra, Ash07, KuGiOtt2015}} and figure \ref{fig8} below for related ideas.
Note how each one-man rebellion $s_\ell\!=\!\pm$ induces an opposite drift by $-s_\ell\,\alpha_2\pi\!=\!\mp\pi/N$, in the position $x_1$ of the fat majority cluster. See corollary \ref{corxshift} and remark \ref{resetTheta}.}
}
\label{fig22}
\end{figure}

Then the 3-cluster heteroclinic orbit $\bolth_*(t)\!\in\!\bT_\mathbf{J}^2$ is supported on the index partition \eqref{J123} associated to $\boTh_\pm$\,.
But let us consider $\bolth_*(t)$ in \emph{covering space} $\boXi$.
Then $\lim x_1(t)\!=\!x_\pm^b$ and $\lim x_3(t)\!=\!x_\pm^c\!=\!x_\pm^b\!+\!\pi$, for $t \rightarrow \pm \infty$ respectively.
The phase  $x_2(t)$ of the \emph{rebel cluster} $J_2$\,, however, \emph{switches allegiance}, as follows.
Starting from $\lim x_2(t)\!=\!x_-^c\!=\!x_-^b\!+\!\pi$ at the source $\boTh_-$\,, for $t\rightarrow -\infty$, we encounter the alternative
\begin{equation}
\label{x2pm}
 \lim_{t\rightarrow +\infty}\, x_2(t)=
\begin{cases}
         \,  x_+^b\,, &\text{ or else} \\
  \,x_+^b+2\pi\,,&
\end{cases} 
\end{equation}
at the target $\boTh_+$\,.
Correspondingly, the 3-cluster heteroclinic orbit $\bolth_*(t)\!\in\!\boXi$ satisfies
\begin{align}
\label{right}
    \textrm{either}\quad x_1<x_2& <x_3<x_1+2\pi\quad\textrm{\emph{(right rebellion)}}\,, \\
\label{left}
    \textrm{or} \quad x_1<x_3& <x_2<x_1+2\pi\quad\textrm{\emph{(left rebellion)}}\,, 
\end{align}
for all $t\!\in\!\bR$. 
These options are realized by different heteroclinic orbits in $\bT^{N\!-\!1}\!=\!\boXi/\Lambda'$.
See figure \ref{fig21}(b) and (c), respectively, as well as the two fundamental triangles in figure \ref{fig52}.

In the case $R_+\!=\!\alpha_+^b\!=\!1$ of convergence towards total synchrony $\boTh_+\!=\!0,\ J_+^c\!=\!\emptyset$, the heteroclinic orbit $\bolth_*(t)$ can be chosen to be of 2-cluster type; see section \ref{2Clu}.
\end{thm}

In $\bT^{N\!-\!1}\!=\!\boXi/\Lambda'$, theorems \ref{thmmon} and \ref{thm3clu} state that any heteroclinic relation $\bolth{:}\,\ \boTh_-\!\leadsto\!\boTh_+$ can be realized by a \emph{heteroclinic orbit} $\bolth_*$ \emph{of maximal isotropy} $\bK_*\!=\!\bK_{\bolth_*(t)}$ in $\bK\!=\!S_N$\,.
Indeed, isotropy invariance \eqref{Ktht}, lower semicontinuity \eqref{Kthn} of isotropy, for $t\!\rightarrow\!\pm\infty$, and partitioning \eqref{J123} imply
\begin{equation}
\label{SNmax}
\bK_{\bolth(t)}\ \leq\ \bK_{\boTh_-}\!\cap\, \bK_{\boTh_+}\ =\ \bK_{\bolth_*(t)}\,.
\end{equation}

\begin{figure}[t]
\centering \includegraphics[width=0.5\textwidth]{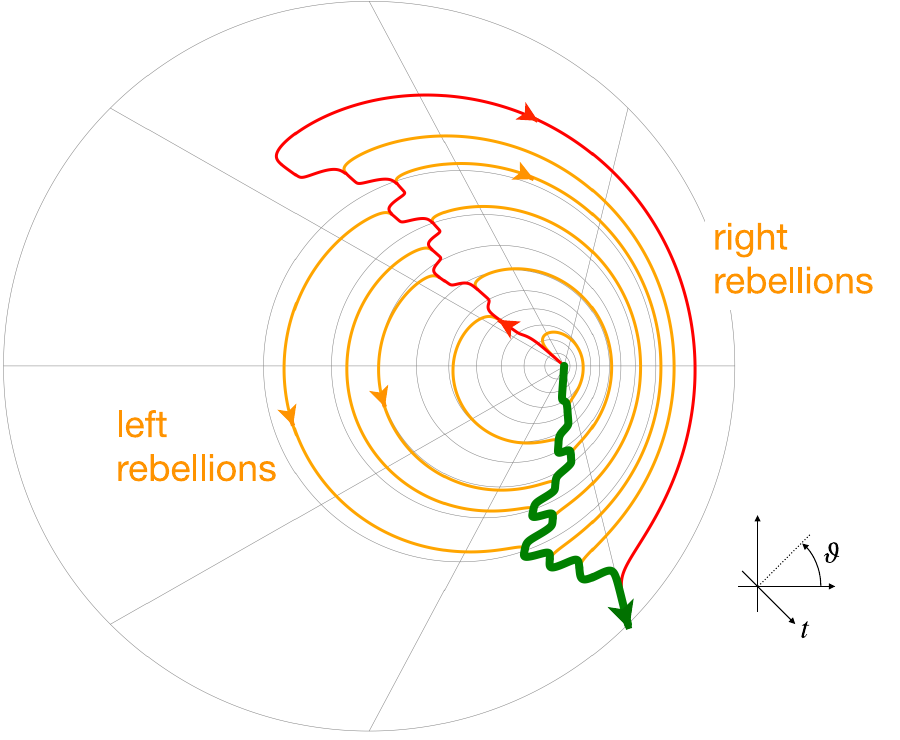} 
\caption[Concatenation of rebellions II]{\emph{
Phase plots of $n\!=\!10$ successive transient one-man rebellions $\bolth_{\ell} (t)\!\in\!\boXi$  towards synchrony, concatenated to a heteroclinic orbit $\bolth_*(t){:}\, \ \boTh_0\!\leadsto\!\boTh_n$\,.
The plot is based on the same data as figure \ref{fig22}, with the same color coding.
Phases $x_m\!\in\!\mathbb{S}^1$ are drawn on a cylinder $(x_m,t)\!\in\!\mathbb{S}^1\times\bR$ with vanishing point of perspective at $t \!=\! -\infty$.
Thin gray lines mark a cylindrically Cartesian grid with angle sectors of width $\pi/4$.
Clockwise orange right rebellions, right, follow figure \ref{fig21}(c).
Counter-clockwise orange left rebellions, left, follow figure \ref{fig21}(b).
Their tangencies near equilibria $\boTh_\ell$ point to the front, in the direction of time.
The outermost red clockwise arc indicates the final tenth rebellion (right) which vacates all opposition and submits to total synchrony, at last.}}
\label{fig23}
\end{figure}

For a proof of theorem \ref{thm3clu} and of the next corollary \ref{corxshift}, see section \ref{3Clu}.
In covering space $\boXi$, figure \ref{fig22} illustrates successive transient one-man rebellions $\alpha_2\!=\!1/N$ and their influence on the cluster positions $x_1, x_3$\,.
Here $x_1\!=\!x_1(t)$ keeps track of the angles $\vth_j$ of $j\!\in\! J_-^b$\,, and $x_3$ refers to $j\!\in\! J_3\!=\!J_+^c$\,; see \eqref{J123}.
To track the rebel position $x_2(t)$ in $\boXi$, we define the shift symbols
\begin{equation}
\label{spm}
s =
\begin{cases}
     -\,, & \textrm{for a right rebellion \eqref{right}}, \\
     +\,, & \textrm{for a left rebellion \eqref{left}},
\end{cases}
\end{equation}
as in figure \ref{fig21}(c), (b), respectively.

\begin{cor}[Phase shifts]\label{corxshift}
Consider any 3-cluster heteroclinic orbit $\bolth_*(t){:}\, \ \boTh_-\!\leadsto\! \boTh_+$ in $\boXi$, as in theorem \ref{thm3clu}.
Then the heteroclinic rebellion of $x_2(t)$ causes a shift of $x_1(t)$ for $t\rightarrow\pm\infty$ by 
\begin{equation} 
\label{xshift}
x_1(+\infty) - x_1(-\infty)=x_+^b-x_-^b=-s\,\alpha_2\pi\,,
\end{equation}
with the sign $s$ from \eqref{spm}, and similarly for the limits $x_3(\pm\infty)=x_1(\pm \infty)\!+\!\pi$.
\end{cor}

\begin{rem}[Phase reset at target equilibria]\label{resetTheta}
The alert reader may already have noticed:  the condition \eqref{xm}, of identical variables $\vth_j\!=\!x_m$ on each cluster, has been violated in \eqref{x2pm}, \eqref{xshift}, i.e.\ in covering space $\boXi$.
The violation does not occur at right rebellions $s\!=\!-$ with target $\boTh_+$\,.
It only occurs at left rebellions $s\!=\!+$, within the fat cluster $J_+^b$ at target equilibrium $\tilde{\boTh}_+$\,.
See $\bolth_{**}$ versus $\bolth_*$ in figure \ref{fig52}.
For left rebellions $\bolth_*\!\in\!\boXi$, indeed, the $x_2$-asymptotics \eqref{x2pm} forces a discrepancy of $2\pi$ among the phases $\vth_j$ of the fat target cluster $J_+^b=J_-^b\cup J_2$\,.
The discrepancy tracks whether $j\!\in\! J_+^b$ originated from the majority cluster $J_-^b$ of the source $\boTh_-$\,, or from the rebel cluster $J_2$\,; see \eqref{x2pm} and figure \ref{fig22}.
In particular, the targets $\boTh_+$ and $\tilde{\boTh}_+$ of a right and left rebellion, respectively,
\emph{belong to different fundamental parallelograms}, under $\Lambda'$.
When passing to the torus $\bT^2\!=\!\boXi/\Lambda'$ of theorem \ref{thm3clu}, this subtlety requires a reset of left rebellion phases $x_2$ to return to the original parallelogram.
In $\bT^2\!=\!\boXi/\Lambda'$, the discrepancy then disappears and \eqref{xm} is restored by the phase reset $\tilde{\boTh}_+\mapsto\boTh_+$\,.

Consider the phase orderings \eqref{right}, \eqref{left} of transients in theorem \ref{thm3clu}, and \eqref{xorder} in theorem \ref{thmMclu} below, as well as the phase increments \eqref{xshift} in corollary \ref{xshift}, and \eqref{sumxshift} in corollary \ref{shift-type}.
As is appropriate when working in covering space $\boXi$, these are all invariant under simultaneous translations of $\boTh_\pm$ and $\bolth_*$ by any $\bopsi\!\in\!\Lambda'$.
\end{rem}

We address heteroclinic orbits $\bolth_*{:}\,\ \boTh_-\!\leadsto\! \boTh_+$ between hyperbolic equilibria $\boTh_\pm$ next.
We claim that the associated intersections $W_-^u\cap W_+^s$ between their unstable and stable manifolds are transverse, automatically, in the sense of section \ref{Het}.
See section \ref{Trans} for a proof.

\begin{thm}[Transverse heteroclinics]
\label{thmtrv}
Let $\bolth_*(t){:}\, \ \boTh_-\!\leadsto\! \boTh_+$ denote any heteroclinic orbit, in covering space $\boXi$, between two equilibria $\boTh_\pm$ with order parameters $0\!<\!R_-\!<\!R_+\!\leq\! 1$.
Then the stable and unstable manifolds $W_-^u$ and $W_+^s$ intersect transversely, in $\boXi$, along $\bolth_*(t)$\,:
\begin{equation}
\label{hettrv}
\boTh_-\mathrel{\mathop{\leadsto}_{\transv}^{\bolth_*(t)}}\boTh_+\,.
\end{equation}
\end{thm}

The significance of such \emph{transverse heteroclinicity} is two-fold. 
\emph{First}, Morse--Smale systems are $C^0$ structurally stable, in compact settings. 
The \emph{Morse--Smale property} requires hyperbolic equilibria and hyperbolic periodic orbits, as the only recurrence, and transversality of all heteroclinic orbits among them.
\emph{Structural stability} then asserts that $C^1$-small perturbations lead to $C^0$ orbit equivalent flows, i.e., (flow)-orbits are mapped to orbits under a global near-identity homeomorphism.
See \cite{Palis,PalisSmale,PalisdeMelo,HMO} for details. 
We discuss perturbations of the Kuramoto model \eqref{kurODE} in section \ref{Pert}.

In co-rotating coordinates \eqref{kurODE}, and under the normalization \eqref{S1} in the flow-invariant covering space $\boXi\!=\!\mathbf{1}_N^\bot$ of the torus $\bT^{N\!-\!1}\!=\!\boXi/\Lambda'$, the Morse--Smale property on the positively invariant subset of $R \!>\! 0$ follows from our results above.
For hyperbolicity see proposition \ref{Equil} and theorem \ref{thmmon}.
Theorem \ref{thmtrv} asserts transversality of all heteroclinic orbits within the \emph{maximal compact invariant} subset $\cA $ of $\{\bolth\!\in\!\bT^{N\!-\!1}\,|\, R(\bolth)\!>\!0\}$, commonly called the \emph{global attractor}.
Projecting unstable manifolds $W^u(\boTh)\subset\boXi$, of equilibria $\boTh$ with $R>0$, to $\bT^{N\!-\!1}$, we may equivalently define 
\begin{equation}
\label{cA}
\cA =\bigcup_{\mathbf{f}(\boTh)\!=\!0} W^u(\boTh)\ \subset\ \bT^{N\!-\!1}\,.
\end{equation}
For our finite-dimensional setting of the normalized Kuramoto model \eqref{kurODE} in $\boXi$, we borrow the attractor terminology from infinite-dimensional dynamics.
The Morse--Smale property of \eqref{kurODE} then implies structural stability on $\cA $ in $\bT^{N\!-\!1}$.
See \cite[chapter 6]{HMO}.
The attractor $\cA$ and the repeller $\{R\!=\!0\}$ also form an \emph{attractor-repeller pair}, in the sense of Conley \cite{Conley, Mischaikow02}.

In full space $\bR^N$ and $\bT^N=\bR^N/\Lambda,\ \Lambda=2\pi\bZ^N$, however, the equilibria become relative equilibria \eqref{releq} (here: \emph{frozen waves}) under the group action of $\bR\mathbf{1}_N$\,.
They form normally hyperbolic stationary lines, alias circles $\mathrm{mod}\, 2\pi$. 
In principle, this violates the Morse--Smale paradigm of structural stability. 
Nevertheless, we may return to the original Kuramoto model \eqref{omegaj} with identical nonzero rotation frequencies $\omega_j\!=\!\omega\!\neq\!0$, as we did at the end of section \ref{Equi}.
Then the equilibrium circles remain relative equilibria, but revert to their original manifestation as hyperbolic periodic orbits, for $\omega\!\neq\!0$ (here: \emph{rotating waves}).
See \cite{FieHopf} for this terminology. 
Unstable and stable manifolds of the periodic orbits in \eqref{omegaj} then are just the group orbits, under $\bR\mathbf{1}_N$\,, of their intersections with $\boXi$.
In particular, transversality lifts to $\bR^N$. Therefore, structural stability extends to the lifted attractor 
\begin{equation}
\label{cAlift}
(\cA +\bR\mathbf{1}_N)/\Lambda\ \subset\ \bT^N\,.
\end{equation}
This proves the following corollary.

\begin{cor}[Morse--Smale structural stability]\label{strucstab}\quad
The classical Kuramoto model \eqref{omegaj} with identical frequencies $\omega_j\!=\!\omega\!\neq\! 0$ is Morse--Smale, within the open, positively invariant subset of $\{\bolth\!\in\!\bT^N\,|\, R(\bolth)\!>\!0\}$\,. 
In particular, the flow on the lifted compact global attractor \eqref{cAlift} of \eqref{omegaj} is $C^0$ structurally stable under $C^1$-small perturbations.
\end{cor}

We now return to the gradient setting of the co-rotating Kuramoto model \eqref{kurODE} with angle normalization \eqref{S1}.
The \emph{second} consequence of transversality theorem \ref{thmtrv} is \emph{transitivity} of the heteroclinicity relation. 
More precisely, suppose we have a sequence of transverse heteroclinic orbits
\begin{equation}
\label{trans}
\boTh_0\mathrel{\mathop{\leadsto}_{\transv}^{\bolth_1(t)}}\boTh_1\mathrel{\mathop{\leadsto}_{\transv}^{\bolth_2(t)}}\ldots   \mathrel{\mathop{\leadsto}_{\transv}^{\bolth_n(t)}}\boTh_n
\end{equation}
between successive hyperbolic equilibria $\boTh_0,\ldots,\boTh_n$ in covering space $\boXi$.
Then, arbitrarily close to the union of $\bolth_\ell(t)$ and their rest stops $\boTh_\ell$\,, there exists a direct transverse heteroclinic orbit
\begin{equation}
\label{concat}
\boTh_0\mathrel{\mathop{\leadsto}_{\transv}^{\bolth_*(t)}}\boTh_n\,.
\end{equation}
This is a consequence of the $\lambda$-Lemma \cite{PalisdeMelo}.
Actually, Lin's method provides a refinement: for any sufficiently large $T_* \!>\! 0$\,, \emph{passage times} $T_1,\ldots,T_{n\!-\!1}\geq T_*$ of $\bolth_*(t)$ near the rest stop equilibria $\boTh_1,\ldots,\boTh_{n\!-\!1}$ can be prescribed arbitrarily. 
See \cite{Lin,Linschol} for details.

Let us describe this transitivity, in the language of shadowing.
For some background, including the relation to structural stability, see \cite{Sauer, Tikhomirov} and the many references there.
We first represent each heteroclinic orbit $\bolth_\ell(t)$ by a uniformly $\delta$-close trajectory piece.
We call such a sequence of trajectory pieces, with jumps of at most $2\delta$ near each equilibrium, a $\delta$ \emph{pseudo-orbit}.
An $\eps$-\emph{shadow} is a single trajectory, which remains uniformly $\eps$-close to the pseudo-orbit.
See also \cite{Palmer} for detailed definitions.
Then the \emph{shadowing property} holds, in our setting: 
for every $\eps\!>\!0$, there exists some $\delta\!>\!0$, such that every $\delta$ pseudo-orbit possesses an $\eps$-shadow.
We call any such $\eps$-shadow $\bolth_*(t)$ a \emph{concatenation} of the successive heteroclinic steps $\bolth_1(t),\ldots,\bolth_n(t)$\,.

As a consequence, we obtain the following corollary to theorem \ref{thmtrv}.
See figures \ref{fig22} and \ref{fig23} for illustrations.

\begin{cor}[Concatenation of heteroclinic orbits]
\label{corconcat}
Consider any sequence of 1- and 2-cluster equilibria $\boTh_\ell\!\in\!\boXi$ with ascending order parameters
\begin{equation}
\label{Rseq}
0<R(\boTh_0)<\ldots<R(\boTh_n)\leq 1\,.
\end{equation}
As in \eqref{trans}, assume any sequence of heteroclinic orbits $\bolth_\ell\!\in\! \boXi$ between them.
The $\bolth_\ell$ are transverse, by theorem \ref{thmtrv}, but not necessarily of 2- or 3-cluster type.

Then, for any $\eps\!>\!0$, there exists an $\eps$-shadow of the sequence by a single transverse concatenated heteroclinic orbit
\begin{equation}
\label{Rdirect}
\boTh_0\mathrel{\mathop{\leadsto}_{\transv}^{\bolth_*(t)}}\boTh_n\,.
\end{equation}
For $n\!\geq\! 2$ and small $\eps\!>\!0$, the concatenations $\bolth_*(t)\!\in\! \boXi$ cannot be of 3-cluster type.
\end{cor}

The converse phenomenon of \emph{cascading} is a consequence of theorems \ref{thmmon}--\ref{thm3clu}, as follows.
Consider any two 1- or 2-cluster equilibria $\boTh_\pm\!\in\!\bT^{N\!-\!1}=\boXi/\Lambda'$ with Morse indices $0\!\leq\! i_+\!=\!i_-\!-\!n$ and any heteroclinic orbit $\bolth_*(t){:}\, \ \boTh_-\!\leadsto\!\boTh_+$ between them.
We do not require $\bolth_*(t)$ to be of 3-cluster type.
The fat clusters of $\boTh_\pm$ satisfy $J_+^b\supsetneq J_-^b$\,; see theorem \ref{thmmon}(iii).
Theorem \ref{thm3clu} asserts that there also exists a 3-cluster heteroclinic orbit between the same equilibria.

Next consider any sequence
\begin{equation}
\label{ascent}
J_-^b=J_0^b\,\subsetneq\,\ldots\,\subsetneq\, J_n^b=J_+^b
\end{equation}
of $n\!+\!1$ fat clusters $J_\ell^b$\,, which strictly ascend by adding one single oscillator, at a time.
Let $\boTh_\ell\!\in\! \bT^{N\!-\!1}\!=\!\boXi/\Lambda'$ denote the unique 1- or 2-cluster equilibria with those fat clusters, and note their descending Morse indices
\begin{equation}
\label{thj}
i_\ell:=i(\boTh_\ell)=i_--\ell\,.
\end{equation}
By the same reasoning as above, this induces a sequence \eqref{trans} of $n$ transverse heteroclinic orbits $\bolth_\ell(t)\!\in\!\boXi/\Lambda'$, with $J_{\ell-1}^b$ as fat cluster.
Each transverse heteroclinic orbit $\bolth_\ell(t)$ increases the size of the fat cluster by 1, and hence is a one-man rebellion with size $N_2\!=\!1$ of the rebellious cluster $J_2$\,.
In particular, each one-man rebellion $\bolth_\ell(t)$ is a 3-cluster.
By transversality and \eqref{ipm}, each rebellion orbit $\bolth_\ell(t)$ is isolated.
The term ``cascading'' describes the above ``interpolation'' phenomenon by a \emph{cascade} $\bolth_\ell(t)$ of one-man rebellions dropping the Morse index by 1, each.

As a first caveat, we recall remark \ref{resetTheta}: the condition \eqref{xm} of identical variables $\vth_j\!=\!x_m$ on each cluster, which holds on the torus $\bT^{N\!-\!1}\!=\!\boXi/\Lambda'$, is violated for the concatenations in covering space $\boXi$, as soon as left rebellions $\bolth_\ell$ with $s_\ell\!=\!+$ occur; see \eqref{x2pm}, \eqref{spm}.
See also \eqref{xshift} and the drift of variables $x_1(t), x_3(t)$ in figure \ref{fig22}.
Continuity of the left rebellion $\bolth_*(t)\rightarrow\tilde{\boTh}_+$ in figure \ref{fig52}, for example, forces us to leave the fundamental parallelogram of $\Lambda'\!\leq\!\boXi$.

As a second caveat, we add that cascading may fail in general Morse systems.
The height flow on the standard $n$-sphere $\mathbb{S}^n\subset\bR^{n\!+\!1}$ is a counterexample.
Rather, cascading is a peculiar feature of the Kuramoto model \eqref{kurODE}.
In the technically very different context of scalar one-dimensional PDEs of reaction-drift-diffusion type, cascading has first been observed in \cite{cascading}; see also \cite{ThomSmale,firoSFB,firoFusco,rofiln25} for later illustrations.

Combining theorems \ref{thmmon}--\ref{thmtrv}, we can now describe the \emph{connection graph} $\cC$ of the Kuramoto model \eqref{kurODE} for order parameters $0\!<\!R\!\leq\! 1$ and under the phase constraints \eqref{xm}, \eqref{S1}.
By proposition \ref{Equil}(i)--(ii), all equilibria $\boTh_{J^b}\!\in\!\bT^{N\!-\!1}\!=\!\boXi/\Lambda'$ with $R\!>\!0$ are hyperbolic, and uniquely characterized by their fat cluster set $J^b$ of size fraction $\alpha^b=|J^b|/N=(R\!+\!1)/2>1/2$; see section \ref{2Clu}.
These $\boTh_{J^b}$\,, alias $J^b$, serve as the \emph{vertices} of $\cC$.
By the Lyapunov property \eqref{Lyap}, all other orbits $\bolth(t)\!\in\!\bT^{N\!-\!1}$, whose closure remains in the region of $0 \!<\! R \!\leq\! 1$, are heteroclinic among those equilibria.
\emph{Edges} of $\cC$ represent such heteroclinicity.

The description of global dynamics by connection graphs has been common in variational settings of Morse type.
For many examples involving PDEs and delay equations with nodal properties see \cite{brfi88,brfi89, cascading, raugel, ThomSmale, firoSFB, firoFusco, rofiln25} and the references there.
Automatic transversality \cite{henrytrv}, due to nodal properties, also played a major role there, as well as in certain delay equations.
For general, not necessarily variational settings related to Conley index, see \cite{Conley,franzosa,Mischaikow02,Mischaikow}.
For recent adaptations to more general settings of global dynamics see also \cite{Yorke-a, Yorke-b, deleo}.

\begin{cor}[Connection graph]
\label{corconngraph}
The connection graph $\cC$ among the equilibria $\boTh_{J^b}\!\in\!\bT^{N\!-\!1}$ with $R\!>\!0$ is given by the partial order of the fat subsets $J^b\subseteq\{1,\ldots,N\}$ with respect to inclusion, as edges.
Heteroclinic orbits orient each edge towards the larger index set $J^b$.
\end{cor}

\begin{figure}
\centering\includegraphics[width=0.8\linewidth]{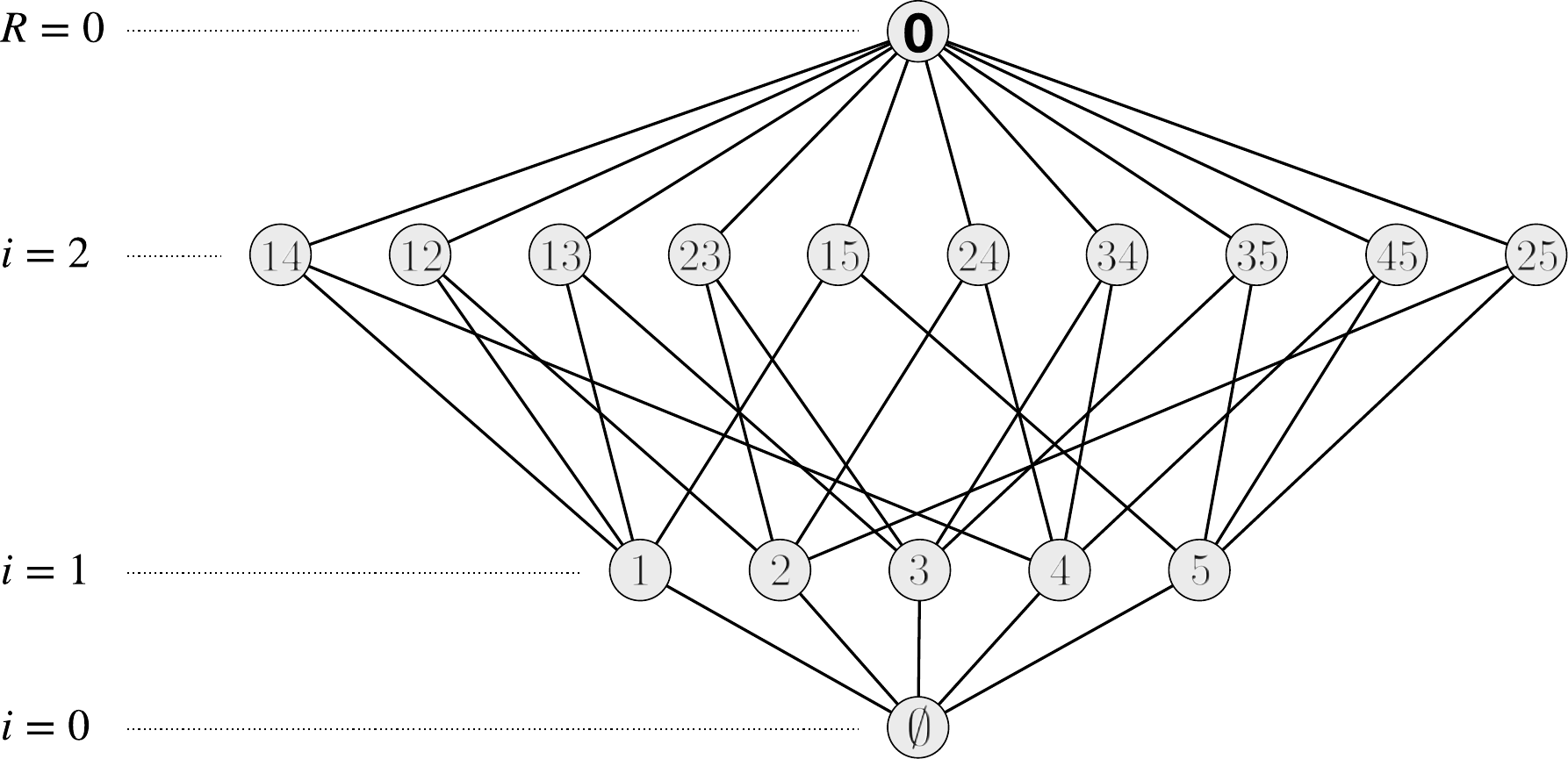}
\caption[Connection graph of five oscillators]{\emph{
Connection graph $\cC$ for the Kuramoto model \eqref{kurODE} with $N \!=\! 5$ oscillators.
Equilibrium vertices $\boTh_{J^b}\!\in\!\boXi/\Lambda'$ are represented by their complementary slim cluster sets $J^c$ of size $|J^c|\!=\!5\!-\!|J^b|$.
Note order parameters $0<R=2\alpha^b\!-\!1=1\!-\!2|J^c|/5\leq1$, and Morse indices $i\!=\!|J^c|$.
The 5-bar linkage equilibrium repeller $R\!=\!0$ of \eqref{R=0} is represented by the artificial top vertex $\mathbf{0}$.
The lowest level vertex $\emptyset$, of the empty slim cluster at $R\!=\!1$, denotes locally attracting total synchrony.
Directed edges represent heteroclinicity between adjacent Morse levels.
In particular, edges are directed downward, towards lower Morse index $i$ and larger order parameter $0 \!<\! R \!\leq\! 1$.
In view of cascading and transverse transitivity, directed paths in $\cC$ capture all heteroclinicity.
}}
\label{fig24}
\end{figure}

We follow Conley's concept of attractor-repeller pairs \cite{Conley, Mischaikow02}, and adjoin an artificial top vertex $\mathbf{0}$ to represent the repelling $N$-bar linkage equilibrium set $R\!=\!0$ of \eqref{R=0}.
For even $N$, the artificial vertex $\mathbf{0}$ includes the half-size pairs $|J^b|\!=\!|J^c|\!=\!N/2$.
For general $N$, we can thus describe the connection graph $\cC$ based on the well-known \emph{subset lattice} of $\{1,\ldots,N\}$, defined by set inclusion.
We just have to identify complementary pairs $J^b,J^c$, and direct the remaining edges downward, towards smaller slim representatives $J^c$.
By cascading, transversality theorem \ref{thmtrv}, and the resulting transitivity corollary \ref{corconcat} on heteroclinicity, it is sufficient to represent the connection graph by its edges between fat cluster sets of adjacent size. 
These are given by the one-man rebellion variants $|J_2|\!=\!1$ from 3-cluster theorem \ref{thm3clu}.
The number $V$ of vertices in the lattice is therefore given by
\begin{equation}
\label{vertexcount}
V=1+\sum_{N_1=\lfloor N/2\rfloor+1}^{N} \binom{N}{N_1}=
\begin{cases}
     2^{N-1}+1\,, & \textrm{for } N \textrm{ odd}, \\
     2^{N-1}+1-\tfrac{1}{2}\tbinom{N}{N/2}\,, & \textrm{for } N \textrm{ even}.
\end{cases}\end{equation}
Here $N_1\!=\!|J^b|\!=\!\alpha^b N\!>\!N/2$, and $\lfloor\cdot\rfloor$ denotes the floor function.
See figure \ref{fig24} for the example $N \!=\! 5$.
To shorten vertex annotation, we have replaced fat vertex labels $J^b$ by their slim complements $J^c$ there.
Heteroclinic orientation is downward, towards increasing order parameter $R=2\alpha^b\!-\!1$ and decreasing Morse index $i\!=\!|J^c|\!=\!(1\!-\!\alpha^b)N$.

Note how corollary \ref{corconngraph} also follows from theorem \ref{thm3clu}, directly, without the long detours on transversality, shadowing, and concatenation.
Indeed, the elementary building brick of 3-cluster heteroclinicity realizes \emph{all} heteroclinicity within $R\!>\!0$.
In particular, it realizes any downwards path in the part $R\!>\!0$ of the connection graph $\cC$, by a direct heteroclinic 3-cluster ``express'' from the start vertex to the terminal vertex of the path.
The 3-cluster express skips all Lin loitering \cite{Lin,Linschol} at intermediate rest stop equilibria along the path.

Transversality, on the other hand, comes with the added bonus of shadowing, alias heteroclinic concatenation.
As an immediate consequence of the 3-cluster heteroclinic left/right alternatives (b), (c) of figure \ref{fig21}, the concatenations $\bolth_*(t)$ of corollary \ref{corconcat} can be composed from an arbitrary sequence of left/right one-man rebellions $N_2\!=\!|J_2|\!=\!1$. 
See theorem \ref{thm3clu}, corollary \ref{corxshift}, and figures \ref{fig22} and \ref{fig23}.

\begin{cor}[Finite $2$-shift of heteroclinic concatenations]
\label{shift-type}
For $n\!\geq\! 1$, let $\boTh_\pm\!\in\!\bT^{N\!-\!1}\!=\!\boXi/\Lambda'$ denote any two 1- or 2-cluster equilibria with sizes $N_+^b\!=\!N_-^b\!+\!n$ of their fat clusters $J_-^b\subset J_+^b$\,, and corresponding order parameters $0\!<\!R_-\!<\!R_+ \!\leq\! 1$.
Prescribe any finite sequence 
\begin{equation} 
\label{sequence}
\mathbf{s}=s_1s_2\ldots s_n
\end{equation}
of $n$ sign symbols $s_\ell\!\in\!\{\pm\}$.
Consider any associated cascading sequence \eqref{trans} of heteroclinic orbits $\bolth_\ell(t)\,{:}\,\ \boTh_{\ell-1}\!\leadsto\!\boTh_{\ell}$\,, between equilibria $\boTh_-\!=\!\boTh_0\,,\ldots,\boTh_n\!=\!\boTh_+$ with ascending fat clusters of sizes $N_-^b\!+\!\ell$.
Each heteroclinic orbit $\bolth_\ell(t)$, considered in $\boXi/\Lambda'$, increases the size of the fat cluster by 1, and hence is an isolated one-man rebellion with size $N_2\!=\!1$ of its rebellious cluster $J_2$\,.
In particular, each $\bolth_\ell(t)$ is a 3-cluster orbit.
Each $\bolth_\ell(t)$ can be chosen of left or right type (b) or (c), as prescribed by symbol $s_\ell\!=\!+$ or $s_\ell\!=\!-$, respectively; see \eqref{spm} and corollary \ref{corxshift}.

Corollary \ref{corconcat} then provides a total of $2^n$ distinct transverse concatenations $\bolth_*(t)\!\in\!\boXi/\Lambda'$, one for each left/right sign sequence $\mathbf{s}$.
For small $\eps\!>\!0$, the shadowing concatenations are not related by any permutation or phase-shift symmetry in $S_N\!\ltimes\!\Lambda'$; compare \eqref{groupS}, \eqref{groupT}.
Locally, each concatenation generates an $(n\!-\!1)$-dimensional family of heteroclinic orbits from $\boTh_-$ to $\boTh_+$\,.
According to \eqref{xshift}, the positions $x_1\!=\!\vth_j$ in the initial fat cluster $J_0$ shift by a cumulative total of
\begin{equation}
\label{sumxshift}
x_+^b-x_-^b= -\tfrac{\pi}{N}\sum_{\ell=1}^{n} s_\ell 
\end{equation}
in covering space $\boXi$; see remark \ref{resetTheta}.
\end{cor}

For $n\!\geq\! 2$ and small $\eps\!>\!0$, the concatenations $\bolth_*(t)$ are not realizable by direct, 3-cluster, heteroclinic express orbits, which would require a rebellion of cluster size $N_2\!=\!n$\,. 
Mapping $\bolth\mapsto-\bolth$ by $\sigma\!=\!-\mathbf{id}_N$, as in \eqref{groupS}, swaps cases (b), (c) of figure \ref{fig21}.
For example, we will encounter the swap between $\bolth_*$ and $\bolth_{**}$ in section \ref{3Clu}; see the proof of theorem \ref{thm3clu} and figure \ref{fig52}.
Upon concatenation \eqref{sequence}, however, the swap reverses all signs $s_\ell$\,, in unison.
This reduces the total count $2^n$ of corollary \ref{shift-type} to $2^{n\!-\!1}$, but does not affect the exponential growth of available shift-sequences $\mathbf{s}$, with $n$.

Of course we may formulate variants which involve $n'\!<\!n$ symbols $s_\ell$ and concatenate 3-cluster rebellions which also admit larger rebellious left/right clusters $J_{\ell,2}$ of general sizes $N_{2,1}+\ldots+N_{2,n'}=n$\,.
This provides a total count of $2\cdot 3^{n\!-\!1}$ heteroclinic concatenations $\bolth_*(t){:}\, \ \boTh_0\!\leadsto\!\boTh_n$ which are not related by any permutation or phase-shift symmetry. Based on the following theorem \ref{thmMclu}, we may also split each rebel cluster set $J_{\ell,2}$ further, into swarms of subclusters, left/right subclusters, and their concatenations, as sizes permit.

\begin{figure}[t]
\centering \includegraphics[width=\textwidth]{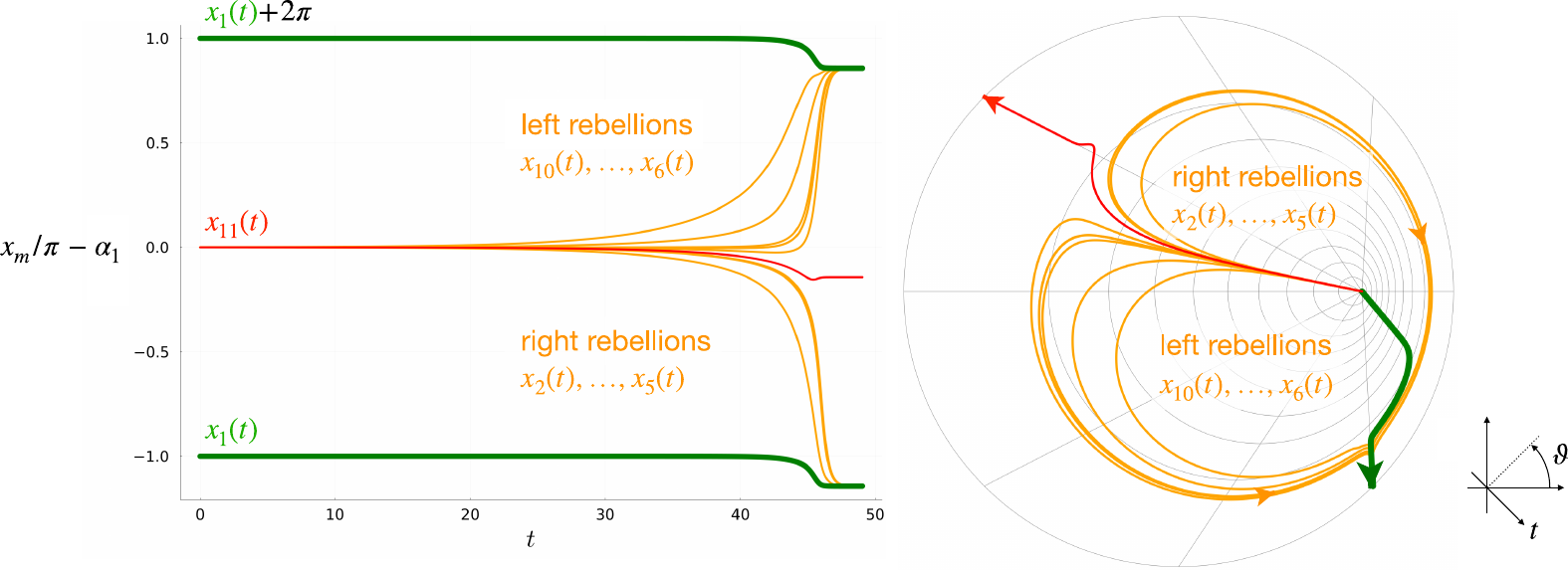}
\caption[Swarm rebellion]{\emph{
Time plot of a swarm rebellion, in covering space $\boXi$, from a 2-cluster equilibrium $\boTh_-$ at $t \to -\infty$ towards an unstable 2-cluster equilibrium $\boTh_+$ with larger fat cluster size $\mathrm{mod}\ \Lambda'$. 
See theorem \ref{thmMclu} and section \ref{NumSwarm}.
Color coding and scaling are the same as in figures \ref{fig22}, \ref{fig23}.
Fat cluster sizes are $N_-^b\!=\!11,\ N_+^b\!=\!20$ for $N\!=\!21$ oscillators $\bolth(t)$, moving in $M\!=\!11$ clusters.
Morse indices are $i(\boTh_-)\!=\!10$ and $i(\boTh_+)\!=\!1$.
Limiting cluster phases $x_1/\pi\!-\!\alpha_1\!=\!\pm1$ and $x_2/\pi\!-\!\alpha_1\!=\!\ldots\!=\!x_{11}/\pi\!-\!\alpha_1\!=\!0$ for $t\rightarrow -\infty$ indicate the same source equilibrium $\boTh_-$ as in figures \ref{fig22}, \ref{fig23}, with the same initial size fraction $\alpha_1:=\alpha_-^b\!=\!11/21$.
Of the $M\!=\!11$ clusters, $M\!-\!2\!=\!9$ are a rebellious swarm of one-man ``clusters'', preserving their ordering by phases; see \eqref{xorder}. 
The swarm of nine is parted into $m_*\!-\!1\!=\!4$ right (downwards) and $M\!-\!1\!-\!m_*\!=\!5$ left (upwards)  rebels. 
The final phases at (shifted) target $\boTh_+$ satisfy $x_1/\pi\!-\!\alpha_1\!=\!\pm1\!-\!1/21$ and $x_{11}/\pi\!-\!\alpha_1\!=\!-1/21$; see corollary \ref{corxshift} and remark \ref{resetTheta}.
Also compare cumulative shifts \eqref{sumxshift} of $x_1$ in concatenation corollary \ref{shift-type}.}}
\label{fig25}
\end{figure}

In fact, the general heteroclinic orbit $\bolth_*(t){:}\, \ \boTh_-\!\leadsto\! \boTh_+$ is neither a 3-cluster, nor a concatenation of one-man rebellions, for small $\eps$.
We therefore refine theorems \ref{thmmon}, \ref{thm3clu}, and describe the \emph{swarm behavior} of general rebellions.
See section \ref{kurGroup} for notation, and figure \ref{fig25} for illustration.

\begin{thm}[Realization by $M$-cluster swarm rebellions]\label{thmMclu}
Let $\bolth_*(t){:}\, \ \boTh_-\!\leadsto\! \boTh_+$ denote any heteroclinic orbit in $\bT^{N\!-\!1}\!=\!\boXi/\Lambda'$ of the co-rotating Kuramoto model \eqref{kurODE}, between two equilibria $\boTh_\pm$ with order parameters $0\!<\!R_-\!<\!R_+\!<\!1$ and fat cluster sets $J_\pm^b$\,, respectively.
Then $\bolth_*(t)$ possesses $M\!\geq\! 3$ clusters.
More precisely,
\begin{equation}
\label{hetFix}
\bolth_*(t)\in\mathrm{Fix}(S_{J_-^b})\cap\mathrm{Fix}(S_{J_+^c})
\end{equation}
inherits the fat cluster $J_1:=J_-^b$ of the source $\boTh_-$ and the slim cluster $J_M:=J_+^c$ of the target $\boTh_+$\,, of sizes $N_-^b\!=\!\alpha_1N\!=\!\alpha_-^bN$ and $N\!-\!N_+^b\!=\!\alpha_+^cN$\,, respectively.
Their complementary rebel set $\hat{\mathbf{J}}:=J_2\dot\cup\ldots\dot\cup J_{M\!-\!1}$ of size $N_+^b\!-\!N_-^b\!\geq\! 1$ collects the remaining $1\!\leq\!M\!-\!2\!\leq\! N_+^b\!-\!N_-^b$ clusters.
In other words, $\bolth_*(t)\in \bT_\mathbf{J}^{M\!-\!1}=\boXi_\mathbf{J}/\Lambda_\mathbf{J}'$\,; see \eqref{TJM-1}.

The locations $x_1\!<\!x_m\!<\!x_1\!+\!2\pi$ of those remaining clusters $J_m\,,\ 1\!<\!m\!<\!M$, can be separated by $x_M$\,. For some $m_*$\,, their phases $x_m$ can therefore be labeled by $m$ such that
\begin{equation}
\label{xorder}
x_1<x_2<\ldots<x_{m_*}<x_M<x_{M\!-\!1}<\ldots<x_{m_*\!+\!1}<x_1+2\pi\,.
\end{equation}
See \eqref{right}, \eqref{left}, and figure \ref{fig21}.
We call $\bolth_*$ an $M$\emph{-cluster swarm rebellion}.

Conversely, any such $M$-cluster swarm rebellion and ordering of phases is realized by certain (families of) heteroclinic orbits $\bolth_*(t)\in \bT_\mathbf{J}^{M\!-\!1}=\boXi_\mathbf{J}/\Lambda_\mathbf{J}'$\,.
\end{thm}

See section \ref{Trans} for a proof.
In figure \ref{fig25} we have chosen initial data $\bolth_*(0)\ \in\ \mathrm{Fix}(S_{J_-^b}) \cap \mathrm{Fix}(S_{J_+^c})$ near $\boTh_+$\,.
Picking suitable perturbations supported on the swarm indices $\hat{\mathbf{J}}$, the numerical solution converges to a source 2-cluster equilibrium $\boTh_-$ in backwards time.
See section \ref{NumSwarm} for further comments on the precise choice \eqref{part21eps} of perturbations, the assignment of left and right rebels by $m_*$\,, and actual numerical implementation.

\section{Partitioned Coordinates}\label{Part}

We set up some general notation to prove proposition \ref{Equil} and theorems \ref{thmmon}(iii)--\ref{thmMclu}.
From section \ref{kurGroup} we recall the orthogonal decomposition 
\begin{equation}
\label{1+Xi}
\bolth=x\mathbf{1}_N\oplus\boxi\,\in\, \mathbf{X}= \bR\mathbf{1}_N \oplus \boXi\,,
\end{equation} 
of the covering space $\bolth\!\in\!\mathbf{X}\!=\!\bR^N$ of $\bT^N\!=\!\mathbf{X}/\Lambda$, into irreducible representations of the permutation group $\sigma\!\in\! S_N$\,.
Here $x$ is the average of $\bolth$, and $\boxi\!=\!\mathbf{P}\bolth\!\in\!\boXi\!=\!\mathbf{1}_N^\bot$ possesses zero average.
Note $\bR\mathbf{1}_N\!=\!\mathrm{Fix}(S_N)$ is the space of $S_N$-fixed vectors in $\mathbf{X}$.
The representation of $S_N$ on the orthogonal complement $\boXi\!=\!\mathbf{1}^\perp$ is called the \emph{standard irreducible representation} of $S_N$\,.
This representation is known to be \emph{absolutely irreducible}, i.e.~ the only $S_N$-equivariant matrices on $\boXi$ are \emph{real} multiples of the identity matrix. See for example \cite{Elm1,StElm,Elm2}. 
In other words, Schur's Lemma holds over $\bR$ rather than just over $\mathbb{C}$.
The angle normalization \eqref{S1}, for example, leads to $x\!=\!0$ in our coordinates, so that $\boxi\!\in\!\boXi\!\cong\!\bR^{N\!-\!1}$ measures the deviation from total synchrony.

Next, we extend the above notation to arbitrary index partitions $\mathbf{J}\!=\!(J_m)$ of the $N$-oscillator index set into disjoint subsets $J_m$\,, for $m=1,\ldots,M$.
See \eqref{J}.
The index partition induces the following orthogonal direct sum decomposition of $\mathbf{X}\!=\!\bR^N$\,:
\begin{align}
\label{Xj}
   \mathbf{X} &= \mathbf{X}_1\oplus\ldots\oplus\mathbf{X}_M\,,   \\
\label{Xij}
    \mathbf{X}_m &= \bR\mathbf{1}_m \oplus \boXi_m\,.
\end{align}
Here $\mathbf{x}_m\!\in\!\mathbf{X}_m$ is only supported on $J_m$\,, i.e.\ the components of $\mathbf{x}_m$ are zero outside $J_m$\,.
Analogously to \eqref{1+Xi} we have decomposed $\mathbf{x}_m\!\in\!\mathbf{X}_m$ as
\begin{equation}
\label{1+Xim}
\mathbf{x}_m=x_m\mathbf{1}_m\oplus\boxi_m\,\in\, \mathbf{X}_m= \bR\mathbf{1}_m  \oplus \boXi_m\,,
\end{equation}
orthogonally.
The components $\mathbf{1}_{mj}$ are 1 for $j\!\in\! J_m$\,, and 0 otherwise.
The averages of $\mathbf{x}_m$ are therefore the scalar coordinates $x_m$\,, and each $\boxi_m$ averages to zero.
We rearrange the coordinates \eqref{Xj}--\eqref{1+Xim} to represent $\bolth\!\in\!\boXi$ by the $M$\emph{-partitioned coordinates}
\begin{equation}
\label{xxi}
(x_2,\ldots,x_M,\boxi_1,\ldots,\boxi_M)\,,
\end{equation}
with the built-in normalization $x_1:=-\alpha_1^{-1}\sum_{m\geq 2}\alpha_mx_m$ of \eqref{S1}.
For example, the normalized $M$-clusters with coordinates \eqref{xm} correspond to the case where all $\boxi_m$ vanish.
The 2-cluster dynamics \eqref{2cluODEJ} of section \ref{2Clu} correspond to coordinates $(x_2,\mathbf{0})$ with $x_2\!=\!x^c$.

As in \eqref{J}--\eqref{xm}, let $\sigma\!\in\! S_{J_m}\!\leq\! S_N$ denote the subgroup of permutations $\sigma$ which only permute indices  $j\!\in\! J_m$\,, leaving all other indices fixed.
The decomposition \eqref{1+Xim} is the orthogonal decomposition of $\mathbf{X}_m$ into irreducible representations of $S_{J_m}$\,.
Similarly, the decompositions \eqref{Xj}--\eqref{1+Xim} provide the orthogonal decomposition of $\mathbf{X}_m$ into irreducible representations of the direct product $S_\mathbf{J}$ of all $S_{J_m}$ for $m=1,\ldots,M$.
Note that all summands $\boXi_{m}$ are pairwise inequivalent, absolutely irreducible representations of $S_\mathbf{J}$\,. 
The span 
\begin{equation}
\label{FixSJ}
\mathbf{X_J}:=\mathrm{Fix}(S_\mathbf{J})=\mathrm{span}\{\mathbf{1}_m\,|\,1\!\leq\! m\!\leq\! M\}
\end{equation}
of all $\mathbf{1}_m$ features the trivial representation of $S_\mathbf{J}$\,. 

Let $S_{J_m}'$ denote the direct product of all $S_{J_{m'}}$ with $m'\!\neq\! m$.
Similarly, let $\mathbf{X}_m'$ and $\boXi_m'$ denote the direct sums of all $\mathbf{X}_{m'}$ and $\boXi_{m'}$ with $m'\!\neq\! m$, respectively.
Then
\begin{align}
\label{FSm}
\mathrm{Fix}(S_{J_m}) &= \mathrm{Fix}(S_\mathbf{J})\oplus\boXi_m'=\{\boxi_m\!=\!0\}\,;\\
\label{FS'm}
\mathrm{Fix}(S_{J_m}') &= \mathrm{Fix}(S_\mathbf{J})\oplus\boXi_m\,.
\end{align}
By $S_N$-equivariance, any such fix space is a time-invariant subspace of the normalized Kuramoto model \eqref{kurODE}; see \eqref{FixJ}.
In case you find this terminology impenetrable, kindly return to introductory sections \ref{Group}, \ref{kurGroup}, and the references there.

For illustration, we base the proof of proposition \ref{Equil} on linearizations at 1- and 2-cluster equilibria $\boTh_\pm$ on 2-partitioned coordinates $(x_2,\boxi_1,\boxi_2)$; see \eqref{xxi} with $M\!=\!2$.

\emph{\textbf{Proof of proposition} \ref{Equil}}.\quad
It only remains to establish hyperbolicity and the Morse indices $i\!=\!i(\boTh)\!=\!|J^c|$ in $\boXi$, at 1- and 2-cluster equilibria $\boTh\!\in\!\boXi/\Lambda'$ with $0\!<\!R\!\leq\!1$.

We start with the linearization $\mathbf{f}'(\boTh)$ at the unique 1-cluster equilibrium $\boTh\!=\!0$ under the phase constraint \eqref{S1}.
By \eqref{f'} at $\bolth\!=\!\boTh\!=\!0$, the linearization commutes with all permutations $\sigma\!\in\! S_N$\,.
Therefore, $\mathbf{f}'(\boTh)$ diagonalizes over the $S_N$-irreducible decomposition \eqref{1+Xi}.

The synchronous equilibria $\mathbf{f}(\bR\mathbf{1}_N)=0$ imply an eigenvalue zero of $\mathbf{f}'(0)$, due to $\mathbf{f}'(0)\mathbf{1}_N=0$.
Let $A\,{:}\,\ \boXi\rightarrow\boXi$ denote the restriction of $\mathbf{f}'(0)$ to the complementary invariant representation subspace $\boXi\!=\!\mathbf{1}_N^\bot$\,.
Irreducibility of $S_N$ on $\boXi$ implies $A=a\,\mathbf{id}_\boXi$ is a scalar multiple of identity on $\boXi$.
Hence
\begin{equation}
\label{am1}
a\ =\ \frac{1}{\dim \boXi}\cdot\mathrm{trace}\,\mathbf{f}'(0)\ =\ \tfrac{1}{N\!-\!1} \cdot\tfrac{1}{N}\sum_{j=1}^N\sum_{k\neq j} (-1)\ =\ -1\,.
\end{equation}
This proves claim (i) of proposition \ref{Equil}.

We address proposition \ref{Equil}(ii) next, and consider the linearization $\mathbf{f}'(\boTh)$ in $\mathbf{X}\!=\!\bR^N$ at any 2-cluster equilibrium $\boTh$ with $0<R=2\alpha^b\!-\!1<1$.
As in section \ref{2Clu}, let $J^b\!=\!J_1$ and $J^c\!=\!J_2$ indicate the partition $\mathbf{J}$ into the complementary fat and slim cluster sets.
On the invariant section $\boXi\!=\!\mathbf{1}_N^\bot$ of interest, we can invoke 2-partitioned coordinates $(x^c\,,\boxi^b,\boxi^c)$ for $\bolth$; see \eqref{xxi}.
Since the  absolutely irreducible representations of $S_\mathbf{J}=S_{J^b}\times S_{J^c}$ on $x^c\!\in\!\bR$ and on the components $\boXi_1,\boXi_2$ are inequivalent, respectively, we can rewrite the restricted linearization $\mathbf{f}'(\boTh)\vert_\boXi$ as a matrix $A$ in block diagonal form, with respect to these coordinates\,:
\begin{equation}
\label{Ablock}
A=\begin{pmatrix}
     a_0 &  0 & 0  \\
     0 &  a_b\,\mathbf{id}_{\boXi^b}\ & 0  \\
     0 & 0 & a_c\,\mathbf{id}_{\boXi^c}  
\end{pmatrix}\,.
\end{equation}
Here $a_0\!\in\!\bR$, and the real scalars $a_b\,,\,a_c$ multiply identity matrices on $\boXi^b,\boXi^c$\,.

For the linearization $a_0$ within the 2-cluster subspace $\mathbf{X_J}=\mathrm{Fix}(S_\mathbf{J} )$ of vanishing $\boxi_1,\boxi_2$, the 2-cluster dynamics \eqref{2cluODEJ} of $x^c$ on the line $\boXi_\mathbf{J}$ at the midpoint \eqref{midpoint} implies
\begin{equation}
\label{a02}
a_0=+1\,.
\end{equation}

To determine $a_b\,,\,a_c$\,, we fix $\boxi^\iota\!\in\!\boXi^\iota$ with $\xi_j\!=\!1$ for some $j\!\in\! J^\iota$ and superscript $\iota\!=\!b,c$.
Then the block diagonal structure of $A$ implies
\begin{equation}
\label{am2}
\begin{aligned}
a_\iota=f_j'(\boTh)\boxi^\iota&= \tfrac{1}{N} \sum_{k=1}^N \cos(\Th_k-\Th_j)(\xi_k^\iota-\xi_j^\iota)=-\alpha^\iota+(1-\alpha^\iota)\\
       &=
       \begin{cases}
      -2 \alpha^b+1=-R\,<\,0\,, & \text{ for } \iota\!=\!b, \\
      +2\alpha^b-1=+R\,>\,0\,, & \text{ for } \iota\!=\!c.
\end{cases}
\end{aligned}
\end{equation}
Here we have used $\sum\xi_k^\iota=0,\ x^{\iota'}\!-\!x^\iota=\pm \pi$ for $\iota' \!\neq\! \iota$, and $\alpha^b\!=\!1\!-\!\alpha^c$.

Together, \eqref{a02} and \eqref{am2} prove hyperbolicity of $\boTh$ in $\boXi$, and collect a Morse index 
\begin{equation}
\label{i2}
i(\boTh)=1+\dim(\boXi^c)=1+(|J^c|-1)=|J^c|\,.
\end{equation}
Indeed, $\boXi^c$ consists of zero-sum vectors whose components vanish on the complement $J^b$.
This proves claim (ii), and proposition \ref{Equil}. 
\qed

\section{Stable and unstable manifolds}\label{Wsu}
We have introduced stable and unstable manifolds, in section \ref{Het}.
We can now describe the stable and unstable manifolds of hyperbolic 2-cluster equilibria $\boTh$ with order parameter $0<R=2\alpha^b\!-\!1<1$; see lemma \ref{lemW}.
In lemma \ref{lemWpm}, we collect some consequences for general heteroclinic orbits $\bolth_*(t){:}\, \ \boTh_-\!\leadsto\! \boTh_+$  between 2-cluster equilibria $\boTh_\pm$\,.
In particular, lemma \ref{lemWpm} proves strict monotonicity $J_-^b\subsetneq J_+^b$ of fat clusters in $\bT^{N\!-\!1}\!=\!\boXi/\Lambda'$, i.e.\ theorem \ref{thmmon}(iii).
Our presentation relies on proposition \ref{Equil} and the partitioned coordinates \eqref{xxi} in the linear section $\boXi\!=\!\mathbf{1}_N^\bot$\,; see \eqref{P} and \eqref{S1}.

\begin{lem}[Linearization of invariant manifolds]
\label{lemW}
Let $\boTh\!\in\!\boXi$ be any 2-cluster equilibrium of the Kuramoto model \eqref{kurODE} with order parameter $0<R=2\alpha^b\!-\!1<1$.
Let $J^b\!=\!J_1$ and $J^c\!=\!J_2$ denote the fat and slim clusters of $\boTh$ with size fractions $1/2\!<\!\alpha^b\!<\!1$ and $\alpha^c\!=\!1\!-\!\alpha^b$, respectively.
In the $M\!=\!2$ partitioned coordinates $(x_2,\boxi_1,\boxi_2)=(x^c,\boxi^b,\boxi^c)$ of \eqref{xxi}, the stable and unstable manifolds $W^{s,u}$ of $\boTh$ in $\boXi$ then possess the tangent spaces
\begin{align}
\label{TWu}
T_\bolth W^u &= \mathrm{Fix}(S_{J^b})=\mathbf{X}_{J^b}'=\mathbf{X}_{J^c}=\bR\mathbf{1}^c\oplus\boXi^c=\{\boxi^b\!=\!0\}\,,\\
\label{TWsTh}
T_\boTh W^s &= \mathrm{Fix}(S_{J^c})\cap\{x^c\!=\!0\}\; =\boXi^b\,, \\
\label{TWs}
T_\bolth W^s &\leq \mathrm{Fix}(S_{J^c}) =\bR\mathbf{1}^c\oplus\boXi^b=\{\boxi^c\!=\!0\}\,,
\end{align}
at any $\bolth$ in $W^u$ and $W^s$, respectively. 
Since the linear spaces $\mathrm{Fix}(S_{J^\iota})$ are time-invariant, for $\iota\!=\!b,c$, this implies
\begin{align}
\label{Wu}
W^u &\subseteq \mathrm{Fix}(S_{J^b})=\bR\mathbf{1}^c\oplus\boXi^c=\{\boxi^b\!=\!0\}\,,\\
\label{Ws}
W^s &\subsetneq \mathrm{Fix}(S_{J^c})=\bR\mathbf{1}^c\oplus\boXi^b=\{\boxi^c\!=\!0\}\,.
\end{align}
In particular, $W^u$ is open in $\mathrm{Fix}(S_{J^b})$;
but $T_\bolth W^s$ and $W^s$ each are of codimension 1 in $\mathrm{Fix}(S_{J^c})$.
\end{lem}

\begin{proof}[\textbf{Proof}]
We first establish the linear claims \eqref{TWu}, \eqref{TWsTh} at the 2-cluster equilibrium $\bolth\!=\!\boTh$.
Next, we show how \eqref{TWu} at $\bolth\!=\!\boTh$ implies the nonlinear claim \eqref{Wu}.
Linearization at any $\bolth\!\in\! W^u$ then implies \eqref{TWu} everywhere.
We conclude with slight adaptations for \eqref{TWs} and \eqref{Ws}.

By proposition \ref{Equil}(ii), the 2-cluster equilibrium $\boTh$ is hyperbolic.
In \eqref{a02}, \eqref{am2} we have determined the unstable eigenvalues $a_0\!=\!+1$, in the 2-partitioned coordinate $x^c$ of $\bR\mathbf{1}^c$.
Moreover, $a_c=R=2\alpha^b\!-\!1>0>-R=a_b$\,, in the coordinates $\boxi^\iota\!\in\!\boXi^\iota$\,.
For variants of subspace notation, see also \eqref{FSm}, \eqref{FS'm}.
This proves the linear claims \eqref{TWu}, \eqref{TWsTh} at $\bolth\!=\!\boTh$.

To prove the nonlinear claim \eqref{Wu}, we rely on invariance of the \emph{linear space} $\mathrm{Fix}(S_{J^b})$ in \eqref{TWu} under the \emph{nonlinear flow} of \eqref{kurODE}; see \eqref{FixJ}.
We therefore obtain a local unstable manifold $W_\mathrm{loc}^u$ of $\boTh$ within $\mathrm{Fix}(S_{J^b})$, of full dimension.
By uniqueness of the nonlinear local unstable manifold at hyperbolic equilibria, $W_\mathrm{loc}^u$ therefore coincides with its flow-invariant tangent space, locally.
Continuation from $W_\mathrm{loc}^u$ to $W^u$, by the flow $\boPhi^t$, then proves claim \eqref{Wu}.
Linearization and invariance of dimension, at last, imply claim \eqref{TWu}, at any $\bolth\!\in\! W^u$.

To prove nonlinear claim \eqref{Ws}, we first have to confess that nonlinear flow-invariance of the representation subspace $\boXi^b$ may fail: it is not a fix space \eqref{FixJ} under a suitable subgroup $\mathrm{\bJ}$.
Moreover, the tangent space $T_\boTh W^s$ needs to be attached to $\boTh\!\in\!\bR\mathbf{1}^c$.
Working in the larger, but nonlinearly time-invariant, linear space $\mathrm{Fix}(S_{J^c})$\,, proves claim \eqref{Ws} by uniqueness of $W^s$.
Linearization then implies claim \eqref{TWs} and shows codimension 1.
This proves the lemma.
\end{proof}

\begin{lem}[Invariant manifolds along heteroclinic orbits]
\label{lemWpm}
Let  $\bolth_*(t){:}\, \ \boTh_-\!\leadsto\! \boTh_+$ denote any heteroclinic orbit in $\boXi/\Lambda'$ between equilibria $0<R(\boTh_-)=2\alpha_-^b\!-\!1<2\alpha_+^b\!-\!1=R(\boTh_+)<1$.
Then $\boTh_\pm\!\in\!\boXi/\Lambda'$ are 2-cluster equilibria, by theorem \ref{thmmon}(i), with fat clusters $J_\pm^b$ of size fractions $\alpha_\pm^b$\,, and complementary slim clusters $J_\pm^c$\,.
Moreover, the fat clusters grow strictly\,:
\begin{enumerate}[(i)]
  \item \qquad$\bolth_*(t)\in\mathrm{Fix}(S_{J_-^b})\cap\mathrm{Fix}(S_{J_+^c})$\,,\quad for all\, $t\!\in\!\bR$\,; 
  \item \qquad $\alpha_-^b<\alpha_+^b$ and $J_-^b\subsetneq J_+^b$\,.
\end{enumerate}
Let $(x_2,x_3,\boxi_1,\boxi_2,\boxi_3)$ denote the 3-partitioned coordinates \eqref{xxi} of $\bolth_*(t)$ for the resulting index partition $\mathbf{J}$ given by $J_1:=J_-^b\,,\ J_2:=J_+^b\setminus J_-^b \,,\ J_3:=J_+^c$\,; see \eqref{J123}.
Then
\begin{equation}
\label{xi130}
\boxi_1=0\ \textrm{and}\ \boxi_3=0.
\end{equation}
The tangent spaces $T_*W_-^u$ and $T_*W_+^s$\,, in $\boXi$, to the unstable and stable manifolds of $\boTh_\pm$ at any $\bolth_*\!=\!\bolth_*(t)$ satisfy
\begin{align}
\label{T*Wu}
T_*W_-^u &= \mathrm{Fix}(S_{J_-^b})=\mathbf{X}_2\oplus\mathbf{X}_3=\mathrm{span}\{\mathbf{1}_2, \mathbf{1}_3\}\oplus\boXi_2\oplus\boXi_3=\{\boxi_1\!=\!0\}\,,\\
\label{T*Ws}
T_* W_+^s &\leq \mathrm{Fix}(S_{J_+^c})=\mathbf{X}_1\oplus\mathbf{X}_2\oplus\bR\mathbf{1}_3=\mathrm{span}\{\mathbf{1}_2, \mathbf{1}_3\}\oplus\boXi_1\oplus\boXi_2=\{\boxi_3\!=\!0\}\,.
\end{align}
Here $T_*W_+^s$ is of codimension 1 in $\mathrm{Fix}(S_{J_+^c})$.

In the case $R(\boTh_+)\!=\!1$ of a totally synchronous 1-cluster target $\boTh_+\!=\!0$, the results hold with $\alpha_+^b\!=\!1,\ J_3\!=\!J_+^c\!=\!\emptyset$.
Moreover, $T_* W_+^s = \mathrm{Fix}(S_{J_+^c})=\boXi$, with vanishing codimension.
\end{lem}

\begin{proof}[\textbf{Proof}]
We prove the claims, in order of appearance.
The trivial adaptations to the synchronous 1-cluster case $\boTh_+\!=\!0$ will be omitted.

By section \ref{2Clu}, both equilibria $\boTh_\pm\!\in\!\boXi/\Lambda'$ are 2-clusters.
Claim (i) then follows from \eqref{Wu} and \eqref{Ws} in lemma \ref{lemW}, applied to $\boTh\!=\!\boTh_\pm$\,.

To prove the strict monotonicity claim (ii), we first note $\alpha_-^b\!<\!\alpha_+^b$, by strict increase of the order parameter $0\!<\!R\!<\!1$ along $\bolth_*(t)$\,.
This mere increase in size already excludes the homoclinic case $\boTh_+\!=\!\boTh_-$ and the case $J_+^b\!=\!J_-^b$\,.
It remains to prove the weak monotonicity claim $J_+^b\supseteq J_-^b$\,.
We apply \eqref{Wu} with $\boTh:=\boTh_-$ to $\bolth_*(t)\!\in\! W_-^u$\,, in the closed time-invariant subspace
\begin{equation}
\label{ThJ}
\boTh_+=\lim_{t\rightarrow\infty}\,\bolth_*(t)\in \mathrm{Fix}(S_{J_-^b})\,.
\end{equation}
In particular, the components $\Theta_{+,j}$ with $j\!\in\! J_-^b$ all coincide.
Since the cluster $J_-^b$ is fat, already, this implies weak monotonicity $J_-^b\subseteq J_+^b$ for the fat cluster $J_+^b$ of $\boTh_+$\,.

This proves claim (ii) of the lemma, and strict monotonicity of fat clusters, theorem \ref{thmmon}(iii).
It also justifies the $M\!=\!3$ partition $\mathbf{J}$ by $J_1:=J_-^b\,,\ J_2:=J_+^b\setminus J_-^b\,,\ J_3:=J_+^c$\,.
We proceed via the associated 3-partitioned coordinates $(x_2,x_3,\boxi_1,\boxi_2,\boxi_3)$ of \eqref{xxi}. 

Claim \eqref{xi130} follows by definition of the 3-partitioned components $\boxi_1, \boxi_3$ associated to $J_1\!=\!J_-^b$ and $J_3\!=\!J_+^c$\,.
Indeed, \eqref{Wu} and \eqref{Ws} apply to $\bolth_*(t)\in W_-^u\cap W_+^s$\,.
With \eqref{FS'm} we obtain
\begin{equation}
\label{thSpm}
\bolth_*(t)\, \in\, \mathrm{Fix}(S_{J_-^b})\cap\mathrm{Fix}(S_{J_+^c})\, = \, 
\mathrm{Fix(S_{J_2}')}\, =\, 
\mathrm{span}\{\mathbf{1}_2, \mathbf{1}_3\}\oplus\boXi_2\, =\, 
\{\boxi_1\!=\!0,\ \boxi_3\!=\!0\}\,.
\end{equation}

For the remaining two claims \eqref{T*Wu} and \eqref{T*Ws} in covering space $\boXi$, we recall equivariance \eqref{TWg} of tangent spaces under $\sigma\!=\!\pm\mathbf{id}_N$ and $\bopsi\!\in\!\Lambda'$.
This allows us to choose representatives $\boTh_\pm\!\in\!\boXi/\Lambda'$ with the above 2-clustering, separately, and with unchanged tangent spaces.
Then the two claims follow from \eqref{TWu} and \eqref{TWs}, similarly and respectively, in 3-cluster coordinates $(x_2,x_3,\boxi_1,\boxi_2,\boxi_3)$.
Here we use \eqref{FSm} to prove the last equality in \eqref{T*Wu}.
To prove claim \eqref{T*Ws}, we start from \eqref{TWs} to proceed in the same way:
\begin{equation}
\label{T*Ws+}
T_*W_+^s\, \leq\, \mathrm{Fix}(S_{J_3}) \, =\, 
\mathrm{span}\{\mathbf{1}_2,\mathbf{1}_3\}\oplus\boXi_1\oplus\boXi_2 \, =\, 
\{\boxi_3\!=\!0\}\,.
\end{equation}
Recall how the subspace inclusion $\leq$ is of codimension 1.
This proves the lemma, and theorem \ref{thmmon}(iii).
\end{proof}

\section{Three-cluster orbits}\label{3Clu}

In this section, we prove theorem \ref{thm3clu} on 3-cluster heteroclinic orbits, as an elementary building brick of all heteroclinicity.
This generalizes section \ref{2Clu} on 2-cluster dynamics.

Consider any pair $\boTh_\pm\!\in\!\bT^{N\!-\!1}=\boXi/\Lambda'$  of 2-cluster equilibria with order parameters $0<R_-=2\alpha_-^b\!-\!1<2\alpha_+^b\!-\!1=R_+<1$, fat size fractions $\alpha_\pm^b$, and associated $M\!=\!3$ index partition $\mathbf{J}$ given by $J_1:=J_-^b\,,\ J_2:=J_+^b\setminus J_-^b \,,\ J_3:=J_+^c$\,; see theorem \ref{thmmon}(iii) and \eqref{J123}.
Theorem \ref{thm3clu} requires the construction of a 3-cluster heteroclinic orbit $\bolth_*(t){:}\, \ \boTh_-\!\leadsto\! \boTh_+$ in the torus $\bT^{N\!-\!1}\!=\!\boXi/\Lambda'$ .
We therefore start from any such 3-cluster partition with size fractions $0<\alpha_2, \alpha_3<1/2<\alpha_1<1$.
The proof of theorem \ref{thm3clu} will then follow from a detailed analysis of the flow of all 3-cluster solutions $\bolth(t)$.
In the 3-cluster coordinates $(x_1,x_2,x_3)$ of \eqref{cluODE} with $M\!=\!3$, we have to discuss
\begin{equation}
\label{3cluODE}
   \dot x_m\,=\,\alpha_{m-1} \sin(x_{m-1}-x_m)+\alpha_{m+1} \sin(x_{m+1}-x_m)\,,
\end{equation}
for indices $m\ \mathrm{mod}\,3$.
This is the dynamics on $\mathbf{X_J}=\mathrm{Fix}(S_\mathbf{J})=\bR^3$\,; see \eqref{XJ}, \eqref{FixSJ}.

\begin{figure}[t]
\centering \includegraphics[width=0.9\textwidth]{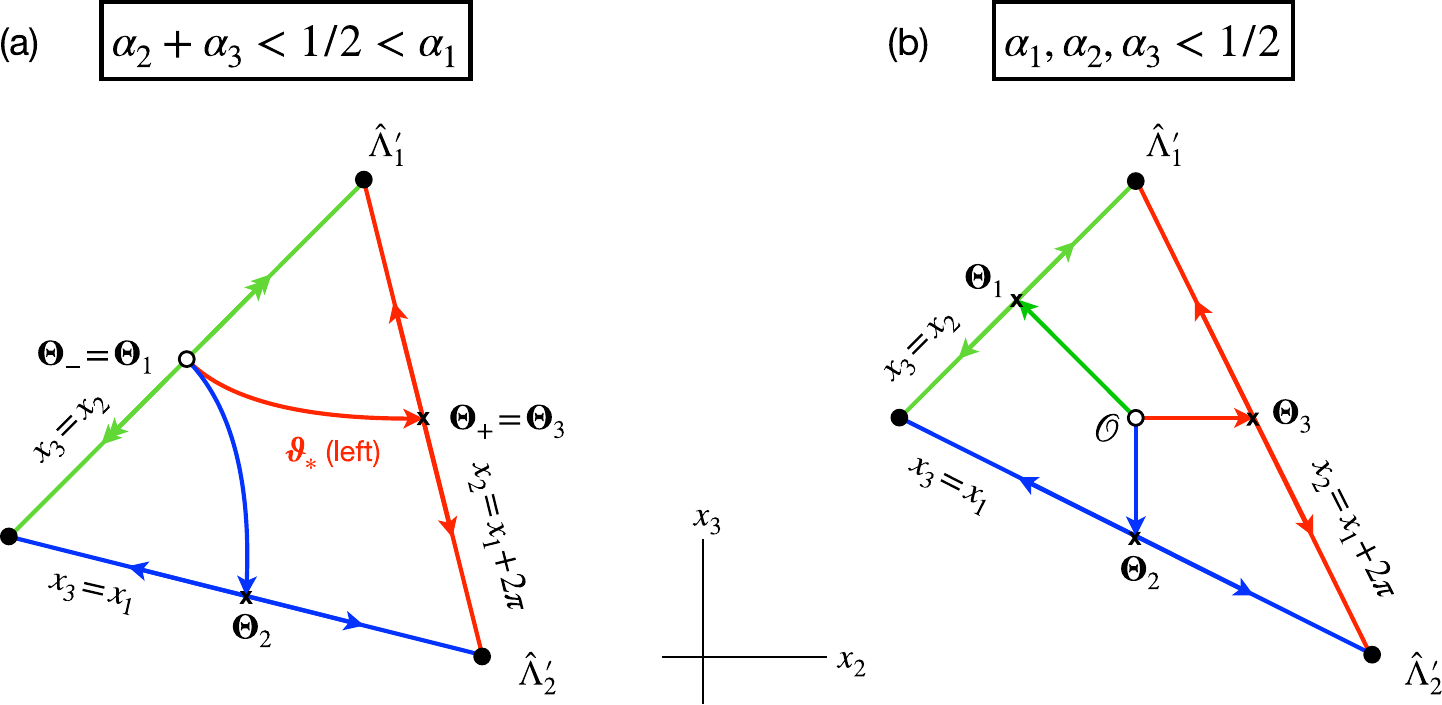}
\caption[Fundamental triangles of 3-cluster dynamics]{\emph{
Qualitative dynamics on fundamental triangles \eqref{triangle} of the 3-cluster flows \eqref{3cluODE} under the phase constraint \eqref{S1}, in the coordinate plane $(x_2, x_3)$ of $\boXi_\mathbf{J}$\,.
See \eqref{TN-1}--\eqref{XiJ} and \eqref{1+Xim}, \eqref{xxi} for notation.
The vertices $\boTh\!=\!0$ (solid dots) are torus copies of the 1-cluster equilibrium $x_1\!=\!x_2\!=\!x_3\!=\!0$.
The three sides are the time-invariant lines $x_{m-1}\!=\!x_{m+1}$\,, with indices $m\ \mathrm{mod}\ 3$; a torus copy is used for $m\!=\!3$.
Each side contains a unique nontrivial 2-cluster equilibrium $\boTh_m$\,, at its midpoint.
Compare section \ref{2Clu} and proposition \ref{Equil}.\\
\textbf{Case (a):} $1/2<\alpha_1<1$.\quad  
For illustration, we have used size fractions $(\alpha_1,\alpha_2,\alpha_3)=(3,1,1)/5$.
The 2-cluster equilibrium $\boTh_-\!=\!\boTh_1$ (circle) possesses a fat cluster of size fraction $1/2\!<\!\alpha_-^b\!=\!\alpha_1\!=\!3/5\!<\!1$.
It is unstable, with fast unstable eigenvalue $\mu_1 \!=\! +1$ along its side (green) and slow unstable eigenvalue $0<\mu_2=2\alpha_1\!-\!1<1$ towards the interior; see \eqref{mu2}. 
The equilibria $\boTh_2$ and $\boTh_3\!=\!\boTh_+$\,, marked by $\times$, possess fat clusters of larger size fractions $\alpha_1\!+\!\alpha_3$ and $\alpha_+^b=\alpha_1\!+\!\alpha_2$\,, respectively.
By \eqref{mu2}, again, they are hyperbolic saddles with $\mu_2\!<\!0\!<\!\mu_1=+1$.
The stable and unstable manifolds of the saddles are colored blue for $\boTh_2$\,, and red for $\boTh_3$\,.
Their stable manifolds establish heteroclinic orbits $\boTh_1\leadsto\boTh_2$ and the left rebellion $\bolth_*{:}\,\ \boTh_1\!\leadsto\!\boTh_3$\,.
Note their vertical tangents $x^c\!=\!x_2\!=\!0$, at $\boTh_2$\,, and horizontal tangents $x^c\!=\!x_3\!=\!0$, at $\boTh_3$\,.
\\
\textbf{Case (b):} $0<\alpha_1, \alpha_2, \alpha_3<1/2$.\quad
For schematic illustration, we have used size fractions $(\alpha_1,\alpha_2,\alpha_3)=(1,1,1)/3$.
All three 2-cluster equilibria $\boTh_m$ are hyperbolic saddles $\times$,  now, with their sides as unstable manifolds; see \eqref{mu2}.
Their time-invariant manifolds are colored green, blue, and red, respectively.
The interior of the fundamental triangle features the unique 3-bar linkage equilibrium $\mathcal{O}$ with order parameter $R\!=\!0$, marked by $\circ$.
Here, $\mathcal{O}$ is the equilateral triangle according to \eqref{R=0}.
The increasing Lyapunov function $R$ implies total instability of $\mathcal{O}$.
Hence, the green, red, and blue stable manifolds of the 2-cluster saddles $\boTh_m$ emanate from the interior 3-cluster equilibrium.\\
}}
\label{fig51}
\end{figure}

The phase constraint \eqref{S1} on the time-invariant linear section $\boXi\!=\!\mathbf{1}_N^\bot$ eliminates uniform phase shifts, putting $x_1:=-\alpha_1^{-1}(\alpha_2x_2\!+\!\alpha_3x_3)$\,.
In \eqref{XiJ}--\eqref{TJM-1}, now for $M\!=\!3$, we have already seen how the dynamics \eqref{3cluODE} on the covering space $\boXi_\mathbf{J}=\mathbf{X_J}\cap\boXi=\bR^2$ reduces to the 2-torus $\bT_\mathbf{J}^2=\boXi_\mathbf{J}/\Lambda_\mathbf{J}'$.
The restriction $\mathbf{P_J}$ of the projection $\mathbf{P}$, to $\mathbf{X_J}$\,, defines the discrete restricted lattice
\begin{equation}
\label{LJ'}
\Lambda_\mathbf{J}'\ :=\ \mathbf{P_J}\Lambda_\mathbf{J}\ \leq\ \boXi\cap\mathbf{X_J}\ =\ \bR^2\,,
\end{equation}
analogously to \eqref{P}--\eqref{TN-1}.
In coordinates $(x_1,x_2,x_3)$:
\begin{equation}
\label{PJx}
\mathbf{P_J}=
\begin{pmatrix}
  1-\alpha_1    & -\alpha_2  & -\alpha_3  \\
   -\alpha_1   & 1-\alpha_2  & -\alpha_3  \\
   -\alpha_1   & -\alpha_2  &  1-\alpha_3 
\end{pmatrix}\,.
\end{equation}
The first two columns $\Lambda_1', \Lambda_2'$ already define the 2-torus $\bT_\mathbf{J}^2=\boXi_\mathbf{J}/\Lambda_\mathbf{J}'$\,. 
Explicitly, $\Lambda_\mathbf{J}'=\mathbf{P_J}2\pi\bZ^3=\bZ\Lambda_1'\oplus\bZ\Lambda_2'$ becomes the free $\bZ$-module generated by
\begin{equation}
\label{LJ'x}
\Lambda_1'=
-\begin{pmatrix}
      2\pi        \\
      0            \\
      0  
\end{pmatrix} +
2\pi\alpha_1\begin{pmatrix}
    1       \\
    1       \\
    1    
\end{pmatrix},\qquad
\Lambda_2'=
\begin{pmatrix}
      0        \\
      2\pi    \\
      0  
\end{pmatrix} -
2\pi\alpha_2 \begin{pmatrix}
    1       \\
    1       \\
    1    
\end{pmatrix}.
\end{equation}
For the planar coordinates $(x_2,x_3)$ of $\boXi_\mathbf{J}$ in \eqref{xxi}, we obtain generators 
$\hat{\Lambda}_\iota'$ by just omitting the first components of $\Lambda_\iota'$ in \eqref{LJ'x}.
See figure \ref{fig51} for two examples.

Equivariance \eqref{groupS} under $\sigma\!=\!-\mathbf{id}_N\!\in\!\mathbf{S}$ allows us to describe the dynamics in the (closed) \emph{fundamental triangle} of figure \ref{fig51}:
\begin{equation}
\label{triangle}
\{\tau_1\hat{\Lambda}_1'+\tau_2\hat{\Lambda}_2'\,\vert\, 0\leq\tau_1,\tau_2,\ \tau_1+\tau_2\leq1\}\,.
\end{equation}
In redundant cluster phases $(x_1,x_2,x_3)$ of $\boXi_\mathbf{J}$\,, the fundamental triangle reads
\begin{equation}
\label{xtriangle}
x_1\leq x_3\leq x_2\leq x_1\!+\!2\pi\,.
\end{equation}

The synchronous 1-cluster equilibrium copies arise as $x_1\!=\!x_2\!=\!x_3\!=\!0 \ \mathrm{mod}\, \Lambda_\mathbf{J}'$.
Any 2-cluster equilibria $\boTh_m$ must arise, in this system \eqref{3cluODE}, via $x_{m-1}\!=\!x_{m+1} \ \mathrm{mod}\, \Lambda_\mathbf{J}'$\,, with indices $m\ \mathrm{mod}\,3$.
In our fundamental triangle, we have to replace the line $x_2\!=\!x_1$ by $x_2\!=\!x_1\!+\!2\pi$, to meaningfully satisfy \eqref{xtriangle}.
Our description of 2-cluster dynamics, in section \ref{2Clu}, applies to each triangle side, with properly adapted notation.
In particular, each side of the fundamental triangle is time-invariant.
Lifting to the covering space $\boXi_\mathbf{J}$\,, via $\Lambda_\mathbf{J}'$\,, three lines intersect where any two do, because any such intersection represents a copy of the totally synchronous 1-cluster equilibrium $\boTh\!=\!0\!\in\!\bT^2$.
From \eqref{2cluODE}, \eqref{2line}, and proposition \ref{Equil} we recall that the nontrivial 2-cluster equilibria $\boTh_m$ are located at the midpoints of each side of the fundamental triangle \eqref{triangle}. 

In variables $(x_2,x_3)$, figure \ref{fig51}(a) sketches the case where the cluster $J_1\!=\!J_-^b$ is fat, i.e.\ of size fraction $\alpha_1\!=\!\alpha_-^b\!>\!1/2$.
Then the other two clusters remain slim, even when combined as $J_2\cup J_3$\,, i.e.
\begin{equation}
\label{alpha3a}
0<\alpha_2+\alpha_3=1-\alpha_1<1/2<\alpha_1<1\,.
\end{equation}

The only 3-cluster equilibria, which are not degenerate 2- or 1-clusters in disguise, 
arise as distinct 3-bar linkages with $R\!=\!0$, i.e.
\begin{equation}
\label{R30}
0=\alpha_1\exp(\mi x_1)+ \alpha_2\exp(\mi x_2)+ \alpha_3\exp(\mi x_3)\,;
\end{equation}
see \eqref{RPsi}, \eqref{R=0}. Such triangles with distinct angles $x_1,x_2,x_3$ arise if, and only if
\begin{equation}
\label{alpha3b}
0<\alpha_1,\alpha_2, \alpha_3<1/2\,.
\end{equation}
Indeed all other cases lead to linearly degenerate triangles, or violate the triangle inequality, altogether, because $\alpha_1\!+\!\alpha_2\!+\!\alpha_3\!=\!1$.
In figure \ref{fig51}(b) we sketch case \eqref{alpha3b}, again in variables $(x_2,x_3)$.

For linearization of \eqref{3cluODE}, we invoke block diagonalization \eqref{Ablock} in $\boXi$.
In cluster notation for $\bT^2$, properly adapted to each $\boTh_m$\,, we recover eigenvalues $\mu_1\!=\!a_0\!=\!+1$ along each side of the fundamental triangle.
By \eqref{am2}, the eigenvalue towards the interior is
\begin{equation}
\label{mu2}
    \mu_2=a_b=-2(\alpha_{m-1}+\alpha_{m+1})+1=2\alpha_m\!-\!1, \quad\textrm{at } \boTh_m\,,
\end{equation}
for $m\!=\!1,2,3$.
Of course, this also follows from \eqref{3cluODE}, directly.
In the 3-bar linkage triangle case \eqref{alpha3b} this implies a stable eigenvalue $\mu_2\!<\!0$  towards the interior, and hyperbolicity, at all boundary saddles $\boTh_m$\,; see figure \ref{fig51}(b).
The unique interior rational 3-bar linkage chain equilibrium $\mathcal{O}$ is totally unstable, by proposition \ref{Equil}(iii).

In the opposite case $\alpha_1 \!>\! 1/2>\alpha_2\!+\!\alpha_3$ of figure \ref{fig51}(a), with a fat majority cluster $J_1$\,, we see that only the boundary 2-cluster equilibrium $\boTh_1$ at $x_2\!=\!x_3$ is two-dimensionally unstable.
The other two boundary 2-cluster equilibria $\boTh_2$ and $\boTh_3$ remain hyperbolic saddles, as before. 
Since $0<\mu_2=2\alpha_1\!-\!1<1$, the 2-cluster side $0<x_2=x_3<2\pi\alpha_1$\,, with eigenvalue $\mu_1\!=\!1$ at the midpoint $\boTh_1$, constitutes its fast unstable manifold. 
All other unstable orbits emanate from the stationary 2-cluster tangentially to the eigenvector of the slower unstable eigenvalue $\mu_2$\,, as part of its ``mustache'' \cite{ArnoldODE}.
We have omitted the case $\alpha_1=1/2=\alpha_2\!+\!\alpha_3$ here, because it generates a degenerate 2-cluster equilibrium $\boTh_1$ at $R\!=\!0$, with a zero eigenvalue $\mu_2\!=\!0$ at the boundary $x_2\!=\!x_3$.
In the absence of other equilibria or periodic orbits, this establishes the planar flows of figure \ref{fig51} in the fundamental triangle \eqref{triangle} of the 3-cluster flow on the constrained 2-torus $\bT^2$; see for example \cite{HSD}.

\emph{\textbf{Proof of theorem \ref{thm3clu}}}.\quad
Given any two equilibria $\boTh_\pm\!\in\!\bT^{N\!-\!1}\!=\!\boXi/\Lambda'$ with order parameters $0\!<\!R_-\!<\!R_+\!\leq\! 1$, we have to construct a 3-cluster heteroclinic orbit $\bolth_*{:}\, \ \boTh_-\!\leadsto\! \boTh_+$ between them.
By theorem \ref{thmmon}, the source $\boTh_-$ and the target $\boTh_+$ are 2-cluster equilibria, with ascending fat clusters on $J_-^b\subsetneq J_+^b$ and descending complementary slim clusters on $J_-^c\supsetneq J_+^c$\,.
Therefore, \eqref{J123} uniquely determines an $M\!=\!3$ index partition $J_1:=J_-^b\,,\ J_2:=J_+^b\setminus J_-^b\,,\ J_3:=J_+^c$ of our 3-cluster candidate $\bolth_*(t)$; see theorem \ref{thmmon}.
This also determines the respective size fractions $\alpha_1=\alpha_-^b\,,\ \alpha_2=\alpha_+^b-\alpha_-^b\,,\ \alpha_3=1\!-\!\alpha_+^b$.
In particular, our assumption $1/2<\alpha_\pm^b=(R_\pm+1)/2<1$ implies $1/2\!<\!\alpha_1\!=\!\alpha_-^b\!<\!1$.

Our analysis of the 3-cluster flow \eqref{3cluODE} therefore leads to a 3-cluster heteroclinic orbit $\bolth_*(t){:}\, \ \boTh_1 \!\leadsto\! \boTh_3$\,.
In the fundamental triangle of figure \ref{fig51}(a), that heteroclinic orbit $\bolth_*(t)$ is located in the red stable manifold of the saddle $\boTh_3$\,.
Moreover, our choice of $J_1,J_2,J_3$ implies $\boTh_-\!=\!\boTh_1$ and $\boTh_+\!=\!\boTh_3$\,, up to translation by $\Lambda'$.
This proves that the 3-cluster orbit $\bolth_*(t){:}\, \ \boTh_-\!\leadsto\! \boTh_+$\,, in $\bT_\mathbf{J}^2=\boXi_\mathbf{J}/\Lambda_\mathbf{J}'$\,, is indeed heteroclinic between the prescribed source and target equilibria $\boTh_\pm$\,.

It remains to prove claim \eqref{x2pm}, i.e.\ to determine the limiting values $x_m(t)$ at $t\!=\!\pm\infty$, i.e.\ at $\boTh_\pm$ in covering space $\boXi_\mathbf{J}$\,.
We invoke \eqref{LJ'x} and find
\begin{equation}
\label{Th3}
\boTh_-=\boTh_1=\Lambda_1'/2= 
\pi\begin{pmatrix}
 \alpha_1-1        \\
 \alpha_1        \\
 \alpha_1      
\end{pmatrix}, \quad
\boTh_+=\boTh_3=(\Lambda_1'+\Lambda_2') /2=
\pi\begin{pmatrix}
 \alpha_1-\alpha_2-1       \\
 \alpha_1-\alpha_2+1       \\
\alpha_1-\alpha_2       
\end{pmatrix}.
\end{equation}
The bulk phases of $\boTh_\pm$ are $x_\pm^b:=x_1(\pm\infty)$.
Then \eqref{Th3} implies the limits $x_3(\pm\infty)=x_\pm^b\!+\!\pi$.
Also $x_2(-\infty)=x_3(-\infty)=x_-^b\!+\!\pi$.
At $t\!=\!+\infty$, however, we obtain $x_2(+\infty)=x_+^b\!+\!2\pi$.
This identifies our specific construction of $\bolth_*(t){:}\, \ \boTh_-\!=\!\boTh_1\!\leadsto\! \boTh_3\!=\!\boTh_+$ as a left rebellion, in the sense of \eqref{x2pm}, \eqref{left}.

\begin{figure}[t]
\centering 
\includegraphics[width=0.5\textwidth]{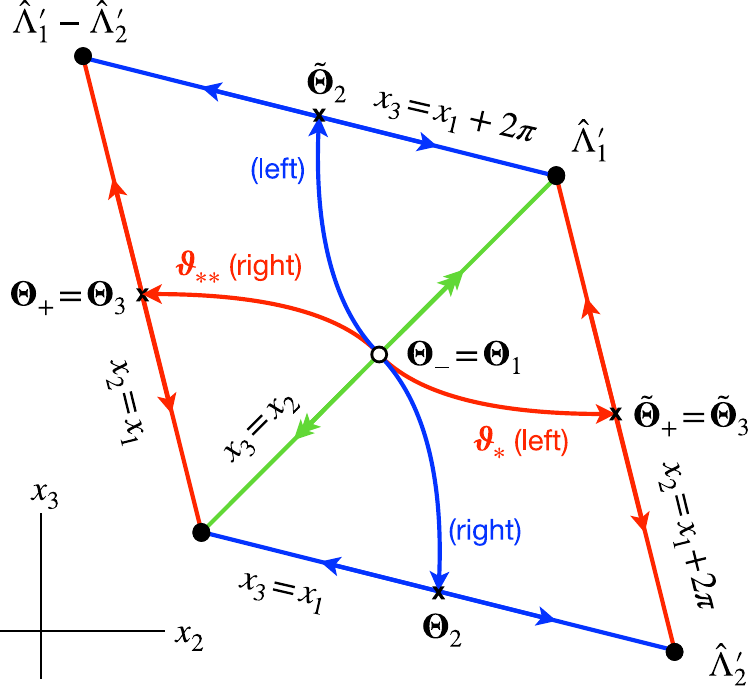}
\caption[Fundamental parallelogram of 3-cluster dynamics]{\emph{
We consider the fat/slim case of figure \ref{fig51}(a), in the same notation and colors.
The fundamental parallelogram \eqref{parall} of the 3-cluster dynamics consists of two adjacent copies of the fundamental triangle \eqref{triangle} from figure \ref{fig51}(a).
We chose to attach the second copy to the upper left of the fundamental triangle.
The two red branches $\bolth_*(t)$ and $\bolth_{**}(t)$ of the stable manifold $W^s(\boTh_3)$ arrive at their lattice-equivalent target saddles $\tilde{\boTh}_+\!=\!\tilde{\boTh}_3$ and $\boTh_+\!=\!\boTh_3$ horizontally, but from opposite sides in $x_2$\,.
They originate from the same source $\boTh_-\!=\!\boTh_1$\,.
In theorem \ref{thm3clu}, this generates the two options \eqref{right} and \eqref{left} for the right rebellion $\bolth_{**}$ and the left rebellion $\bolth_*$\,.
}}
\label{fig52}
\end{figure}

It remains to construct a second heteroclinic orbit $\bolth_{**}(t){:}\, \ \boTh_-\!=\!\boTh_1\!\leadsto\! \boTh_3\!=\!\boTh_+$ which constitutes a right rebellion \eqref{right} between instantiations of the same  equilibria $\boTh_\pm \textrm{ mod } \Lambda_\mathbf{J}'$\,.
We work in the following \emph{fundamental parallelogram} of the planar lattice $\Lambda_\mathbf{J}'$\,:
\begin{equation}
\label{parall}
\{\tau_1\hat{\Lambda}_1'+\tau_2\hat{\Lambda}_2'\,\vert\, 0\leq\tau_1,\tau_2<1\}\,;
\end{equation}
see figure \ref{fig52} and remark \ref{resetTheta}.
Here we have attached a $180^\circ$ rotated copy
\begin{equation}
\label{xright}
x_1\leq x_2\leq x_3\leq x_1+2\pi
\end{equation}
of the fundamental triangle \eqref{triangle}, \eqref{xtriangle} of figure \ref{fig51}(a), along its left side $x_2\!=\!x_3$ (green).
More precisely, note
\begin{equation}
\label{xparall}
x_1\leq x_2, x_3 < x_1+2\pi
\end{equation}
in the parallelogram \eqref{parall}; compare \eqref{right}, \eqref{left}.
We have relabeled the target $\boTh_3$ of the left rebellion $\bolth_*(t)$ from figure \ref{fig51}(a) as $\tilde{\boTh}_3$\,, in figure \ref{fig52}, because $x_2=x_1\!+\!2\pi$ excludes it from the fundamental parallelogram \eqref{xparall}.
See also remark \ref{resetTheta}.

We use equivariance of the Kuramoto model \eqref{kurODE} under $\sigma\!=\!-\mathbf{id}_N$ and apply the rotation $\bolth_{**}:=-\bolth_*+\Lambda_1'$\,; see \eqref{groupS}, \eqref{groupT}.
This automorphism rotates the fundamental parallelogram of figure \ref{fig52} by $180^\circ$ around the fixed source $\boTh_-\!=\!\boTh_1\!=\!\Lambda_1'/2$; see \eqref{Th3}.
It also maps the target $\tilde{\boTh}_+\!=\!\tilde{\boTh}_3\!=\!(\Lambda_1'\!+\!\Lambda_2')/2$ to its translated copy $\boTh_+\!=\!\boTh_3\!=\!(\Lambda_1'\!-\!\Lambda_2')/2\!=\!\tilde{\boTh}_+\!-\!\Lambda_2'$\,.
Since the $180^\circ$ rotation swaps right and left rebellions, $\bolth_{**}{:}\, \ \boTh_-\!\leadsto\! \boTh_+$ is a right rebellion in the fundamental parallelogram \eqref{parall} of $\bT_\mathbf{J}^2=\boXi_\mathbf{J}/\Lambda_\mathbf{J}'$ between the prescribed equilibria $\boTh_\pm$\,.
In figure \ref{fig52}, the two rebellions constitute  the two branches of the stable manifold $W_+^s$ (red) of $\boTh_+=\tilde{\boTh}_+\ \mathrm{mod}\ \Lambda_\mathbf{J}'$\,.
The two branches arrive at their target equilibrium $\boTh_+\!=\!\tilde{\boTh}_+$ along the same horizontal tangent, but from opposite sides in $x_2$. 
Indeed \eqref{TWsTh} asserts horizontal tangents in coordinates $(x_2, x_3)$.
Likewise, the two branches emanate from the unstable source equilibrium $\boTh_-\!=\!\boTh_1$ and $\tilde{\boTh}_-$ along the same slow unstable eigenvector tangent of $0<\mu_2=2\alpha_1\!-\!1< 1=\mu_1$\,, but in opposite directions.
This distinguishes, and realizes, both options \eqref{left} and \eqref{right} for 3-cluster heteroclinic orbits $\bolth_*,\,\bolth_{**}{:}\, \ \boTh_-\!\leadsto\! \boTh_+$\,. 

The remaining 2-cluster heteroclinic orbit $\boTh_1\!=\!\boTh_-\!\leadsto\!\boTh_+\!=\!0$ towards total synchrony follows from invariance of the boundary line $x_2\!=\!x_3$.
This completes the proof of theorem \ref{thm3clu}.   
\qed

\emph{\textbf{Proof of corollary \ref{corxshift}}}.\quad
We address the left rebellion $s\!=\!+$, first.
In the proof of theorem \ref{thm3clu}, the left rebellion is represented by the heteroclinic orbit $\bolth_*{:}\, \ \boTh_-\!=\!\boTh_1\!\leadsto\!\tilde{\boTh}_3\!=\!\tilde{\boTh}_+$\,.
See the right fundamental triangle in figure \ref{fig52}.
Then \eqref{Th3} implies
\begin{equation}
\label{xshiftl}
\tilde{\boTh}_+-\boTh_-=\tilde{\boTh}_3-\boTh_1=(\Lambda_1'+\Lambda_2'-\Lambda_1')/2=\Lambda_2'/2\,.
\end{equation}
The first component $-\pi\alpha_2$ of $\Lambda_2'$ in \eqref{LJ'x} proves claim \eqref{xshift} for $s\!=\!+$.

For right rebellions $s\!=\!-$ by $\bolth_{**}{:}\,\ \boTh_-\!\!\leadsto\!\!\boTh_+$\,, we obtain
\begin{equation}
\label{xshiftr}
\boTh_+-\boTh_-=-\Lambda_2'/2\,,
\end{equation}
which now proves claim \eqref{xshift} for $s\!=\!-$. 
This proves corollary \ref{corxshift}.
\qed

\section{Transversality}\label{Trans}

In this section, we prove transversality theorem \ref{thmtrv}, and theorem \ref{thmMclu} on  $M$-cluster rebellions.

Let $\bolth_*(t){:}\, \ \boTh_-\!\leadsto\! \boTh_+$ denote any heteroclinic orbit in covering space $\boXi$ between two equilibria $\boTh_\pm$ with order parameters $0\!<\!R_-\!<\!R_+\!<\!1$.
As in the proof of lemma \ref{lemW}, equivariance \eqref{Wg}, \eqref{TWg} under $g=(\mathbf{id}_N,\bopsi),\ \bopsi\!\in\!\Lambda'$, implies $\Lambda'$-invariance of tangent spaces, in $\boXi$, at the same $\bolth\in\boXi/\Lambda'$.
This allows us to work in the setting of section \ref{2Clu}, for $\boTh_\pm\in\boXi/\Lambda'$.
By cluster monotonicity of theorem \ref{thmmon}(iii), we may therefore invoke the partition $\mathbf{J}=(J_1,J_2,J_3)$ of \eqref{J123} and the 3-partitioned coordinates
\begin{equation}
\label{J123c}
(x_2,x_3,\boxi_1,\boxi_2,\boxi_3)\ \in\ \bR^2\oplus\boXi_1\oplus\boXi_2\oplus\boXi_3
\end{equation}
associated to $\boTh_1:=\boTh_-$ and $\boTh_3:=\boTh_+$\,; see \eqref{xxi}. 
Here $\bR^2:=\bR\mathbf{1}_2\oplus\bR\mathbf{1}_3$\,.
We denote the unstable manifold of $\boTh_-$ and the stable manifold of $\boTh_+$ by $W_-^u$ and $W_+^s$\,, respectively.

\emph{\textbf{Proof of theorem} \ref{thmtrv}}. \quad
In case of total synchrony $\boTh_+\!=\!0,\ R_+\!=\!1$, linear stability $i(\boTh_+)\!=\!0$ leaves nothing to be proved.

In the 2-cluster cases $0\!<\!R_-\!<\!R_+\!<\!1$, we have to show the left equality of
\begin{equation}
\label{trv13}
T_{\,\bolth_*(t)}W_-^u\, +\, T_{\,\bolth_*(t)}W_+^s\ =\ \boXi\ =\ 
\bR^2\oplus\boXi_1\oplus\boXi_2\oplus\boXi_3\,;
\end{equation}
compare \eqref{trv}.
The claim is invariant under small perturbations of the tangent spaces.
Since the stable manifold $W_+^s$ is $C^1$, and $\bolth_*(t)\rightarrow \boTh_+$ for $t\rightarrow\infty$, we may equivalently show
\begin{equation}
\label{trv13infty}
T_{\,\bolth_*(t)}W_-^u\, +\, T_{\,\boTh_+}W_+^s\ =\ 
\bR^2\oplus\boXi_1\oplus\boXi_2\oplus\boXi_3\,,
\end{equation}
for some large $t \gg 1$.

For any $t\!\in\!\bR$, lemma \ref{lemWpm} with $\bolth:=\bolth_*(t)$ asserts
\begin{equation}
\label{TW1u*}
T_{\bolth_*(t)}W_-^u=\mathrm{Fix}(S_{J_1})=\bR^2\oplus\boXi_2\oplus\boXi_3\,;
\end{equation}
see \eqref{T*Wu}.
At the 2-cluster target equilibrium $\boTh\!=\!\boTh_3\!=\!\boTh_+$\,, \eqref{TWsTh} of lemma \ref{lemW} asserts 
\begin{equation}
\label{TW3s+}
T_{\,\boTh_+}W_+^s\ =\ 
\mathrm{Fix}(S_{J_+^c})\cap\{x_+^c\!=\!0\}\ =\ \boXi_+^b\ =\ 
\bR\mathbf{1}_2\oplus\boXi_1\oplus\boXi_2\,.
\end{equation}
Indeed, $J_+^c\!=\!J_3$ and $\mathrm{Fix}(J_3)=\bR^2\oplus\boXi_1\oplus\boXi_2$\,, in the 3-partitioned coordinates \eqref{xxi}.
See \eqref{FSm}.
Combined, \eqref{TW1u*} and \eqref{TW3s+} prove transversality claim \eqref{trv13infty}, and theorem \ref{thmtrv}.   
\qed

\emph{\textbf{Proof of theorem} \ref{thmMclu}}. \quad
We can write any heteroclinic orbit $\bolth_*(t){:}\, \ \boTh_-\!\leadsto\! \boTh_+$ as an $M$-cluster orbit with coordinates $x_m$\,; see \eqref{J}--\eqref{cluODE}.
By permutation equivariance \eqref{sigma}, distinct $\vth_j(t)$ cannot cross in finite time.
We may therefore choose $M$ minimal so that all $x_m$ are, and remain, distinct for all finite times $t$.
By theorem \ref{thmmon}, we have fat equilibrium cluster sets $J_1:=J_-^b\subsetneq J_+^b$ and slim $J_M\!=\!J_+^c$\,.
Let $\mathbf{J}$ denote this $M$-clustering partition.
Then $\bolth_*(t)\in \bT_\mathbf{J}^{M\!-\!1}=\boXi_\mathbf{J}/\Lambda_\mathbf{J}'$\,; see \eqref{TJM-1}.

Note $M\!\geq\! 3$, because $R_+\!<\!1$ excludes synchronous targets $\boTh_+\!=\!0,\ R_+\!=\!1$.
Moreover, $M\!-\!2\leq |J_+^b\setminus J_-^b|=N_+\!-\!N_-$\,, because each of the $M\!-\!2$ remaining clusters contains at least one element.

Lemma \ref{lemWpm}(i) implies claim \eqref{hetFix} of the theorem.
Minimality of $M$, invariance of phase ordering under time, and application of suitable translations in $\Lambda'$, imply the ordering \eqref{xorder}, possibly after relabeling of indices $m$ by a permutation $\sigma\!\in\! S_M$.
See section \ref{kurGroup}.

The ``converse'' claim of the theorem requires the construction of clustered heteroclinic orbits $\bolth_*(t){:}\, \ \boTh_-\!\leadsto\! \boTh_+$ in $\bT_\mathbf{J}^{M\!-\!1}=\boXi_\mathbf{J}/\Lambda_\mathbf{J}'$\,, with prescribed order \eqref{xorder} of the cluster coordinates, for any $M\!\geq\!3$, any clustering $\mathbf{J}\!=\!(J_m)$\,, and any prescribed $1\!\leq\! m_*\!<\! M$.
For $M\!=\!3$, the 3-cluster theorem \ref{thm3clu} proves just that: see \eqref{right} for the right rebellion $m_*\!=\!2$, and \eqref{left} for the left rebellion $m_*\!=\!1$.

To address $M\!\geq\!4$, we first construct a heteroclinic orbit $\bolth_*(t){:}\, \ \boTh_-\!\leadsto\!\boTh_+$ in  $\bT_\mathbf{J}^{M\!-\!1}$, via theorem \ref{thm3clu} and concatenation corollaries \ref{corconcat}, \ref{shift-type}.
In a second step, we show that $\bolth_*(t)$ realizes the prescribed ordering \eqref{xorder} of phases $x_m(t)$.

To construct $\bolth_*(t)$\,, let $n:=N_+^b\!-\!N_-^b$\,.
As in \eqref{sequence}, prescribe any sequence $\mathbf{s}$ of symbols $s_j\!=\!\pm$.
Then corollary \ref{shift-type} provides a cascading sequence $\bolth_j(t)\,{:}\,\ \boTh_{j\!-\!1}\!\leadsto\!\boTh_j$ of one-man rebellions by $j=N_-^b\!+\!\ell,\ 1\!\leq\!\ell\!\leq\! n$, between equilibria $\boTh_-\!=\!\boTh_{N_-^b}\,,\ldots,\boTh_{N_+^b}\!=\!\boTh_+$\,, with ascending fat clusters of size $j$ at $\boTh_j$\,.
To pass to prescribed clusters $J_m$\,, as in \eqref{xm}, we now choose identical symbols $s_j\!=\!\mathfrak{s}_m$\,, for $\vth_j\!=\!x_m$ in the same cluster $j\!\in\! J_m\,,\ 1\!<\!m\!<\!M$.
Specifically, we choose right rebellions $\mathfrak{s}_m\!=\!-$, \eqref{right}, for the left interval $2\!\leq\! m \!\leq\! m_*$ of $m$, and left rebellions $\mathfrak{s}_m\!=\!+$, \eqref{left}, for the right interval $m_*\!<\!m\!<\!M$.
Let $j_m$ denote the last $j\!\in\! J_m$\,, and let $\boTh_m^*:=\boTh_{j_m}$\,.
For each $1\!<\!m\!<\!M$, consider the path of heteroclinic $\bolth_j(t)\,, \ j\in J_m$\,, in the connection graph $\cC$.
By our comments on corollary \ref{corconngraph}, each such path can be represented by an express 3-cluster heteroclinic $\bolth_m^*(t){:}\, \ \boTh_{m-1}^*\!\leadsto\!\boTh_m^*$\,, with associated right/left sign $\mathfrak{s}_m$\,.
Here we define $\boTh_1^*:=\boTh_-$\,; also note $\boTh_M^*:=\boTh_+$\,.
By corollary \ref{corconcat}, restricted to the invariant subspace $\mathrm{Fix}(\mathbf{J})$ of the prescribed partition $\mathbf{J}\!=\!(J_m)$, we may concatenate the sequence of heteroclinic 3-clusters $\boTh_m^*$ to a direct $\eps$-shadow $\bolth_*(t){:}\,\ \boTh_-\!\leadsto\!\boTh_+$\,.

It remains to show that $\bolth_*(t)$ realizes the prescribed ordering \eqref{xorder} of phases, for sufficiently small $\eps\!>\!0$.
All along the first heteroclinic 3-cluster $\bolth_2^*(t)$, with $\mathfrak{s}_2\!=\!-$ and right rebel cluster $J_2$, theorem \ref{thm3clu} and \eqref{right} imply $x_1<x_2<x_3=\ldots=x_M<x_1\!+\!2\pi$.
Of course we assume the presence of $m\!=\!2\!\leq\! m_*$, here.
Let us iterate to any rebellion $\bolth_m^*(t)$ with rebel cluster $J_m$, for any $1\!<\! m\!\leq\! m_*\!<\!M$, i.e.\ with all $\mathfrak{s}_m\!=\!-$.
All along each $\bolth_m^*(t)$, allegiance switching \eqref{x2pm} of the rebel phases $x_m$ then successively implies
\begin{equation}
\label{xorderright}
x_1=\ldots=x_{m-1}<x_m<x_{m+1}=\ldots=x_M<x_1+2\pi\,.
\end{equation}
For $m_*\!<\!m \!<\! M$, analogously, successive left rebellions by $J_m$ with $\mathfrak{s}_m\!=\!+$ imply
\begin{equation}
\label{xorderleft}
x_1=\ldots=x_{m_*}<x_{m+1}=\ldots=x_M<x_m<x_{m_*\!+\!1}=\ldots=x_{m-1}=x_1+2\pi\,, 
\end{equation}
all along $\bolth_m^*(t)$.
(For $m\!=\!m_*\!+\!1$, the rightmost string of equalities still collapses to $x_1\!+\!2\pi$, only.)
For small $\eps\!>\!0$, this implies the phase ordering \eqref{xorder}, all along the $\eps$-shadow $\bolth_*(t)$ of the concatenated heteroclinic 3-clusters $\bolth_m^*(t)$.
Indeed, we only have to choose $\eps$ small enough, such that each strict inequality in \eqref{xorder} remains satisfied, even after perturbation by $\eps$, at \emph{some} time $t\!=\!t_m$ along \emph{each} heteroclinic 3-cluster step $\bolth_m^*(t)$.
Strict inequality is then inherited all along the concatenation $\bolth_*$\,.
This completes the proof of theorem \ref{thmMclu}.
\qed

Invariance of the pair $\boTh_\pm$ under all permutations of indices $j\in\hat{\mathbf{J}}=J_2\dot\cup\ldots\dot\cup J_{M\!-\!1}$\,, of course, generates $\bolth_*(t){:}\, \ \boTh_-\!\leadsto\! \boTh_+$ in $\bT^{N\!-\!1}\!=\!\boXi/\Lambda'$ for arbitrary rearrangements of the clustering phases $x_m\,,\ 1\!<\!m\!<\!M$ in \eqref{xorder}.
Moreover, $1\!\leq\! M\!-\!2\!\leq\! N_+^b\!-\!N_-^b$ and the clustering partition of $\hat{\mathbf{J}}$ can be chosen arbitrarily.
In \eqref{part21eps} of the next section, we will numerically illustrate another construction of $\bolth_*(t)$\,, which avoids $\eps$-shadowing.

\section{Numerical illustrations}\label{Numerics}

In this section we comment on, and illustrate, concatenations of transverse heteroclinic rebellions among 2-clusters and towards the 1-cluster of total synchrony.
See concatenation corollaries \ref{corconcat}, \ref{shift-type} and figures \ref{fig21}--\ref{fig23}, for successive one-man rebellions.
We also illustrate the general swarm rebellions of theorem \ref{thmMclu} and figure \ref{fig25}.

\subsection{Concatenation of one-man rebellions}\label{NumConcat}

We address one-man rebellions first.
These are special cases $0\!<\!\alpha_3\!<\!1/2\!<\!\alpha_1\!<\!1,\ \alpha_2\!=\!1/N$, of the 3-cluster rebellions studied in theorem \ref{thm3clu} and section \ref{3Clu}.

Prior to their concatenation, we sketch how to approximate heteroclinic saddle-saddle connections.
The numerical literature has established numerical methods and codes for such tasks; see for example \cite{Beyn, KuzCode, DoedelSurvey}.
For the Kuramoto model \eqref{kurODE}, however, the 3-cluster phase plane of figure \ref{fig52} offers a much more elementary (and cheaper) approach.
Indeed, the concatenation corollaries \ref{corconcat}, \ref{shift-type} only require us to approximate each 3-cluster step $\bolth_\ell(t){:}\, \ \boTh_{\ell-1}\!\leadsto\!\boTh_\ell$ in our heteroclinic cascade \eqref{trans}, individually, by a $\delta$ pseudo-orbit.
To construct a concatenation $\bolth_*(t){:}\, \ \boTh_-\!=\!\boTh_0\!\leadsto\!\boTh_n\!=\!\boTh_+$\,, we then invoke uniform $\eps$-shadowing.
Moreover, total instability of $\boTh_-$ within the 3-cluster plane of figure \ref{fig52} suggests tracing trajectories \emph{backwards in time}.

To approximate any heteroclinic 3-cluster step $\bolth_\ell(t)$\,, in its 3-cluster coordinates $(x_1,x_2,x_3)$ of section \ref{3Clu}, we start from a phase perturbation $\boldsymbol{\delta} = (\delta_1, \delta_2, 0)$ of the target $\boTh_\ell$\,, such that $\delta_2=-s_\ell 10^{-2}$.
Here $s_\ell \!\in\! \{-, +\}$ denotes the symbols \eqref{spm}, \eqref{sequence} of right/left rebellions in corollary \ref{shift-type}. 
The vanishing third component of $\boldsymbol{\delta}$ ensures that the perturbation is horizontally tangent to the stable manifold $W^s(\boTh_\ell)$\,, at the target $\boTh_\ell$\,; see \eqref{TWsTh}. 
To perturb $\boTh_\ell$ within the flow-invariant space $\boXi\!=\!\mathbf{1}_N^\bot$ of phase constraint \eqref{S1}, we require
\begin{equation}
\label{numalpha}
    \alpha_1\delta_1 + \alpha_2\delta_2  = 0\,.
\end{equation}
For our choice of $\delta_2=-s_\ell10^{-2}$\,, this determines $\delta_1$.
In particular, the 3-cluster solution $\tilde{\bolth}_\ell(t)$ of \eqref{3cluODE} which starts from
\begin{equation}
\label{numpert}
    \tilde{\bolth}_\ell(0):=\boTh_\ell+ \boldsymbol{\delta}\,,
\end{equation}
in the fundamental parallelogram \eqref{parall} of figure \ref{fig52}, approximates $\bolth_\ell(t)$ and preserves the phase constraint $\boXi$. 
Here we have to choose target representatives $\boTh_\ell\!=\!\boTh_3$\,, for $s_\ell\!=\!-$ with right rebellions \eqref{right}, and $\boTh_\ell\!=\!\tilde{\boTh}_3$\,, for $s_\ell\!=\!+$ with left rebellions \eqref{left}.

For $\ell=1,\ldots,n$, we then perform the following loop.
\begin{description}
    \item[Step 1:]  Solve the 3-cluster system \eqref{3cluODE} with initial condition \eqref{numpert} numerically \emph{in backward time}, until reaching a \emph{stopping time} $t \!=\! T_\ell\!<\!0$ such that $x_2(T_\ell) - x_3(T_\ell) \ \mathrm{mod}\, 2\pi$ is smaller than $10^{-2}$, in absolute value.  
    By our rigorous analysis of the dynamics on the fundamental triangles in figure \ref{fig51}, we conclude that $0\!>\!T\!>\!-\infty$ is finite and we are numerically close to a copy of the source 2-cluster equilibrium $\boTh_{\ell-1}$\,.    
    \item[Step 2:] Decrease $\ell$ by 1.
\end{description}

By this loop, we obtain numerical approximations $\tilde{\bolth}_\ell(t)$ of the rebellious sequences $\bolth_\ell(t)$ shown in corollary \ref{corconcat}, backwards in time and in reversed order $\ell$. 
\textbf{Step 2} jumps from a backwards solution $\tilde{\bolth}_\ell(T)\in W_\mathrm{loc}^u(\boTh_{\ell-1})$ to a new initial condition $\tilde{\bolth}_{\ell-1}(0)$ defined by \eqref{numpert} for index $\ell\!-\!1$, instead of index $\ell$.
I.e.\ we solve \eqref{3cluODE}, with updated size fraction parameters $\alpha_m$ and in the updated coordinates $(x_1,x_2,x_3)$ of a different fundamental domain, adapted to the previous heteroclinic step $\bolth_{\ell-1}(t)$. 
The theoretical justification for this jump is the shadowing property which we have introduced just before concatenation corollary \ref{corconcat}.
To produce continuous $\eps$-shadows, however, we do not reset phase angles according to remark \ref{resetTheta}.
Rather, we accumulate phase shifts, as in \eqref{sumxshift} of concatenation corollary \ref{shift-type}.

How large to choose $\delta_2$\,, such that the $\delta$ pseudo-orbits $\tilde{\bolth}_\ell(t)$ actually produce a concatenated $\eps$-shadow $\bolth_*(t)$, within image resolution of figures \ref{fig22}, \ref{fig23}?
In fact, we chose totally synchronous $\boTh_+$ as a target, with $n\!=\!i(\boTh_-)\!=\!|J_-^c|\!=\!10$ steps of heteroclinic one-man rebellions.
By inspection, we found our fixed choice of modestly small $|\delta_2|$ and thresholds, at $10^{-2}$, up to the task, all the way.

\subsection{Cluster swarm rebellions}\label{NumSwarm}

In figure \ref{fig25} we have illustrated theorem \ref{thmMclu} on $M$-cluster swarm rebellions $\bolth_*(t){:}\, \  \boTh_-\!\leadsto\! \boTh_+$ from a 2-cluster equilibrium $\boTh_-$ at $t \to -\infty$ towards a 2-cluster equilibrium $\boTh_+$ with a larger fat cluster. 
We have partitioned the $N\!=\!21$ oscillators of $\bolth_*(t)$ into $M\!=\!11$ clusters of sizes $N_1\!=\!N_-^b\!=\!11,\ N_2\!=\!\ldots\!=\!N_M\!=\!1$, and $N_+^b\!=\!20$.
Morse indices $i_\pm:=i(\boTh_\pm)\!=\!|J_\pm^c|$ are therefore $i_-\!=\!10$ and $i_+\!=\!1$.
Specifically, the cluster decomposition \eqref{J} reads $J_-^b\!=\!J_1:=\{1,\ldots,11\}$, and $J_m:=\{10\!+\!m\}$ for $2\!\leq\! m\!\leq\! 11$.
Note $J_M\!=\!J_+^c\!=\!\{21\}$ and $\hat{\mathbf{J}}=J_2\dot\cup\ldots\dot\cup J_{M\!-\!1}=\{12,\ldots,20\}$.
We aim for a total of 9 rebels $m\!=\!2,\ldots,10$.
Their positions $x_m$ should realize the prescribed order \eqref{xorder}, with $m_*\!-\!1\!=\!4$ right rebels $s_m\!=\!-$, at $m\!=\!2,\ldots,5$, versus $M\!-\!1\!-\!m_*$ left rebels $s_m\!=\!+$, at $m\!=\!6,\ldots,10$.
Consider partitioned coordinates
\begin{equation}
\label{part21}
(\hat{x},x_M,\boxi_1,\hat{\boxi})
\end{equation}
with respect to the index partition $J_1\dot\cup\hat{\mathbf{J}}\dot\cup J_M$; compare \eqref{xxi}.
We omit $\boxi_M\!=\!0$, because the one-man cluster $J_M\!=\!\{21\}$ is already represented by its ``average'' $x_M$\,.
In coordinates \eqref{part21}:
\begin{equation}
\label{part21th}
\begin{aligned}
   \bolth_*(t)\ &\in\ \mathrm{Fix}(S_{J_-^b})\cap\mathrm{Fix}(S_{J_+^c})\ =\ \{\boxi_1\!=\!0\} \,;   \\
   \boTh_-\ & =\ (11/21,11/21,0,0)\pi\,;\\
   \boTh_+\ & =\ (\hat{x}_+\,,\,10/21\,\pi,\,0,\,\hat{\boxi}_+)\,.
\end{aligned}
\end{equation}
At $\boTh_-$\,, our coordinates have omitted $x_1\!=\!x_-^b\!=\!-\alpha_-^c\pi\!=\!-10/21\,\pi$; see \eqref{midpoint}.
At $\boTh_+$ we obtain $x_1\!=\!-11/21\,\pi$, by cumulative $x_1$-shift \eqref{sumxshift} in concatenation corollary \ref{shift-type}.
The nine rebel components $x_m\!=\!\hat{x}_+\!+\!\hat{\boxi}_{+,m}$ of $\boTh_+$ satisfy $x_m\!=\!x_1$\,, for the four right rebels $m\!=\!2,\ldots,5$, and $x_m\!=\!x_1\!+\!2\pi$, for the five left rebels $m\!=\!6,\ldots,10$.
This determines $\hat{x}_+$ and $\hat{\boxi}_+$\,.

As an approximate heteroclinic $\tilde{\bolth}(t)$, from $\boTh_-$ to $\boTh_+$ in the local stable manifold $W_\mathrm{loc}^s(\boTh_+)$, we pick the following initial perturbation of $\boTh_+$ with $\delta\!=\!10^{-2}$\,:
\begin{equation}
\label{part21eps}
\tilde{\bolth}(0)\ :=\ \boTh_++\delta\cdot(\hat{x},0,0,\hat{\boxi})\ \in\  \mathrm{Fix}(S_{J_-^b})\cap\mathrm{Fix}(S_{J_+^c})\,;
\end{equation}
see lemma \ref{lemWpm}(i), \eqref{TWsTh}, and \eqref{T*Ws}.
We then solve \eqref{kurODE} backwards in time, where \eqref{TWu} and \eqref{Wu} ensure $\tilde{\bolth}(0)\!\in\! W_-^u$ and hence convergence of $\tilde{\bolth}(t)$ to $\boTh_-$\,, for $t\rightarrow-\infty$.
 
Indeed, lemma \ref{lemW} asserts that the unstable manifold $W_-^u$ of $\boTh_-$ is open in, and of the same dimension as, $\mathrm{Fix}(S_{J_-^b})$.
To guarantee that the closure of $W_-^u$ extends to $\boTh_+$\,, at all, we also use existence of a heteroclinic orbit $\tilde{\bolth}_*(t){:}\, \  \boTh_-\!\leadsto\! \boTh_+$ with clusters ordered by \eqref{xorder}, e.g.\ as constructed in theorem \ref{thmMclu} and section \ref{NumConcat} by concatenation.
Alternatively, we may invoke the density argument in \eqref{Wud} for $d\!=\!10$ and the torus $\bT^{10}\!=\!\bT_\mathbf{J}^{M\!-\!1}$\,, in our notation \eqref{TJM-1}.

In case the whole rebellion is required to swarm out towards the same side of $x_M$\,, unilaterally right or left, we start from a perturbation $\tilde{\bolth}(0)$ with nonzero average component $\hat{x}$\,, of the appropriate sign, and then a perturbation $|\hat{\boxi}|\ll |\hat{x}|$ of zero average to fix the appropriate order \eqref{xorder}, all along the swarm $\tilde{\bolth}(t),\ T\!\leq\! t\!\leq\!0$.
To realize general prescribed sign configurations \eqref{xorder} of $\hat{\boxi}$, we let $\hat{x}\!=\!0$.
We then choose the first $m_*\!-\!1\!=\!4$ nonzero components $m\!=\!2,\ldots,5$ of the perturbation $\hat{\boxi}$ on $\hat{\mathbf{J}}$ negative, and the remaining $M\!-\!1\!-\!m_*\!=\!5$ nonzero components $m\!=\!6,\ldots,10$ positive.
The zero average condition on $\hat{\boxi}\!\in\!\hat{\boXi}$ is easily respected here.
As with the 3-cluster rebellions $\tilde{\bolth}_\ell(t)$ for figures \ref{fig22}, \ref{fig23} before, we may consider the orbit $\tilde{\bolth}(t)$ of \eqref{part21eps} as a $\delta$ pseudo-orbit, in just one step, with an associated heteroclinic $\eps$-shadow $\bolth_*(t)$ of $\tilde{\bolth}(t)$.
Note that $\bolth_*(t)$ inherits the order \eqref{xorder} from $\tilde{\bolth}(t)$.
Depending on the choice of $\hat{\boxi}$, however, the resulting heteroclinic $\eps$-shadow swarm $\bolth_*(t){:}\, \  \boTh_-\!\leadsto\! \boTh_+$  may deviate widely from the original concatenation $\tilde{\bolth}_*(t)$ of 3-clusters.
See figure \ref{fig25}.

Simulation code has been implemented in \texttt{Julia} \cite{Julia}. 
Plots have been generated with the package \texttt{Plots.jl}. 
We integrated via the standard explicit 4th-order Runge--Kutta method in the package \texttt{DifferentialEquations.jl} with fixed time steps of $10^{-2}$. 
We have made the code more human-readable using the free version of \texttt{Cursor}, as an AI aid to arrange documentations.

\section{Discussion}\label{Disc}

To put our results in perspective, we comment on robustness under small perturbations, including changes of the model type.
See section \ref{Pert}.
Section \ref{SL} compares our global results on the Kuramoto model \eqref{kurODE} with previous local results on the Stuart--Landau variant.
The final section \ref{REV} is a note added in revision, to discuss a very elegant and closely related preprint which has appeared after submission of the first version of our present paper.

\subsection{Rebellions and perturbations} \label{Pert}

We have considered the classical Kuramoto model \eqref{omegaj} with equal frequencies $\omega_j \!=\! \omega$ and order parameter $R \!\geq\! R_0 \!>\! 0$ in \eqref{RPsi}\,:
\begin{equation}
\label{omega}
\dot\vth_j\,=\,\omega+\tfrac{1}{N}\sum_{k=1}^N\,\sin(\vth_k-\vth_j)\,, 
\end{equation}
for $j\!=\!1,\ldots,N$.
Our results address an uncharted, yet fundamental, problem regarding global dynamics.
Co-rotating coordinates reduce \eqref{omega} to the case $\omega\!=\!0$ of \eqref{kurODE}.
We have then explored paths of rebellions from a plethora of partially synchronized equilibria to total synchrony. 
The consequences of our approach and results are at least threefold.

\emph{First}, some rebellions to total synchrony involve concatenations of successive rebellions among unstable, or ``metastable'' 2-cluster configurations.
This is not accessible to the Watanabe--Strogatz transformation, which excludes any heteroclinic orbits between 2-cluster equilibria. 
Our approach is complementary, and is based on the presence of a fat majority cluster. 

\emph{Second}, for $R \!>\! 0$ and $\omega\!\neq\!0$,  corollary \ref{strucstab} has established structural stability of the original Morse--Smale system \eqref{omega} on the global attractor \eqref{cAlift}, in the full angle space $\bT^N$. 
This yields $C^0$ orbit equivalent flows under \textit{all} $C^1$-small perturbations of the Kuramoto model \eqref{omega}.
In particular, we may consider small deviations \eqref{omegaj} from the single frequency regime $\omega_j \!\equiv\! \omega$. 
Any permutation-equivariance from $S_N$ is then lost.
Periodic orbits survive, due to hyperbolicity, but lose their antipodal 2-cluster structure $\Th_j\!\in\!\{\Psi,\Psi\!+\!\pi\}$; see \eqref{Lyap} and section \ref{Equi}.
In particular, total synchrony of the single-frequency model (also called \textit{phase synchronization}) prevails, but not as a cluster of identical $\Th_j(t)$.

Suppose next that the perturbation preserves shift-equivariance \eqref{groupT} under $\bR\mathbf{1}_N$, only, but all equivariance \eqref{groupS} under any permutations or $\sigma\!=\!-\mathbf{id}_N$ from \eqref{groupS} may have been destroyed.
Then, at least, the weaker notion of \textit{frequency synchronization} prevails:
\begin{equation}
\label{fre-syn}
\lim_{t \to \infty}\  \sup_{i, j = 1,\dots,N} |\dot \vth_i(t) - \dot\vth_j(t)| = 0,
\end{equation} 
for any solutions with $R\!\geq\! R_0\!>\!0$. 
This follows from in-phase convergence towards periodic orbits, because all periodic orbits are normally hyperbolic and remain relative equilibria (here: rigidly rotating waves), under small perturbations.
See \eqref{releq}, and \cite{FieHopf} for further details on rotating waves.

Much more generally than the original model \eqref{omega}, consider perturbations like
\begin{equation}
\label{generalODE}
m_j \ddot\vth_j +  \dot\vth_j\,=\,\omega_j+\tfrac{K}{N}\sum_{k=1}^N\,\sin(\vth_k(t\!-\! \tau_{jk})
-\vth_j(t \!-\! \tilde{\tau}_{jk}) - \gamma_{jk})\,;
\end{equation}
see \cite{DoeBu, HKL, HKP,  HsJuKw, HJKU}. 
Assume that frequency deviations $|\omega_j\!-\!\omega|$, inertiae $m_j \!\geq\! 0$, time delays $\tau_{jk}, \tilde{\tau}_{jk} \geq 0$, and phase lags $|\gamma_{jk}|$ are all small.
The resulting $N$-dimensional slow inertial manifold, in the sense of Fenichel theory \cite{Jones}, is then invariant and strongly exponentially attracting.
The slow global dynamics, on it, therefore amounts to a $C^1$-small perturbation of the structurally stable single-frequency model \eqref{omega}, again for order parameters $R\!\geq\! R_0\!>\!0$.
Notably, positive time delays increase the dimensionality of the phase space to become infinite \cite{halu93}. 
Nevertheless, the above approach admits small delays as well; see \cite{Ch}. 
Early results on the elimination of small delays in nonlinear equations go back as far as \cite{Kurz}.
Morse--Smale structural stability then admits small perturbations of the general form \eqref{generalODE}, by corollary \ref{strucstab}.

\emph{Third}, our approach is based on all-to-all $S_N$-equivariance and the existence of a Lyapunov function $R$.
Transversality follows from the somewhat surprising observation that invariant manifolds of cluster equilibria are contained, globally, in linear subspaces characterized by invariance under certain subgroups of $S_N$, and share the same dimension with such linear fix-spaces. 
Many other variants of the Kuramoto model also retain at least some permutation-equivariance, and a Lyapunov function. 
Examples are more general network topologies defined by weighted coupling constants $K_{ij} \!\geq\! 0$ (see \cite{ChHsHs}), or by larger inertial terms.
We have not yet ventured into these vast areas, systematically.

\subsection{Comparison with Stuart--Landau oscillators} \label{SL}
 
For comparison, we recall some previous results on \emph{local} $S_N$ symmetry-breaking among $N$ nonlinear complex Stuart--Landau oscillators with linear all-to-all coupling.
These results have been derived and proved for the dynamics of 3-clusters, in the limit of large $N$.
For complete details, see \cite{Krischer1, Krischer2} and the references there.

The Stuart--Landau system reads
\begin{equation}
\label{slj}
\dot{W}_j=(1-(1+\mathrm{i}\gamma)|W_j|^2)W_j+ \beta \tfrac{1}{N} \sum_{k=1}^N (W_k - W_j)\, .
\end{equation}
Here $W_j \!\in\! \mathbb{C}$ indicates phase and amplitude of the $j$-th oscillator, $j\!=\!1, \ldots, N$.
The papers consider real amplitude dependence $\gamma\!\neq\!0$ of individual periods,  and complex coupling $\beta \!\in\! \mathbb{C}$. 
Note equivariance of \eqref{slj} under the action of $\sigma\!\in\! \mathbf{S}$\,, analogously to \eqref{groupS}.
Translation equivariance \eqref{groupT} by $\bopsi\!\in\!\mathbf{T}$ appears in the phases of $W_j$\,.

For a background and motivation, we recall how \eqref{slj} often serves, in physics, as a ``normal form'' for oscillatory systems close to the onset of oscillation and under the influence of a linear coupling through the mean field \cite{Kurbook,VGM2008}. 
This normal form has been established to be a good approximation in a multitude of contexts from various disciplines, whether it be in physics, chemistry, biology, neuroscience, social dynamics, or engineering. 
For an overview see e.g.~\cite{Pik2003, PiRoSurv} or references 1--15 in \cite{KuGiOtt2015}.

\begin{figure}[t]
\centering
\includegraphics[width=\textwidth]{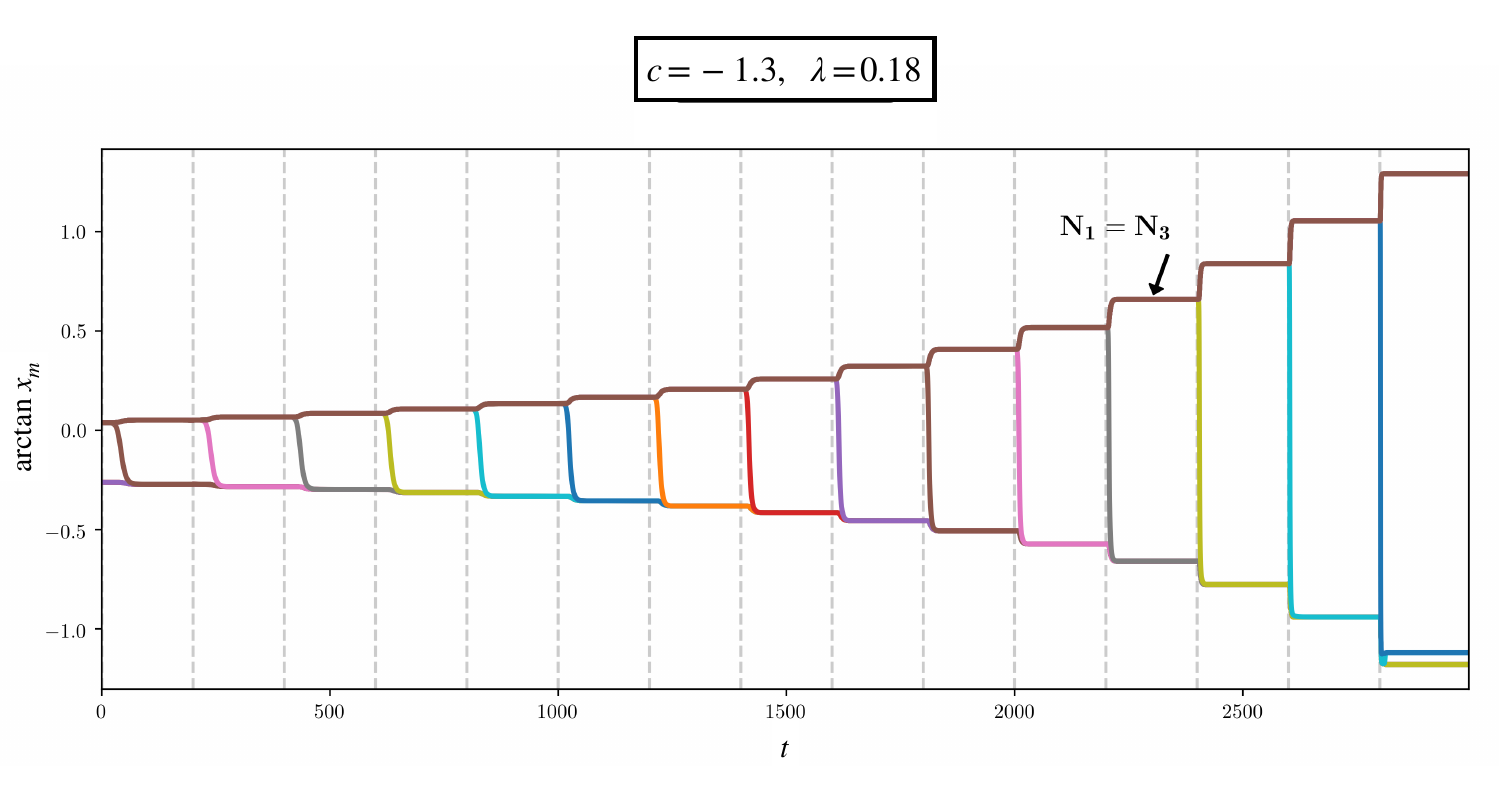}
\caption[Concatenation of Stuart--Landau rebellions]{\emph{
Heteroclinic rebellions obtained from numerical simulations of \eqref{ODEyj} for $c\!=\!-1.3$, $\lambda\!=\!0.18$ and $N\!=\!32$ oscillators. 
See \emph{\cite[section 5]{Krischer2}} for details.
The graphs indicate 3-cluster solutions \eqref{3cluy} with one-man rebellions of $x_2\!=\!y_n\,,\ n\!=\!N_1\!+\!1$.
Starting from $N_1\!=\!4$, a small random perturbation is applied to $y_n\!=\!x_3$\,, at a time indicated by the dashed vertical lines, only.
Then ODE \eqref{ODEyj} with rebel  $x_2=y_n\neq x_1\,,x_3\,$ is integrated forward, until the dynamics settles again.
See the colored rebellions of $x_2\!=\!y_n$ from $x_3$ (top) to $x_1$ (bottom), along which the rebel $x_2\!=\!y_n$ switches its cluster allegiance.
The minority cluster takes over the majority, after $N_1\!=\!N_2$\,.
Updating $N_1\mapsto N_1\!+\!1$ and $n\mapsto n\!+\!1$, this process is repeated until instability to a 3-cluster equilibrium takes over from $N_1\!=\!18$.}
}
\label{fig8}
\end{figure}

The papers \cite{Krischer1, Krischer2} address bifurcation from the totally synchronous 1-cluster solution $W_j(t):=\exp(-\mathrm{i}\gamma t)$.
The Benjamin--Feir instability describes destabilization of synchrony, due to the appearance of a Floquet exponent zero of high (nontrivial) multiplicity $N-1$.
Explicit center manifold reduction, in suitable coordinates, and truncation to polynomial order 3, lead to a reduced real system of the general form
\begin{equation}
	\label{ODEyj}
	\dot{y}_j = f_j(\mathbf{y}) := \left(\lambda+c\cdot \langle y^2\rangle\right) y_j+ \widetilde{y^2_j}+ \widetilde{y^3_j}
\end{equation}
for $\mathbf{y}\in\mathbf{Y}:=\bR^N$. 
Here $\langle y^p \rangle$ and $\widetilde{x^p_j}:= x^p_j\!-\! \langle x^p\rangle$ abbreviate the average of $p$-th powers and the deviation from that average, respectively.
The real parameter $\lambda$ is the linearization at the trivial equilibrium $\mathbf{y}\!=\!0$ of total synchrony, which corresponds to the reference 1-cluster synchronous periodic orbit $W$. 
The Benjamin--Feir instability of $\mathbf{y}\!=\!0$ occurs at $\lambda\!=\!0$.
In fact, the zero sum space $\mathbf{Y}_0:=\{\langle y\rangle\!=\!0\}$ is invariant, not only under the linear action \eqref{groupS} of the group $\mathbf{S}$, but also under the nonlinear dynamics of \eqref{ODEyj}.
Indeed $\langle y \rangle= \langle \widetilde{y_j^p} \rangle=0$ on $\mathbf{Y}_0$\,.
See \cite{Elm1,DiSt,StElm,Elm2,GoStNet} for analysis of equilibrium bifurcations.

The importance of \eqref{ODEyj} reaches far beyond the Stuart--Landau context, or any direct interpretation as a network of $N$ identical scalar ``cells'' with all-to-all coupling via power sums. 
For example,  \cite{Elm1, StElm} arrive at \eqref{ODEyj} in the context of biological evolution of sympatric speciation.
In fact, steady-state bifurcation analysis of \emph{any} fully permutation-symmetric network, at eigenvalue 0, typically leads to irreducible eigenspaces.
Beyond trivial total synchrony $y_1\!=\!\ldots\!=\!y_N$\,, the standard representation on $\mathbf{Y}_0$ provides the simplest interesting case.
Up to scaling and possible time-reversal, in fact, \eqref{ODEyj} on $\mathbf{Y}_0$ constitutes the most general cubic vector field which is equivariant under the standard representation \eqref{g} of $S_N$ on $\mathbf{Y}_0$.
Just based on linearization at total synchrony, rather than the computational intricacies of center manifold reduction, it is therefore not surprising that all-to-all Stuart--Landau coupling \eqref{slj} reduces and truncates to \eqref{ODEyj}.
But conversely, the explicit computations in \cite{Krischer2} also show the Stuart--Landau setting \eqref{slj} to be universal, in the sense that all real higher order coefficients $c$, and all abstract bifurcation parameters $\lambda$ near $\lambda\!=\!0$, actually arise in \eqref{ODEyj}, for suitable Stuart--Landau parameters $\gamma\!\in\!\bR$ and $\beta\!\in\!\mathbb{C}$ in \eqref{slj}.

Surprisingly, the cubic truncation \eqref{ODEyj} possesses an explicit \emph{gradient structure}.
In fact, \emph{all} $S_N\,$-equivariant polynomial vector fields $\mathbf{f}$ on the standard representation $\mathbf{Y}_0 \,$, up to and including order three, are gradients of scalar invariants.
For higher orders, in contrast, that gradient property fails.

The gradient property allowed \cite{Krischer1,Krischer2} to study the dynamics of \eqref{ODEyj}, in terms of clustering and rebellions among cluster equilibria, much in the spirit of our present approach to the Kuramoto model \eqref{kurODE}.
Globally, \cite{Krischer1} classifies the cluster dynamics of \eqref{ODEyj} by seven different nondegenerate cases, depending on the location of
\begin{equation}
\label{cintervals}
c\ \in\ \bR\setminus\{-2,-3/2,-4/3,-5/4,-1,-1/2\}\,.
\end{equation}
For $N\!=\!32$ oscillators, figure \ref{fig8} illustrates the choice $-4/3\!<\!c\!=\!-1.3\!<\!-5/4$ at $\lambda\!=\!0.18$\,. 
In the spirit of figures \ref{fig22}, \ref{fig23}, the 3-cluster dynamics of successive one-man rebellions is documented by their heteroclinic transients between two large clusters near equilibrium.
In cluster notation analogous to \eqref{xm}--\eqref{S1}, the three clusters are
\begin{equation}
\label{3cluy}
\begin{aligned}
   x_1&:=y_1\!=\! \dots \!=\! y_{N_1}\,, \\  
   x_2&:=y_{n}\,,\\
   x_3&:=y_{n\!+\!1}\!=\! \dots \!=\! y_N\,,
   \end{aligned}
\end{equation}
with the abbreviation $n\!=\!N_1\!+\!1$.
Start from $N_1\!=\!4$ and $x_2\!=\!x_3$\,.
This corresponds to a \emph{small} initial size fraction $\alpha_1\!=\!N_1/N \!=\! 0.125$.
After a small random perturbation of the single component $y_n$ in cluster $x_2\,$, the clustered ODE \eqref{ODEyj} is integrated forward in time, until the dynamics appear to settle.
Heteroclinic rebellion ensues, i.e.\ the perturbed oscillator $y_n\!=\!x_2$ changes its cluster affiliation from the large cluster at $x_3$ to the small cluster at $x_1$.
In other words, $N_1\!=\!5$, after the rebellion.
Note that the heteroclinic rebellion, from majority cluster to minority cluster, occurs in the opposite direction of the Kuramoto model \eqref{kurODE}; compare section \ref{Main}.
Updating $N_1\mapsto N_1\!+\!1$ and $n\mapsto n\!+\!1$, rebellions towards the minority cluster repeat, for ever increasing cluster sizes $N_1\,$.

After 12 rebellions, equal cluster parity $N_1\!=\!N_2\!=\!16$ is reached and the former minority cluster $x_1$ reaches and overtakes the majority.
After two further rebellions, the dynamics enter a region with 3-cluster equilibria.
At that stage, the third coexisting cluster of \eqref{ODEyj}, at $x_2\!<\!x_3$ near $x_1\!<\!x_2$\,, consists of just one single rebel.
Such 3-cluster equilibria cannot occur in the present Kuramoto model \eqref{kurODE}.
Similarly, cubic \eqref{ODEyj} does not admit $M$-cluster equilibria with $M\!\geq\!4$. 
Multi-cluster \emph{transients} in the Stuart--Landau case \eqref{slj} have not been analyzed, to our knowledge.

\subsection{Note added in revision}\label{REV}
After posting the first version of our present paper \textbf{P1}\footnote{J.-Y.\ Dai, B.\ Fiedler, and A.\ López-Nieto. Transient rebellions in the Kuramoto oscillator: Morse--Smale structural stability and connection graphs of finite 2-shift type. \url{https://arxiv.org/abs/2512.02937v1}}
on Dec-02-2025, and submission to JEMS, a closely related paper \textbf{P2}\footnote{D.\ Burns, G.\ Malajovich, C.\ Pugh, I.\ Rajapakse, and S.\ Smale. Global analysis of the Kuramoto flow. \url{https://arxiv.org/abs/2605.22983v1}} 
appeared on May-21-2026.
Both preprints address the co-rotating Kuramoto model \eqref{kurODE}.
We summarize the overlap and attempt a brief comparison.

\textbf{P2} works in angles $\vth_j\ \mathrm{mod}\ 2\pi$, i.e.\ in $\bT^N\!=\!\bR^N/2\pi\bZ^N$, exclusively.
Sections \hbox{1--12} of \textbf{P2} largely run parallel to our present paper -- with an important difference in choice of coordinates, and, in part, stronger results.
Our flow-invariant cross section $\boXi\!=\!\mathbf{1}_N^\bot$ to uniform phase shifts $\bR\mathbf{1}_N$ is called ``counterdiagonal coordinates'', in \textbf{P2}.
Instead, \textbf{P2} prefers ``cube face coordinates'', described as follows.
Fixing the last $z$ components of $\bolth$ at $\vth_j\!=\!0$, let $\bT^d$ denote sub-tori on the complementary first $d\!=\!N-z$ components.
For $z\!=\!1,\, d\!=\!N\!-\!1$, this replaces $\boXi$ by a cross sectional torus $\bT^{N\!-\!1}$ which is \emph{not} flow-invariant under the Kuramoto flow \eqref{kurODE}.
Proposition 2.2 of \textbf{P2} therefore replaces \eqref{kurODE} by an ODE flow on the cross-section $\bT^{N\!-\!1}$ which represents the Kuramoto flow $\boPhi^t$ by a flow on \emph{group orbits} of $\bR\mathbf{1}_N$.
See our definition of equivariance in section \ref{Group}.
For general $d$, that orbit flow on $\bT^d$ is explicitly given by
\begin{equation}
\label{kurODEd}
\dot\vth_j=\tfrac{1}{N}\sum_{k=1}^d\sin(\vth_k-\vth_j) - \alpha\,\sin\vth_j -\tfrac{1}{N}\sum_{k=1}^d \sin\vth_k\,.
\end{equation}
Here $j\!=\!1,\ldots,d$ and $\alpha:=z/N\!=\!1\!-\!d/N$\,.
The first sum copies \eqref{kurODE} for $k\!\leq\! d$.
For $j\!\leq\!d\!<\!k$, the terms $\vth_k\!=\!0$ in $\bT^d$ contribute $-\alpha\,\sin\vth_j$ to $\dot{\vth_j}$\,.
The last sum projects \eqref{kurODE} back to the cross section $\bT^d$, along $\mathbf{1}_N$\,, to determine the induced flow of group orbits under $\bR\mathbf{1}_N$\,.
This ensures invariance $\dot\vth_j\!=\!0$ of $\vth_j\!=\!0$, for $j\!>\!d$, which would fail for the original Kuramoto ODE \eqref{kurODE}.
See section 10 of \textbf{P2}.
In our setting, $\bT^d$ is represented by $\bT_\mathbf{J}^d=\boXi_\mathbf{J}/\Lambda_\mathbf{J}'$\,, for the partition
\begin{equation}
\label{TJd}
\mathbf{J}\ =\ \{1\}\dot\cup\ldots\dot\cup\{d\}\dot\cup\{d+1,\ldots,N\}
\end{equation}
into $M=d+1$ clusters; see \eqref{TJM-1}.

We explain the advantages of ``cube face coordinates'', next.
In the following, assume
\begin{equation}
\label{dN2}
0<d<N/2\,,\quad\mathrm{i.e.}\quad 1>\alpha>1/2\,.
\end{equation}
The 2-cluster equilibrium $\boTh\!\in\!\bT^d$ with fat cluster size $N^b\!=\!z,\,\alpha\!=\!\alpha^b$ and slim cluster size $N^c\!=\!d\!=\!i(\boTh)$, for example, satisfies $\Theta_j\!=\!\pi$, for $1\!\leq\! j\!\leq\! d$, and $\Theta_j\!=\!0$ for $j\!>\!d$.
At such equilibria $\boTh$, and under assumption \eqref{dN2}, \textbf{P2} then shows
\begin{equation}
\label{Wud}
\mathrm{clos}\ W^u(\boTh)=\bT^d\,.
\end{equation}
The reason is simple:  $\boTh$ is the only equilibrium in $\bT^d$ with maximal unstable dimension $i(\boTh)\!=\!d$.
The boundary of $W^u(\boTh)$ in $\bT^d$ therefore consists of the unstable manifolds of the finitely many \emph{other} equilibria in $\bT^d$, all of which possess Morse index $i<d$.
In particular,  $W^u(\boTh)$ is open and dense in $\bT^d$.
In our setting $\mathbf{P1}$, uniqueness of $\boTh$ was observed in proposition \ref{Equil}(ii).
Indeed, its slim cluster $J^c=\{1,\ldots,d\}$ is maximal and unique, in $\bT_\mathbf{J}^d$\,, for the partition $\mathbf{J}$ in \eqref{TJd}.

Remarkably, section 10 in \textbf{P2} establishes $C^0$ orbit equivalence of \eqref{kurODEd} to the decoupled flow
\begin{equation}
\label{smaODE}
\dot{\vth}_j=-\tfrac{1}{N}\sin\vth_j\,,
\end{equation}
on $\bT^d$, i.e.\ on the first components $j\!=\!1,\ldots,d$.
What a beautiful template for the dynamics on any closure $\bT^d$ of the unstable manifolds \eqref{Wud}\,!
Most elegantly, this is achieved by a standard homotopy on $\bT^d$, from \eqref{kurODEd} to \eqref{smaODE}.
The flows remain gradient with respect to the associated homotopy of potentials.
Equilibria, their hyperbolicity and Morse indices, and transversality of their stable and unstable manifolds within $\bT^d$ all persist, throughout the homotopy.
Compare our transversality results and, in particular, \eqref{Wu}.

A main disadvantage of the ``cube face coordinates'' \textbf{P2}, however, is their lack of equivariance under the full permutation group $S_N$ in \eqref{groupS}, which the ``counterdiagonal'' cross section $\boXi\!=\!\mathbf{1}_N^\bot$ of \textbf{P1} carefully preserves.
Although the correspondence between \eqref{kurODEd} and \eqref{smaODE} is temptingly intuitive, it is therefore not completely clear how clustering notions, like 3-cluster heteroclinicity or $M$-cluster swarming, would fare under the abstract $C^0$ orbit equivalence of \textbf{P2}.
The dimension restriction $d\!<\!N/2$ of \eqref{dN2}, moreover, restricts transversality considerations to reside within the invariant tori $\bT^d$, and their variants under coordinate permutations.
Simultaneous structural stability of the Morse--Smale dynamics on the full attractor $\cA $ in $R\!>\!0$, cf.\ \eqref{cA}, has not been proved in \textbf{P2}.
In fact, structural stability remains restricted to tori $\bT^d,\ d\!<\!N/2$, and the \emph{a priori} fixed presentation as a CW cell complex, separately for each.
This limits the class of admissible (small) perturbations, geometrically, and is in fact less ``complete'', than suggested by the concluding section 19 of \textbf{P2}.
Compare our broader conclusions in section \ref{Pert}.

The remaining sections 13--18 of \textbf{P2} go beyond \textbf{P1}. 
They explore the very interesting problem of attaching the equilibrium variety $R\!=\!0$ of codimension 2, a manifold for odd $N$, to the Morse--Smale flow in $0\!<\!R\!\leq\!1$, via its two-dimensional strong unstable manifolds.
The cell decomposition, in sections 15,16, of the equilibrium variety, alias the $N$-bar linkage chains of unit joint length \cite{Nbar}, is not an inherently dynamical problem.
In section 17, that topological exploration culminates in the computations of homology dimensions.
The discussion of normal hyperbolicity of $R\!=\!0$, in sections 13 and 14, might benefit from Watanabe--Strogatz variables; see \cite{PiRoBunch, PiRoSurv} again.
Global dynamical consequences in section 18 are mainly substantiated and illustrated for the example $N\!=\!5$, under the notion of ``imprints''.
We have already mentioned in section \ref{Equi} how that part is beyond reach of the Watanabe--Strogatz trick.

Our present description \textbf{P1}, in contrast, has emphasized detailed descriptions of heteroclinic orbits $\bolth(t){:}\,\ \boTh_-\!\leadsto\!\boTh_+$ in terms of clustering, flow invariance of spaces fixed under group actions (the ``equality principle'' of section 3 in \textbf{P2}), rebellions, and their shift-encoded, $\eps$-shadowing concatenations.
We have also tried to highlight the order parameter $R$ of \eqref{R=}, throughout \textbf{P1}, to remain closer to the very extensive physics literature on the subject.
As a technical device, the absolute value $R$ of the ``centroid'' $Z$ appears in section 4 of \textbf{P2}.
The explicit relation to the Lyapunov potential $V$ in \textbf{P2} is $V=\tfrac{1}{2}(1\!-\!R^2)$.
Although our transversality and perturbation results are stronger, we admire the beauty and technical simplicity of the torus template \eqref{Wud}.

Numerical simulations in \textbf{P2} are planar.
Our figures \ref{fig51}(a) and \ref{fig52} correspond to \textbf{P2}, figure 7, including the choices of $N\!=\!5$ and the cluster ratio $\alpha_1\!=\!3/5$.
Each fundamental triangle in \textbf{P2}, figure 6, corresponds to our figure \ref{fig51}(b).
Concatenation towards synchrony of $N\!=\!21$ oscillators, like figures \ref{fig22}, \ref{fig23}, and the swarming saddle-saddle heteroclinic orbits like figure \ref{fig25}, all representable in $\bT^{10}$, have not been illustrated in \textbf{P2}.
The saddle-saddle connections, we confess, had not been quite as hard to find as suggested in section 11 of \textbf{P2}; see our section \ref{NumSwarm}.

In summary, we are delighted to suddenly have both viewpoints available, \textbf{P1} \emph{and} \textbf{P2}, some 51 years after Kuramoto's inspiring original note \cite{KurOsc}.
Even more remarkably than such a coincidence, this demonstrates how even the most elementary case of identical, unforced, all-to-all sinusoidal coupling still deserves our intense mathematical interest, now and towards future mathematical exploration.
May the present results foreshadow, and motivate, further development of a future, and more embracing, \emph{equivariant} Morse--Smale theory.



\newpage


\begin{thebibliography}{9999)999}

{\footnotesize{
\bibitem[ABPRS05]{ABPRS}
J.A.\ Acebrón, L.L.~Bonilla, C.J.~Pérez Vicente, F.~Ritort, and R.~Spigler.
The Kuramoto model: A simple paradigm for synchronization phenomena.
\emph{Rev.\ Mod.\ Phys.} \textbf{77} (2005), 137--185;
\url{https://doi.org/10.1103/RevModPhys.77.137}

\bibitem[ADLY25]{Yorke-b}
Ch.\ Adwani, R.\ De Leo, and J.A.\ Yorke.
What is the graph of a dynamical system? 
\emph{Nonlin.\ Dyn.}\ (2025);
\url{https://doi.org/10.1007/s11071-025-11466-9}

\bibitem[Arn88]{ArnoldGeoODE}
V.I.~Arnold.
\emph{Geometrical Methods in the Theory of Ordinary Differential Equations.}
Springer-Verlag, Berlin 1988.

\bibitem[Arn92]{ArnoldODE}
V.I.\ Arnold.
\emph{Ordinary Differential Equations.} 
Springer-Verlag, Berlin 1992. 

\bibitem[AOWT07]{Ash07}
P.~Ashwin, G.~Orosz, J.~Wordsworth, and S.~Townley.
Dynamics on Networks of Cluster States for Globally Coupled Phase Oscillators.
\emph{SIAM J.\ Appl.\ Dyn.\ Syst.} \textbf{6} (2007), 728--758;
\url{https://doi.org/10.1137/070683969}

\bibitem[Beyn90]{Beyn}
W.-J.\ Beyn.
The numerical computation of connecting orbits in dynamical systems. 
\emph{IMA J.\ Numer.\ Analysis} \textbf{10} (1990), 379--405. 

\bibitem[BEKS17]{Julia}
J.~Bezanson, A.~Edelman, S.~Karpinski, and V.B.~Shah. 
Julia: a fresh approach to numerical computing. 
\emph{SIAM Review} \textbf{59} (2017), 65--98.

\bibitem[BroStr18]{Kuramotocycle}
D.~Brockmann and S.~Strogatz.
\emph{Ride my Kuramotocycle.}
Web app (2018).
\url{https://www.complexity-explorables.org/explorables/ride-my-kuramotocycle/}

\bibitem[BrFie88]{brfi88}
P.~Brunovsk\'y and B.~Fiedler.
 Connecting orbits in scalar reaction diffusion equations.
 \emph{Dynamics Reported} \textbf{1} (1988), 57--89.

\bibitem[BrFie89]{brfi89}
P.~Brunovsk\'y and B.~Fiedler.
Connecting orbits in scalar reaction diffusion equations {II}: The complete solution.
 \emph{J.~Differ.~Eqs.} \textbf{81} (1989), 106--135.

\bibitem[ChHsHs24]{ChHsHs}
S.-H.~Chen, C.-H.~Hsia, and T.-Y.~Hsiao.
Complete and partial synchronization of two-group and three-group Kuramoto oscillators.
\emph{SIAM J.\ Appl.\ Dyn.\ Syst.} \textbf{23} (2024), no.~3, 1586--1613;
\url{https://doi.org/10.1137/23M1586227}

\bibitem[Chib15]{Chib}
H.~Chiba.
A proof of the Kuramoto conjecture for a bifurcation structure of the infinite-dimensional Kuramoto model.
\emph{Ergod.\,Th.\,Dynam.\,Sys.}\ \textbf{35} (2015), 762--834;
\url{https://doi.org/10.1017/etds.2013.68}

\bibitem[Chic03]{Ch}
C.~Chicone.
Inertial and slow manifolds for delay equations with small delays.
\emph{J.\ Differ.\ Eqs.} \textbf{190} (2003), 364--406;
\url{https://doi.org/10.1016/S0022-0396(02)00148-1}

\bibitem[ChoLau00]{ChoLau}
P.~Chossat and R.~Lauterbach.
\emph{Methods in Equivariant Bifurcations and Dynamical Systems.} 
World Scientific, Singapore 2000. 

\bibitem[ChoHa82]{ChowHale}
S.-N.~Chow and J.K.~Hale.
\emph{Methods of Bifurcation Theory.}
Springer-Verlag, New York 1982. 

 \bibitem[Con78]{Conley}
 C.C.~Conley.
 \emph{Isolated Invariant Sets and the Morse Index.} 
 AMS, Providence R.I.\ 1978. 
 
 \bibitem[DeLeo26]{deleo}
 R.\ De Leo.
 The chain graph of a general iterated function system.
 \url{https://arxiv.org/abs/2607.07112}

 \bibitem[DLY24]{Yorke-a}
 R.\ De Leo and J.A.\ Yorke.
 Streams and graphs of dynamical systems.
 \emph{Qual.\ Theory Dyn.\ Syst.} \textbf{24} (2024);
 \url{https://doi.org/10.1007/s12346-024-01112-x}

\bibitem[DGK03]{KuzCode}
A.~Dhooge, W.~Govaerts, and Yu.A.~Kuznetsov.
MATCONT: A MATLAB package for numerical bifurcation analysis of ODEs. 
\emph{ACM Trans.\ Math.\ Softw.} \textbf{29} (2003), 141--164. 

\bibitem[DiSt03]{DiSt}
A.P.S.~Dias and I.~Stewart.
Secondary bifurcations in systems with all-to-all coupling.
\emph{Proc.\ R.~Soc.\ Lond.} \textbf{A 459} (2003), 1969--1986, 
\url{https://doi.org/10.1098/rspa.2002.1103}

\bibitem[DoXu13]{DoXu}
J.-G.~Dong and X.~Xue.
Synchronization analysis of Kuramoto oscillators.
\emph{Comm.\ Math.\ Sci.} \textbf{11} (2013), no.~2, 465--480;
\url{https://dx.doi.org/10.4310/CMS.2013.v11.n2.a7}

\bibitem[DöBu12]{DoeBu}
F.~Dörfler and F.~Bullo.
Synchronization and transient stability in power networks and nonuniform Kuramoto oscillators.
\emph{SIAM J.\ Control Optim.} \textbf{50} (2012), no.~3, 1616--1642;
\url{https://doi.org/10.1109/ACC.2010.5530690}

\bibitem[Elm01]{Elm1}
T.~Elmhirst.
\emph{Symmetry and Emergence in Polymorphism and Sympatric Speciation.}
PhD Thesis, Warwick 2001.

\bibitem[Elm04]{Elm2}
T.~Elmhirst.
$S_N$-equivariant symmetry-breaking bifurcations. 
\emph{Int.\ J.\ Bifurcation Chaos} \textbf{14} (2004), 1017--1036. 

\newpage
\bibitem[Fie88]{FieHopf} 
B.~Fiedler. 
\emph{Global Bifurcation of Periodic Solutions with Symmetry.} 
Lect.\ Notes Math.\ \textbf{1309}, Springer-Verlag, Heidelberg 1988.

\bibitem[Fie94]{FieTatra} 
B.~Fiedler.
Global attractors of one-dimensional parabolic equations: sixteen examples.
\emph{Tatra Mountains Math.\ Publ.} \textbf{4} (1994), 67--92.

\bibitem[Fie02]{fi02}
B.~Fiedler (ed.).  
\emph{Handbook of Dynamical Systems} \textbf{2}. 
Elsevier, Amsterdam 2002.

\bibitem[FKHK21]{Krischer2}
B.~Fiedler, F.P.~Kemeth, S.W.~Haugland, K.~Krischer.
Global heteroclinic rebel dynamics among large 2-clusters in permutation equivariant systems.
\emph{SIAM Appl.\ Dyn.\ Systems} \textbf{20} (2021), 1277--1319;
\url{https://doi.org/10.1137/20M1361493}

\bibitem[FieRo96]{cascading}
B.~Fiedler and C.~Rocha.
Heteroclinic orbits of semilinear parabolic equations.
\emph{J.~Differ.~Eqs.} \textbf{125} (1996), 239--281. 

\bibitem[FieRo20]{ThomSmale} 
B.~Fiedler and C.~Rocha.
Boundary orders and geometry of the signed Thom--Smale complex for Sturm global attractors.
\emph{J.~Dyn.~Differ.~Eqs.} (2020), 2787--2818;
\url{https://doi.org/10.1007/s10884-020-09836-5}

\bibitem[FieRo23]{firoSFB}
B.~Fiedler and C.~Rocha.
Design of Sturm global attractors 1: Meanders with three noses, and reversibility.
\emph{Chaos} \textbf{33}, 083127 (2023); \url{https://doi.org/10.1063/5.0147634}

\bibitem[FieRo24]{firoFusco}
B.~Fiedler and C.~Rocha.
Design of Sturm global attractors 2: Time-reversible Chafee--Infante lattices of 3-nose meanders.
\emph{São Paulo J.\ Math.\ Sciences.} (2024), 975--1014;
\url{https://doi.org/10.1007/s40863-023-00385-5}

\bibitem[FrM88]{franzosa}
R.D.~Franzosa and K.~Mischaikow.
The connection matrix theory for semiflows on (not necessarily locally compact) metric spaces. 
\emph{J.\ Differ.\ Eqs.} \textbf{71} (1988), 270--287. 

\bibitem[GMK08]{VGM2008}
V.~Garc\'{\i}a-Morales and K.~Krischer.
Normal-form approach to spatiotemporal pattern formation in globally coupled electrochemical systems.
\emph{Phys.\ Rev.} \textbf{E 78} (2008), 057201.

\bibitem[GoStSch88]{GoStSymm}
M.~Golubitsky, I.~Stewart, and D.G.~Schaeffer.
\emph{Singularities and groups in bifurcation theory.\ II.}
Springer-Verlag, New York 1988. 

\bibitem[GoSt02]{GoStPersp}
M.~Golubitsky and I.~Stewart.
\emph{The Symmetry Perspective.} 
Birkh\"auser, Basel 2002.

\bibitem[GoSt23]{GoStNet}
M.~Golubitsky and I.~Stewart.
\emph{Dynamics and Bifurcation in Networks. Theory and Applications of Coupled Differential Equations.} 
SIAM, Philadelphia 2023.

\bibitem[GPS02]{Sauer}
C.~Grebogi, L.~Poon, T.~Sauer et al.
Shadowability of chaotic dynamical systems.
In \cite{fi02}, 313--344. 

\bibitem[HKP15]{HKP}
S.-Y.~Ha, H.~Kim, and J.~Park.
Remarks on the complete frequency synchronization of Kuramoto oscillators.
\emph{Nonlinearity} \textbf{28} (2015), 1441--1462;
\url{https://doi.org/10.1088/0951-7715/28/5/1441}

\bibitem[HKL14]{HKL}
S.-Y.~Ha, Y.~Kim, and Z.~Li.
Large-time dynamics of Kuramoto oscillators under the effects of inertia and frustration.
\emph{SIAM J.\ Appl.\ Dyn.\ Syst.} \textbf{13} (2014), no.~1, 466--492;
\url{https://doi.org/10.1137/130926559}

\bibitem[Hale69]{Hale}
J.K.\ Hale.
\emph{Ordinary Differential Equations.}
Wiley-Interscience, New York 1969. 

\bibitem[HMO02]{HMO}
J.K.~Hale, L.T.~Magalhães, and W.M.~Oliva.
\emph{Dynamics in Infinite Dimensions.} 2nd ed.
Springer-Verlag, New York 2002. 

\bibitem[HaVL93]{halu93}
J.K.~Hale and S.M.~Verduyn Lunel.
\emph{Introduction to Functional Differential Equations}.
Springer-Verlag, New York 1993.

\bibitem[vHeWre93]{vanHemmen}
J.L.~van Hemmen and W.F.~Wreszinski.
Lyapunov function for the Kuramoto model of nonlinearly coupled oscillators.
\emph{J.~Stat.~Phys.} \textbf{27} (1993), 145--166.

\bibitem[Hen85]{henrytrv}
D.B.\ Henry.
Some infinite-dimensional Morse-Smale systems defined by parabolic partial differential equations.
\emph{J.\ Differ.\ Eqs.} \textbf{59} (1985), 165--205. 

\bibitem[HiPuSh77]{HPS}
M.W.\ Hirsch, C.C.\ Pugh, and M.\ Shub.
\emph{Invariant Manifolds.}
Lect.\ Notes Math.\ \textbf{583}.
Springer-Verlag, New York 1977. 

\bibitem[HiSmDe13]{HSD}
M.W.\ Hirsch, S.\ Smale, and R.L.\ Devaney.
\emph{Differential Equations, Dynamical Systems, and an Introduction to Chaos.} 3rd ed.
Academic Press, Amsterdam 2013. 

\bibitem[HsJuKw19]{HsJuKw}
C.-H.~Hsia, C.-Y.~Jung, and B.~Kwon.
On the synchronization theory of Kuramoto oscillators under the effect of inertia.
\emph{J.\ Differ.\ Eqs.} \textbf{267} (2019), 742--775;
\url{https://doi.org/10.1016/j.jde.2019.01.024}

\bibitem[HJKU20]{HJKU}
C.-H.~Hsia, C.-Y.~Jung, B.~Kwon, and Y.~Ueda.
Synchronization of Kuramoto oscillators with time-delayed interactions and phase lag effect.
\emph{J.\ Differ.\ Eqs.} \textbf{268} (2020), 7897--7939; \url{https://doi.org/10.1016/j.jde.2019.11.090}

\newpage
\bibitem[Jo95]{Jones}
C.~Jones.
\textit{Geometric singular perturbation theory.}
In R.~Johnson (Ed.), \emph{Dynamical Systems}, Lecture Notes in Mathematics, \textbf{1609}, 44--118.
Springer, Berlin 1995.

 \bibitem[KaMM04]{Mischaikow}
 T.~Kaczynski, K.~Mischaikow, and M.~Mrozek.
 \emph{Computational Homology.} 
 Springer-Verlag, New York 2004. 

\bibitem[KFHK20]{Krischer1}
F.P.~Kemeth, B.~Fiedler, S.W.~Haugland, and K.~Krischer.
2-cluster fixed-point analysis of mean-coupled Stuart--Landau oscillators in the center manifold of the Benjamin--Feir instability. 
\emph{J.~Phys.\ Complex.} (2020); 
\url{https://doi.org/10.1088/2632-072X/abd0da}

\bibitem[KOD05]{DoedelSurvey}
B.~Krauskopf, H.M.~Osinga, E.J.~Doedel, et al.
A survey of methods for computing (un)stable manifolds of vector fields. 
\emph{Int.\ J.\ Bif.\ Chaos Appl.\ Sci.\ Eng.} \textbf{15} (2005), 763--791.

\bibitem[KGO15]{KuGiOtt2015}
W.L.~Ku, M.~Girvan, and E.~Ott.
Dynamical transitions in large systems of mean field-coupled Landau--Stuart oscillators: Extensive chaos and cluster states.
\emph{Chaos} \textbf{25} 123122 (2015); 
\url{https://doi.org/10.1063/1.4938534}

\bibitem[Kur75]{KurOsc}
Y.~Kuramoto. 
Self-entrainment of a population of coupled non-linear oscillators.
\emph{Int.\ Symp.\ Math.\ Probl.\ Theor.\ Phys.}, 
H.\ Arakai (ed), Kyoto 1975, 420--422.

\bibitem[Kur84]{Kurbook}
Y.~Kuramoto. 
\emph{Chemical Oscillations, Waves, and Turbulence.}
Springer-Verlag, Berlin, 1984.

\bibitem[KurB02]{KBatt}
Y.\ Kuramoto and D.\ Battogtokh. 
Coexistence of coherence and incoherence in nonlocally coupled phase oscillators.
\emph{Nonlin.\ Phenom.\ Complex Syst.} \textbf{5} (2002), 380--385.

\bibitem[Kurz67]{Kurz}
J.~Kurzweil.
Invariant manifolds for flows.
\emph{Differ.\ Equations Their Appl.} (1967),  Bratislava: Slovenské pedagogické nakladateľstvo, 89--92.

\bibitem[Lin90]{Lin}
X.-B.~Lin.
Using Melnikov's method to solve Silnikov's problems.
\emph{Proc.\ R.\ Soc.\ Edinb.\ A: Math.} \textbf{116} (1990), 295--325;
\url{https://doi.org/10.1017/S0308210500031528}

 \bibitem[Lin08]{Linschol}
 X.-B.~Lin.
 Lin's method.
 \emph{Scholarpedia} (2008);
\url{https://doi.org/10.4249/scholarpedia.6972}

\bibitem[McCSoh11]{Nbar}
J.M.~McCarthy and G.S.~Soh.
\emph{Geometric Design of Linkages.}
Springer-Verlag, New York 2011.

\bibitem[MM02]{Mischaikow02}
K.~Mischaikow and M.~Mrozek.
Conley index.
In \cite{fi02}, 393--460. 

\bibitem[Pali69]{Palis}
J.~Palis.
 On Morse--Smale dynamical systems.
 \emph{Topology} \textbf{8} (1969), 385--404.

\bibitem[PdM82]{PalisdeMelo}
J.~Palis and W.~de~Melo.
\emph{Geometric Theory of Dynamical Systems. An Introduction.}
Springer-Verlag, New York 1983. 

\bibitem[PaSm70]{PalisSmale}
J.~Palis and S.~Smale.
 Structural stability theorems.
 In \emph{Global Analysis},  S.\ Chern, S.\ Smale (eds.). Proc.\ Symp.\ in
 Pure Math.~vol.~XIV.~AMS, Providence 1970. 
 
 \bibitem[Palm09]{Palmer}
 K.J.\ Palmer.
 Shadowing lemma for flows.
 Scholarpedia (2009), 7918;
 \url{https://doi.org/10.4249/scholarpedia.7918}

\bibitem[PPT12]{Tikhomirov}
K.J.~Palmer, S.Yu.~Pilyugin, and S.B.~Tikhomirov.
Lipschitz shadowing and structural stability of flows.
\emph{J.\ Differ.\ Eqs.} \textbf{252} (2012), 1723--1747. 

\bibitem[PiRo09]{PiRoBunch}
A.~Pikovsky and M.~Rosenblum.
Self-organized partially synchronous dynamics in populations of nonlinearly
coupled oscillators.
\emph{Physica D} \textbf{238} (2009), 27--37;
\url{https://doi.org/10.1016/j.physd.2008.08.018}

\bibitem[PiRo15]{PiRoSurv}
A.~Pikovsky and M.~Rosenblum.
Dynamics of globally coupled oscillators: Progress and perspectives.
\emph{Chaos} \textbf{25} (2015), 097616;
\url{https://doi.org/10.1063/1.4922971}

\bibitem[PRK03]{Pik2003}
A.~Pikovsky, M.~Rosenblum, and J.~Kurths.
\emph{Synchronization: A Universal Concept in Nonlinear Sciences.}
Cambridge University Press, 2003.

\bibitem[Rau02]{raugel}
G.~Raugel.
Global attractors in partial differential equations. 
In \cite{fi02}, 885--982. 

\bibitem[RFLN25]{rofiln25}
C.~Rocha, B.~Fiedler, and A.~L\'{o}pez-Nieto.
A classification of global attractors for $\mathbb{S}^1$-equivariant parabolic equations.
(2025);
\url{https://arxiv.org/abs/2507.10051}

\bibitem[SePaLe07]{IEEE}
R.\ Sepulchre, D.A.\ Paley, and N.E.\ Leonhard.
Stabilization of planar collective motion: All-to-all communication.
\emph{IEEE Trans.\ Aut.\ Control} \textbf{52} (2007), 1--14.

\bibitem[StElCo03]{StElm}
I.~Stewart, T.~Elmhirst, and J.~Cohen.
Symmetry-breaking as an origin of species.
In \emph{Bifurcation, Symmetry and Patterns}, J.~Buescu et al (eds.).
Birkhäuser, Basel 2003, 3--54. 

\bibitem[Str03]{Str}
S.H.~Strogatz.
\emph{Sync: The Emerging Science of Spontaneous Order}. Hyperion, New York 2003.

\bibitem[Van82]{Vander}
A.~Vanderbauwhede.
\emph{Local Bifurcation and Symmetry.}
Pitman, Boston 1982.
}}

\end{thebibliography}
\end{document}